\documentclass[11pt]{article}
\usepackage[a4paper,margin=1in]{geometry}
\usepackage{amsmath,amssymb,amsthm,mathtools}
\usepackage[textsize=footnotesize, colorinlistoftodos, obeyFinal]{todonotes}

\usepackage{booktabs,tabularx,array}
\usepackage{enumitem}
\usepackage[normalem]{ulem}
\usepackage{float}
\theoremstyle{plain}
\newtheorem{thm}{Theorem}[section]

\newtheorem{cor}[thm]{Corollary}
\newtheorem{ex}[thm]{Example}
\newtheorem{Assumption}[thm]{Assumption}
\newtheorem{lem}[thm]{Lemma}

\newtheorem{prop}[thm]{Proposition}

\theoremstyle{definition}

\theoremstyle{remark}
\newtheorem{rem}[thm]{Remark}
\providecommand{\C}{}
\renewcommand{\C}{\mathbb{C}}
\newcommand{\diag}{\mathrm{diag}}

\newcommand{\lawto}{\overset{\mathrm{law}}{\to}}
\usepackage{algorithm}
\usepackage{algpseudocode}
\usepackage[unicode,hidelinks,bookmarks=false]{hyperref}
\setlist[itemize]{itemsep=3pt,topsep=4pt,leftmargin=2.2em}
\setlist[enumerate]{itemsep=4pt,topsep=4pt,leftmargin=2.4em}
\allowdisplaybreaks

\title{Fluctuations of sample eigenvalues in a spiked covariance model with Haar-rotated spikes: a transition at $p^2 \asymp n$}

\author{
Noriyoshi Sakuma\thanks{
Department of Mathematics, Graduate School of Science, The University of Osaka, Japan.
Email: \texttt{sakuma@math.sci.osaka-u.ac.jp}.
}
\and
Jean-Lou Valeau\thanks{
ENSAE Paris, France.
Email: \texttt{jean-lou\_valeau@orange.fr}.
}
}
\date{November 10, 2026}

\begin{document}
\maketitle

\begin{abstract}
We study the outlying eigenvalues of a sample covariance matrix whose population
covariance is itself random. The population covariance is
$\Sigma_p=A_p+O_pB_pO_p^*+\sigma^2I_p$, where $A_p$ and $B_p$ are deterministic positive
semidefinite matrices of fixed ranks $r$ and $s$, and where $O_p$ is Haar distributed on
the orthogonal group in the real case, and on the unitary group in the complex case.
Conditionally on $O_p$, we observe $n$ independent centered Gaussian vectors with
covariance $\Sigma_p$, and we build the sample covariance matrix $\widehat\Sigma_p$ from these vectors. The
fluctuations of the spike eigenvalues of $\widehat\Sigma_p$ come from two different sources: the Haar matrix $O_p$, which displaces the population spikes at the scale $1/p$, and
the sampling, which acts at the scale $n^{-1/2}$. Comparing these two scales amounts to
comparing $p^2$ with $n$, so the ratio which governs the fluctuations is $p^2/n$, while
in the classical spiked model it is the aspect ratio $p/n$. We give the limiting law of
each spike in four regimes where $n,p\to \infty$. When $p/n\to c>0$, and when $p/n\to0$ with $p^2/n\to\infty$,
the sampling term dominates, the limit does not depend on $O_p$ and is Gaussian, as in the classical spiked model. When
$p^2/n\to0$, the Haar randomness dominates, and the limit is a weighted sum of independent
$\chi^2_1$ variables in the real case, and of independent exponential variables in the
complex case. When $p^2/n\to\tau\in(0,\infty)$, the two terms are of the same order, and
the limit is the convolution of Gaussian and $\chi^2_1$ or Exponential variables. We then apply these results
to a linear mixed model with a continuous batch effect, in which the Haar assumption is
what makes the separation of the signal and of the batch possible, and we build a test for
the value of a single spike.
\end{abstract}
\medskip
{\small
\noindent\textbf{2020 Mathematics Subject Classification.}
Primary 60B20; Secondary 60F05, 62H10.
\par\smallskip
\noindent\textbf{Keywords.}
Spiked covariance model; Haar random matrices;
eigenvalue fluctuations; central limit theorem;
low-aspect-ratio regime; spectral inference.
\par}
\medskip
\tableofcontents


\section{Introduction}
\label{s:intro}

Classical spiked covariance estimation theory treats the population covariance matrix as
deterministic. In many situations, however, a part of the covariance structure is itself
produced by a random mechanism: a latent environment, a nuisance effect, or a random
effect drawn once and then shared by the observations. The population covariance then
carries deterministic spikes and randomly oriented spikes at the same time. In this paper, we study sample covariance matrices whose population covariance has randomly rotated spikes. 
In this situation, we have to estimate both population randomness and sampling randomness. To handle such situations, we develop a fluctuation theory for such models.

In our model the classical aspect ratio $p/n$ alone no longer determines the limiting behavior. 
In the low-aspect-ratio regime $p/n\to 0$, the finer ratio $p^{2}/n$, which does not appear in the classical spiked model, turns out to be the relevant parameter.

The reason is that a finite rank spike leaves the limiting spectral distribution unchanged, so that first-order free probability cannot detect it.
The appropriate framework is cyclic-monotone independence, introduced in \cite{chs18}, which is designed precisely to describe the discrete part of the spectrum of the model with a Haar-rotated finite rank perturbation. 
The fluctuations of such discrete eigenvalues were subsequently computed in \cite{collins24} for polynomials in Haar unitary and finite rank matrices.
The present paper extends this line of results to the sample covariance setting and reveals a fluctuation transition governed by $p^{2}/n$.
Depending on their relative size, the limiting law is Gaussian, non-Gaussian, or the convolution of the two.

Building on this framework, we compute the fluctuations of the outlying eigenvalues in each of the above regimes.
A key feature of our setting is that the fluctuation comes not only from the sampling variability of the entries of the data matrix, but also from the randomness of the rotation carried by the spikes. 
These two sources are asymptotically independent, and their relative weight is governed by
$p^2/n$: they are of the same order only when $p^2/n\to\tau\in(0,\infty)$. This transition
cannot be seen in \cite{collins24}, where there is no sampling step, nor in the classical spiked model, where the population matrix is deterministic.
Quantifying this transition is the main contribution of this paper.

\subsection{The classical spiked covariance model}
\label{ss:spiked}

In this work we study a spiked covariance model. The spiked covariance model is a small-rank
perturbation of a multiple of the identity. It was introduced by
Johnstone \cite{johnstone01} and has become one of the standard models of
high-dimensional statistics; the books \cite{bai10,yao15} give a complete account.

We recall foundational results of random matrix theory, which are the starting point of
our work. Let $\Sigma_p$ be a deterministic $p\times p$ covariance matrix whose
eigenvalues are $\sigma^2+\theta_1,\ldots,\sigma^2+\theta_K$ and $\sigma^2$ with
multiplicity $p-K$, where $K$ is fixed and $\theta_1>\cdots>\theta_K>0$ are called the
spikes. Let $\widehat\Sigma_p$ be the sample covariance matrix built from $n$
independent centered real or circular complex Gaussian vectors with covariance
$\Sigma_p$, and let $\widehat\lambda_{p,1}\ge\cdots\ge\widehat\lambda_{p,p}$ be its
eigenvalues, which we call the sample eigenvalues. Assume the proportional regime
$p/n\to c>0$. Four important results describe the behavior of these eigenvalues
in the proportional regime.

First, most of the eigenvalues do not stay near $\sigma^2$. The empirical distribution
of the eigenvalues of $\widehat\Sigma_p$ converges to the Mar\v{c}enko--Pastur
law \cite{mp67}. When $c\le1$ this law is absolutely continuous and spread over the
interval $[\sigma^2(1-\sqrt c)^2,\sigma^2(1+\sqrt c)^2]$, which is called the bulk.

Second, a spike is visible only if it is large enough. Baik and
Silverstein \cite{baik06} showed that there is a threshold at $\sigma^2\sqrt c$ such
that for $\ell\leq K$,
\[
  \widehat\lambda_{p,\ell}
  \xrightarrow{\ \mathrm{a.s.}\ } \begin{cases}
      \phi_c(\theta_\ell)
  :=(\sigma^2+\theta_\ell)\Bigl(1+\frac{c\,\sigma^2}{\theta_\ell}\Bigr)
  \ \text{ if $\theta_\ell>\sigma^2\sqrt c$,} \\
  \sigma^2(1+\sqrt c)^2 \ \text{ if $\theta_\ell\le\sigma^2\sqrt c$.}
  \end{cases}
\]
Notice that $\phi_c(\theta_\ell)>\sigma^2+\theta_\ell$ so
the sample eigenvalue is pushed up if $\theta_\ell>\sigma^2\sqrt c$. If not $\widehat\lambda_{p,\ell}$ converges to
the upper edge of the bulk and carries no information on $\theta_\ell$. The transition itself was discovered
earlier, for complex data, by Baik, Ben Arous, P\'ech\'e \cite{bbp05}: the threshold
is named after them and is called the BBP transition. A spike above the threshold is
called supercritical.

Third, the fluctuations of a simple supercritical spike are real Gaussian. Bai and Yao \cite{bai08,bai12} proved that
\[
  \sqrt n\bigl(\widehat\lambda_{p,\ell}-\phi_{p/n}(\theta_\ell)\bigr)
  \xrightarrow{\ \mathrm{law}\ }
  \mathcal N\bigl(0,\kappa\,V_c(\theta_\ell)\bigr),
  \qquad
  V_c(\theta):=(\sigma^2+\theta)^2\Bigl(1-\frac{c\,\sigma^4}{\theta^2}\Bigr),
\]
where $\kappa=2$ for real data and $\kappa=1$ for complex data. Fourth, the eigenvalues that stay in the bulk behave in a completely
different way. The largest of them converges to the upper bulk edge
$\sigma^2(1+\sqrt c)^2$, and it fluctuates at the scale $n^{-2/3}$; its limit law is
the Tracy--Widom law, of order one for real Gaussian data \cite{johnstone01,tw96} and
of order two for complex Gaussian data \cite{johansson00,tw94}.

These four facts hold in the proportional regime, when $p$ and $n$ are comparable. In
the disproportional regime, i.e. when $p/n\to0$, the behavior of eigenvalues is quite
different. For instance, in the simplest case where there is no spike and
$\Sigma_p=\sigma^2I_p$, the empirical distribution of the rescaled eigenvalues
converges to a semi-circle law, and not to the Mar\v{c}enko--Pastur law \cite{byin88}.
Moreover, when $\Sigma_p$ is a spiked covariance, as $p/n\to0$, all spikes become
supercritical and there is no BBP transition. In particular, instead,
$\forall\ell\leq K$,
$$\widehat\lambda_{p,\ell}\xrightarrow{\ \mathbb P\ }\sigma^2+\theta_\ell,$$
which comes from the standard fact
$\|\widehat\Sigma_p-\Sigma_p\|=\mathrm O_{\mathbb P}(\sqrt{p/n})$: this is for instance
a consequence of \cite[Corollary 2]{kol17}.

We study both  the cases $p/n\to c>0$ and $p/n\to0$, and
we will show that our model behaves as above in the proportional regime, but not in the disproportional one. The reason lies in the definition of the model, which we give now. 
Related central limit theorems for spiked sample eigenvalues without assuming convergence of the aspect ratio were obtained by Liu et al.~\cite{liu2022rds}. For our conditional argument, we establish the required triangular-array result in Proposition~\ref{p:dclt}.

\subsection{The model studied in this paper}
\label{ss:model}

Let $\mathbb F=\mathbb R$ or $\mathbb F=\mathbb C$, and set $\kappa=2$ if $\mathbb
F=\mathbb R$ and $\kappa=1$ if $\mathbb F=\mathbb C$. Let $r,s\ge1$ be fixed integers,
and let
\[
  \alpha_1,\ldots,\alpha_r,\qquad \beta_1,\ldots,\beta_s
\]
be strictly positive numbers, all distinct. For $p\ge r+s$ let $A_p$ and $B_p$ be
deterministic positive semidefinite matrices of rank $r$ and $s$ respectively, whose
nonzero eigenvalues are $\alpha_1,\ldots,\alpha_r$ for $A_p$ and
$\beta_1,\ldots,\beta_s$ for $B_p$, let $\sigma^2>0$ be fixed, and let $O_p$ be a
random matrix distributed according to the Haar measure on the orthogonal group $O(p)$
if $\mathbb F=\mathbb R$, and on the unitary group $U(p)$ if $\mathbb F=\mathbb C$ We define formally these groups in Section~\ref{s:grp}. The
population covariance matrix is
\begin{equation}
  \Sigma_p
  =
  A_p+O_pB_pO_p^* + \sigma^2 I_p,
  \qquad \sigma^2>0,
  \label{e:sigm}
\end{equation}
where $O_p^*$ denotes the conjugate transposed matrix of $O_p$ when $O_p$ is unitary,
and the transposed matrix when $O_p$ is orthogonal. The Haar measure is defined
formally in Section~\ref{s:grp}. Informally, $O_p$ is picked uniformly at random among
all orthogonal, respectively unitary, matrices of size $p$. After $O_p$ is drawn, we
observe i.i.d. centered Gaussian vectors with covariance $\Sigma_p$ conditionally on
$O_p$
\[
  x_{p,1},\ldots,x_{p,n} \ | \ O_p
  \stackrel{\mathrm{i.i.d.}}{\sim}\mathcal N_\mathbb F(0, \Sigma_p).
\]
We write $\mathcal N_\mathbb F$ for the Gaussian law when $\mathbb F=\mathbb R$, and for
the circular complex Gaussian law with covariance $\Sigma$ when $\mathbb F=\mathbb C$.
From these vectors, we form the sample covariance matrix
\[
  \widehat\Sigma_p
  =
  \frac1n\sum_{k=1}^n x_{p,k}x_{p,k}^*.
\]
We write $\lambda_{p,\ell}$ for a spike eigenvalue of the population covariance matrix
$\Sigma_p$, and $\widehat\lambda_{p,\ell}$ for the corresponding spike eigenvalue of
$\widehat\Sigma_p$. We consider these eigenvalues to be sorted, i.e.
$\lambda_{p,1}\geq\dots\geq\lambda_{p,p}$ and
$\widehat\lambda_{p,1}\geq\dots\geq\widehat\lambda_{p,p}$. We write
$\theta_1>\cdots>\theta_{r+s}>0$ for the sorted elements of
$\{\alpha_1,\ldots,\alpha_r,\beta_1,\ldots,\beta_s\}$, and we call them the spikes of
the model.

Three comments explain why this model is worth studying.

(i) The matrix $\Sigma_p$ is the sum of two perturbations of $\sigma^2I_p$, differing
by their eigenvectors. Let $a_1,\ldots,a_r$ be orthonormal eigenvectors of $A_p$
associated with $\alpha_1,\ldots,\alpha_r$, and let $b_1,\ldots,b_s$ be orthonormal
eigenvectors of $B_p$ associated with $\beta_1,\ldots,\beta_s$. The eigenvectors of
$O_pB_pO_p^*$ are the vectors $O_pb_1,\ldots,O_pb_s$, so they are random and uniformly
oriented. Indeed, when $O_p$ is Haar distributed, $(O_pb_1,\ldots,O_pb_s)$ is uniformly
distributed over all orthonormal families of $s$ vectors of $\mathbb F^p$. This is the
simplest way to write two small-rank perturbations that are uniformly oriented with
respect to each other. Hence, this model is relevant in
applications where the data carries an unwanted structure along directions which are
unknown, and which have no particular relation with the directions of interest; this
happens for instance when the unwanted structure comes from an unpredictable
measurement issue during the collection of the data. We discuss such applications in
Section~\ref{ss:lmm}.

(ii) Haar conjugation is a standard source of asymptotic freeness \cite{voiculescu91}. 
For sums of randomly rotated matrices, this framework describes limiting empirical spectral distributions through free additive convolution. In our fixed-rank setting, however, the empirical spectral distributions of $A_p$ and
$O_pB_pO_p^*$ both converge to $\delta_0$. 
The free additive convolution of these limiting measures is therefore also $\delta_0$, which does not determine the finitely many nonzero eigenvalues.

The relevant effect of Haar rotation here is that the two spike subspaces become asymptotically orthogonal. Consequently, the nonzero eigenvalues of $A_p+O_pB_pO_p^*$ converge to the combined set of strengths
$\{\alpha_1,\ldots,\alpha_r,\beta_1,\ldots,\beta_s\}$;
see Proposition~\ref{p:pop}.
Cyclic-monotone independence \cite{chs18} provides a framework for studying such discrete spectral limits.
For finite-rank Haar unitary models with simple limiting
eigenvalues, fluctuation formulas were obtained in \cite{collins24}.
Our work adds a Gaussian sampling step and studies how sampling
fluctuations combine with the fluctuations of the population spikes.

(iii) Finally, and this is the point that separates our model from the classical spiked
model, the population covariance matrix $\Sigma_p$ is itself random. In Johnstone's
model, $\Sigma_p$ is a fixed matrix and the only randomness comes from the $n$ samples.
Here there are two independent sources of randomness: the matrix $O_p$, drawn once, and
the $n$ sample vectors. The eigenvalue $\lambda_{p,\ell}$ of $\Sigma_p$ is therefore
itself a random variable, and the sample eigenvalue $\widehat\lambda_{p,\ell}$
fluctuates around a random target. This is what makes the problem different, and we
describe it in the next section.

\subsection{Two sources of randomness}
\label{ss:why}

By fluctuations, we refer to the asymptotic distribution of the rescaled difference
between the sample eigenvalue $\widehat\lambda_{p,\ell}$ and its deterministic limit.
In the classical spiked model, the exact scaling rate and limiting law were identified
by Bai and Yao \cite{bai08,bai12}: the normalization is $\sqrt n$ and the limit is Gaussian. In our model, while the classical Gaussian
fluctuations remain valid for $p/n\to c>0$, the regime $p/n\to0$ yields several distinct fluctuation regimes
depending on how fast $n$ grows relative to $p$. These unusual
fluctuations come from the two sources of randomness in the model.

\textbf{The Haar fluctuation.} 
Even before any sample is drawn, the
eigenvalue $\lambda_{p,\ell}$ of $\Sigma_p$ is random and does not equal
$\sigma^2+\theta_\ell$. The difference depends on the overlaps between the two spike
subspaces,
\[
  c_{ij}^{(p)}=a_i^*O_pb_j,
  \qquad 1\le i\le r,\quad 1\le j\le s .
\]
Because $O_p$ is Haar distributed, $O_pb_j$ is uniformly distributed on the unit sphere
of $\mathbb F^p$, so $\mathbb E[(O_pb_j)(O_pb_j)^*]=I_p/p$ by invariance. Hence
$c_{ij}^{(p)}$ has mean zero and second moment $1/p$, and it is of size $p^{-1/2}$.

The spike eigenvalues of $\Sigma_p$ depend on $O_p$ only through these $rs$ overlaps.
Indeed, $A_p+O_pB_pO_p^*$ has rank $r+s$, and its range is spanned by the $r+s$ vectors
$a_1,\ldots,a_r,O_pb_1,\ldots,O_pb_s$, which are almost surely linearly independent
because of $O_p$ being Haar. Taking these vectors as a basis of the range, the matrix of
the restriction has the $\alpha_i$ and the $\beta_j$ on its diagonal, and entries proportional to the
$c_{ij}^{(p)}$ outside of it; see Collins et al. \cite{collins24}. If all the overlaps
were zero, this matrix would be $\diag(\alpha_1,\ldots,\alpha_r,\beta_1,\ldots,\beta_s)$,
whose eigenvalues are the $\theta_\ell$, and they are simple because the $\theta_\ell$
are distinct. A simple eigenvalue is an analytic function of the entries of a matrix, so
each $\lambda_{p,\ell}$ is an analytic function of the overlaps in a neighbourhood of the
origin. This function has no linear term. Indeed, $b_j$ is only defined up to a
sign in the real case, and up to a phase in the complex case, and this change multiplies
the $j$-th column of $(c_{ij}^{(p)})_{ij}$ by $-1$, respectively by $\mathrm e^{\mathrm
i\varphi}$, while it changes neither $B_p$ nor $\Sigma_p$. Its constant term is its
value at zero overlap, that is, $\sigma^2+\theta_\ell$. So the leading term of
$\lambda_{p,\ell}-\sigma^2-\theta_\ell$ is quadratic in the overlaps, and it is of size
$1/p$ and not $p^{-1/2}$, as each overlap is of size $p^{-1/2}$.

Proposition~\ref{p:pop}, which is a specialization of the result of Collins et
al. \cite{collins24} and an extension to the orthogonal case, characterizes the
fluctuation. For instance, when $\theta_\ell=\alpha_i$,
\begin{equation}
    \label{heuristic}
      p\bigl(\lambda_{p,\ell}-\sigma^2-\theta_\ell\bigr)
  \xrightarrow{\ \mathrm{law}\ } W_\ell
  =\sum_{j=1}^s\frac{\alpha_i\beta_j}{\alpha_i-\beta_j}\,|g_{ij}|^2 ,
\end{equation}
where the $g_{ij}$ are the limits in law of the $\sqrt p\,c_{ij}^{(p)}$, and are
independent standard Gaussian variables, real or circular complex. In particular the
$|g_{ij}|^2$ are i.i.d. $\chi_1^2$ variables if $\mathbb F=\mathbb R$, and i.i.d.
exponential $\text{Exp}(1)$ variables if $\mathbb F=\mathbb C$. A similar result holds
when $\theta_\ell=\beta_k$ for some $k\leq s$, see Proposition~\ref{p:pop} for details.

\textbf{The sample fluctuation.}
Conditionally on $\mathcal F_p=\sigma(O_p)$, the matrix $\Sigma_p$ is fixed and
the only randomness left is the one of the samples. We are then in a spiked model with
spikes $\lambda_{p,\ell}-\sigma^2$, and the Bai--Yao CLT \cite{bai12} presented in
Section~\ref{ss:spiked} suggests that for any $\ell\leq r+s$,
$\widehat\lambda_{p,\ell}$ fluctuates at the speed $n^{-1/2}$ around
\[
  \phi_{p/n}\bigl(\lambda_{p,\ell}-\sigma^2\bigr)
  =\lambda_{p,\ell}
  +\frac{p}{n}\,\frac{\sigma^2\lambda_{p,\ell}}{\lambda_{p,\ell}-\sigma^2} ,
\]
and that the law of the fluctuation is Gaussian. However, let us emphasize that the
Bai--Yao CLT is stated for fixed spikes, while here the spike $\lambda_{p,\ell}$
changes with $p$, so this reasoning is heuristic. We nevertheless follow this heuristic. Writing
$\phi_{p/n}(\theta)=\sigma^2+\theta+p\sigma^2/n+\frac{p\sigma^4}{n\theta}$, one has the identity
\begin{equation}
  \phi_{p/n}\bigl(\lambda_{p,\ell}-\sigma^2\bigr)-\phi_{p/n}(\theta_\ell)
  =\bigl(\lambda_{p,\ell}-\sigma^2-\theta_\ell\bigr)
   \big(1-\frac{p\,\sigma^4}
   {n(\lambda_{p,\ell}-\sigma^2)\,\theta_\ell}\big)
   = \mathrm O_{\mathbb P}\left(\frac{1}{p}\right),
  \label{e:phidiff}
\end{equation}
which is a $\mathrm O_{\mathbb P}(1/p)$ because it corresponds to the non-scaled Haar
fluctuation which scales at $p$ by (\ref{heuristic}), multiplied by a bounded term when
$\ell\leq r+s$. Therefore we can write the informal, heuristic decomposition
\begin{equation}
  \widehat\lambda_{p,\ell}-\phi_{p/n}(\theta_\ell)
  =\underbrace{\Bigl(\widehat\lambda_{p,\ell}
     -\phi_{p/n}\bigl(\lambda_{p,\ell}-\sigma^2\bigr)\Bigr)}
   _{\mathrm O_{\mathbb P}(1/\sqrt n),\ \text{sample fluctuation}}
  +\underbrace{\Bigl(\phi_{p/n}\bigl(\lambda_{p,\ell}-\sigma^2\bigr)
     -\phi_{p/n}(\theta_\ell)\Bigr)}
   _{\mathrm O_{\mathbb P}(1/p),\ \text{Haar fluctuation}} .
  \label{e:heur}
\end{equation}
When $p/n\to c>0$, the Haar fluctuation is of order $1/p$, hence of order $1/n$, which
is negligible compared to the sample fluctuation, of order $1/\sqrt n$. In this case
the Haar randomness does not impact the limit, and the sample eigenvalues of our model
behave as in the Bai--Yao theorem, because the spike $\lambda_{p,\ell}-\sigma^2$
converges fast enough to $\theta_\ell$. This is Theorem~\ref{t:prop}.

When $p/n\to0$, the Bai--Yao CLT does not apply, since it assumes $c>0$. However, we
prove in Proposition~\ref{p:dclt} a specialized CLT for a deterministic spiked model
with $p/n\to0$. Heuristically, Proposition~\ref{p:dclt} suggests that $\hat\lambda_{p,\ell}$ still fluctuates around $\phi_{p/n}\bigl(\lambda_{p,\ell}-\sigma^2\bigr)$, as in the case where $c>0$. Hence, the identity \eqref{e:phidiff} still holds when $c =0$ and the two terms keep the same orders. The behaviour of
$\widehat\lambda_{p,\ell}-\phi_{p/n}(\theta_\ell)$ is therefore decided by the size of
$1/\sqrt n$ compared to $1/p$.

\textbf{A transition at $p^2\asymp n$.} The two fluctuations have different sizes: the
Haar one is of size $1/p$, the sample one of size $1/\sqrt n$, and they are
independent. To obtain a limit which is neither zero nor divergent, one has to multiply
by the inverse of the larger of the two, that is, by $\min(p,\sqrt n)$. Comparing $1/p$
with $n^{-1/2}$ amounts to comparing $p^2$ with $n$, so the relevant parameter is
$p^2/n$ and not $p/n$. We stress this: the ratio $p/n$ does not distinguish
the three disproportional regimes, since it tends to zero in all three. $p^2/n$ it the
ratio that governs the behavior of fluctuations, and, as we will see, the limiting law
is different for $p^2/n\to\infty$, $p^2/n\to\tau\in(0,\infty)$, and $p^2/n\to0$.

No single normalization covers the three disproportional regimes.  If $p^2/n\to0$ then $\sqrt n/p\to\infty$, so $\sqrt
n(\widehat\lambda_{p,\ell}-\phi_{p/n}(\theta_\ell))$ diverges in probability. If $p^2/n\to\infty$, the same happens with the scaling $p$. Both work only when \mbox{$p^2/n\to\tau\in(0,\infty)$}, and
they then differ by the factor $\sqrt\tau$. In what follows, we work with
the normalization $p$ when $\tau$ is finite, because it is the only one
of the two which gives a single formula for the whole range $\tau\in[0,\infty)$,
 and with $\sqrt n$ when $p^2/n\to\infty$.

\subsection{Main results}
\label{ss:main}

We describe the fluctuations of $\widehat\lambda_{p,\ell}$, for each fixed spike
$\theta_\ell$, in four regimes of $(p,n)$. Recall from Section~\ref{ss:spiked} that
\[
  \phi_c(\theta)=(\sigma^2+\theta)\Bigl(1+\frac{c\,\sigma^2}{\theta}\Bigr),
  \qquad
  V_c(\theta)=(\sigma^2+\theta)^2\Bigl(1-\frac{c\,\sigma^4}{\theta^2}\Bigr),
\]
so that $\phi_0(\theta)=\sigma^2+\theta$ and $V_0(\theta)=(\sigma^2+\theta)^2$. The
centering is the same function in the four regimes: it is $\phi_{p/n}(\theta_\ell)$ in all cases. Let $G_\ell\sim\mathcal
N(0,1)$ be independent of $W_\ell$. The results are collected in Table~\ref{t:four}.

\begin{table}[H]
\caption{The four regimes, at a fixed spike $\theta_\ell$. Here $\kappa=2$ if $\mathbb
F=\mathbb R$ and $\kappa=1$ if $\mathbb F=\mathbb C$; the variable $W_\ell$ is the Haar
fluctuation of Section~\ref{ss:why}  and $G_\ell\sim\mathcal N(0,1)$ is independent of
it. In the first row we assume in addition that $\theta_\ell$ is supercritical, that
is, $\theta_\ell>\sigma^2\sqrt c$; in the three other rows this holds for $p$ large,
since the threshold $\sigma^2\sqrt{p/n}$ tends to zero. The word ``diverges'' means
that the absolute value of the scaled fluctuation tends to infinity in probability, so
that no subsequence of it converges in distribution; see Proposition~\ref{p:diverge}.}
\label{t:four}
\centering
\footnotesize
\setlength{\tabcolsep}{1pt}
\begin{tabular}{@{}lll@{}}
\toprule
Regime
& Scaling by $\sqrt n$, and limit
& Scaling by $p$, and limit \\
\midrule
$p/n\to c>0$
& $\sqrt n\,(\widehat\lambda_{p,\ell}-\phi_{p/n}(\theta_\ell))
   \lawto\sqrt{\kappa V_c(\theta_\ell)}\,G_\ell$
& $p\,(\widehat\lambda_{p,\ell}-\phi_{p/n}(\theta_\ell))$ diverges \\[1mm]
$p/n\to0$, $p^2/n\to\infty$
& $\sqrt n\,(\widehat\lambda_{p,\ell}-\phi_{p/n}(\theta_\ell))
   \lawto\sqrt{\kappa V_0(\theta_\ell)}\,G_\ell$
& $p\,(\widehat\lambda_{p,\ell}-\phi_{p/n}(\theta_\ell))$ diverges \\[1mm]
$p^2/n\to\tau\in(0,\infty)$
& $\sqrt n\,(\widehat\lambda_{p,\ell}-\phi_{p/n}(\theta_\ell))
   \lawto\sqrt{\kappa V_0(\theta_\ell)}\,G_\ell+\tau^{-1/2}W_\ell$
& $p\,(\widehat\lambda_{p,\ell}-\phi_{p/n}(\theta_\ell))
   \lawto W_\ell+\sqrt{\kappa\tau V_0(\theta_\ell)}\,G_\ell$ \\[1mm]
$p^2/n\to0$
& $\sqrt n\,(\widehat\lambda_{p,\ell}-\phi_{p/n}(\theta_\ell))$ diverges
& $p\,(\widehat\lambda_{p,\ell}-\phi_{p/n}(\theta_\ell))\lawto W_\ell$ \\
\bottomrule
\end{tabular}
\end{table}

The first row covers the proportional regime $p/n\to c>0$: in this regime, for a supercritical spike $\theta_\ell$, the classical
Bai--Yao CLT holds in our Haar spiked model the same way as in the classical
deterministic spiked model of \cite{johnstone01}. The second row covers the regime $p/n \to 0$ and $p^2/n\to\infty$; in this regime, $\sqrt n\,(\widehat\lambda_{p,\ell}-\phi_{p/n}(\theta_\ell))$ converges to the same gaussian distribution $\sqrt{\kappa V_0(\theta_\ell)}\,G_\ell$ as in the proportional regime (setting $c=0$ in the first row yields the same result as in the second row). This is because in both regimes, the Haar fluctuation is of size $1/p$, which is much smaller than the sample fluctuation of size $1/\sqrt
n$; the Haar matrix leaves no trace in the limit and only the sample fluctuation impacts the limiting distribution. However, we keep the two rows separated because the two regimes differ in terms of proof techniques that we use to show the results of Table \ref{t:four} and in terms of behavior of the rest of the spectrum. In particular, we show in Theorem \ref{t:prop} that the Marcenko--Pastur bulk, the BBP phase transition, the Tracy--Widom edge properties still hold in the proportional regime.

The three last rows cover the three possible cases in the disproportional regime
$p/n\to0$. Each column has one term which does not change, and one term which
appears at $p^2/n\to\tau\in(0,\infty)$. Under the scaling by $\sqrt n$, the Gaussian
term $\sqrt{\kappa V_0(\theta_\ell)}G_\ell$ is the same in the second and the third
row, and it does not depend on the rate at which $p/n$ goes to zero; what appears in
the third row is the Haar fluctuation $W_\ell$, scaled by $\tau^{-1/2}$.
Under the scaling by $p$, the variable $W_\ell$ is the same in the third and the fourth
row, and what appears in the third row is the Gaussian term, scaled by
$\sqrt\tau$. So the mixed regime $p^2/n \to \tau \in (0,\infty)$ is the only case where the two fluctuations of
Section~\ref{ss:why} have the same size, and the only one where both appear in the
limit.

Notice that the limit law in the mixed regime is a convolution of gaussian and
exponential (if $\mathbb F=\mathbb C$) or chi-squared (if $\mathbb F=\mathbb R$)
variables. Indeed, the variable $W_\ell$ comes from the Haar matrix, the variable
$G_\ell$ comes from the sample vectors, these two sources are independent before the
limit, and they stay independent in the limit.

Each column also interpolates towards its neighbour: under the scaling
by $p$, letting $\tau\to0$ in the third row gives the fourth row, and under the scaling
by $\sqrt n$, letting $\tau\to\infty$ makes the coefficient $\tau^{-1/2}$ vanish and
gives the second row. We come back to this in Section~\ref{ss:scal}. Let us emphasize that the fourth row is the case $\tau=0$ of the third, but we display it separately for two
reasons: the limiting law is different, since it loses its Gaussian part and reduces to
$W_\ell$ alone; and the scaling by $\sqrt n$, which still gives a limit in the third row,
diverges here, so the two rows do not merge under that scaling.

We state these results formally in Section~\ref{s:main results}, once the model and the
labeling of the eigenvalues have been set up in Section~\ref{s:mod}, and we prove them
in Section~\ref{s:prf} and in the appendices.

\subsection{An application: linear mixed models with a continuous batch effect}
\label{ss:lmm}

Our model describes a real-world problem and allows us to relax a constraint in linear
mixed models. Linear mixed models (LMM) are used when the data come in batches \cite{laird82,searle92}, i.e.,
in groups of observations which share something: samples processed on the same day,
measurements taken during the same time window,... One observes $n$ vectors of $\mathbb
F^p$, split into $G$ batches, and writes each sample $k$ as
\[
  x_{p,k}=s_{p,k}+m_{p,z_k}+\varepsilon_{p,k},
\]
where $z_k\in[G]$ is the batch of sample $k$, and the three terms are independent
centered Gaussian vectors: the effect of interest $s_{p,k}\sim\mathcal N(0,A_p)$, drawn
independently for each sample; the batch effect $m_{p,g}\sim\mathcal N(0,\Gamma_p)$, drawn once per batch and
shared by all the samples of that batch; and the measurement noise
$\varepsilon_{p,k}\sim\mathcal N(0,\sigma^2I_p)$. The eigenvalues of $A_p$ and of
$\Gamma_p$ are what one wants to estimate. The covariance of one sample is
$A_p+\Gamma_p+\sigma^2I_p$, so it does not separate $A_p$ from $\Gamma_p$. What makes
the separation possible is that two samples $k\neq k'$ of the same batch share the same batch effect, while their other terms are independent, which
gives
\[
  \mathbb E\bigl[x_{p,k}x_{p,k'}^*\bigr]=\Gamma_p .
\]
One therefore estimates $\Gamma_p$ by averaging
$x_{p,k}x_{p,k'}^*$ over the pairs of samples in the same batch and recovers $A_p$ by subtraction.

A limit of the classical LMM is that it requires batches of size at least two, which does not hold when the batch effect varies slightly from one sample to the next instead of being shared by a group: for instance a slow drift of an instrument
between two measurements. We
then call the batch effect continuous, each batch is of size one, and there is no pair to average over:
the only available object is the sample covariance matrix $\widehat\Sigma_p$, which
estimates the sum $A_p+\Gamma_p+\sigma^2I_p$, so $A_p$ and $\Gamma_p$
cannot be recovered.

Our model keeps the covariance structure of the classical LMM, allows a single point
per batch, and adds in return one assumption on $\Gamma_p$: its eigenvectors have no
preferred orientation, and are unrelated to those of $A_p$. As explained in
Section~\ref{ss:model}, this is the Haar assumption, that is, $\Gamma_p=O_pB_pO_p^*$,
and the model becomes \eqref{e:sigm}. It is what makes the separation of $A_p$ and
$\Gamma_p$ possible again: the eigenvectors of $A_p$ and $\Gamma_p$ are asymptotically
orthogonal, and the spikes of $\Sigma_p$ are close to the union of the $\alpha_i$ and
the $\beta_j$ rather than to a mixture of them. In practice, this assumption requires
that the batch effect be caused by a mechanism that is unrelated to the effect of
interest, and it fails when the two are confounded.

Batch effects of this kind are well documented in genomics \cite{leek10,johnson07},
where $x_{p,k}$ collects the expression levels of $p$ genes for one biological
sample. Typically, $p\approx10^3$--$10^4$ and
$n\approx10^2$, so  $p/n\to c>0$ in this setup. In industrial
quality control, $x_{p,k}$ collects $p$ measured characteristics of one manufactured
item, and the batch effect is the distortion added by the measuring instrument, which
changes from one item to another. Typically, $p\approx10^1$ and $n\approx10^3$, so $p/n\to0$ while $p^2/n\to\tau<\infty$. The four regimes of
Table~\ref{t:four} therefore correspond to different application areas.

Section~\ref{ss:test} develops the industrial quality control application, states the assumptions
under which the model is reasonable, discusses its limits, and derives a test for a spike from Table~\ref{t:four}.

\subsection{Related work}
\label{ss:related}

Our work is related to spiked covariance models, Haar matrix models,
and sample covariance matrices with random populations.

\textbf{The spiked covariance model with a deterministic population matrix.} The proportional regime is by now classical \cite{johnstone01,bbp05,baik06,paul07,bai08,bai12}, with a comprehensive treatment provided by \cite{bai10,yao15}. Bloemendal et al. \cite{bkyy16} relax the assumption of a converging $p/n$ ratio. They establish large deviation bounds on outlier locations and eigenvector overlaps under the weaker condition that $p$ and $n$ are polynomially related, in a generalization of the spiked model. The disproportional regime where $p, n \to \infty$ with $p/n \to 0$ is also central to this model. 
Johnstone and Lu \cite{johnstone2009sparseprincipalcomponentsanalysis} showed that the leading sample principal component is consistent
if and only if $p/n\to0$ in their asymptotic model. 
In the same regime, Donoho and Feldman \cite{df26} derived optimal shrinkage rules.
Liu et al.~\cite{liu2022rds} proved a CLT for spiked sample eigenvalues without assuming convergence of $p/n$, subject to a condition ensuring that the spikes remain distant. These works consider deterministic population matrices. In our model, such results are relevant after conditioning on the Haar matrix; we prove the precise low-aspect-ratio statement needed for our conditional argument in Proposition~\ref{p:dclt}.

\textbf{The spiked invariant model.} The matrix $\Sigma_p =A_p+O_pB_pO_p^*+\sigma^2I_p$ in~\eqref{e:sigm} is a sum of matrices, one of which is conjugated by a Haar matrix. Models involving polynomials of diagonal matrices and Haar unitary/orthogonal matrices have been extensively studied. Voiculescu \cite{voiculescu91} established that the global limiting spectral distribution of such additive models is governed by free convolution.  Bao et al. \cite{bao2016convergenceratespectraldistribution} obtained quantitative convergence rates for the empirical spectral distribution in additive Haar models. For the multiplicative version of this model involving matrices $A_p,B_p,O_p$, Ding and Ji \cite{ding2023spikedmultiplicativerandommatrices} established the limits and optimal convergence rates for both the outlier eigenvalues and their corresponding eigenvectors. For fixed-rank perturbations, the normalized limiting spectral distribution does not describe the individual spikes. Collins et al. \cite{chs18} studied the discrete limiting spectrum in a compact-operator framework. Fluctuations for finite-rank polynomial models in the Haar unitary case, under a simplicity assumption on the limiting eigenvalues, were computed in \cite{collins24}. These results provide the background for the population fluctuation analysis in our model.

\textbf{Fluctuations of outliers.}
Benaych-Georges et al. \cite{benaychgeorges2011fluctuationsextremeeigenvaluesfinite} studied extreme eigenvalues of additive finite-rank deformations of deterministic matrices, with extensions to random matrix ensembles including Wishart matrices. In their orthonormalised Gaussian setting, the perturbation has the form $O_pB_pO_p^*$. In our model, this Haar deformation enters the population covariance before the conditional Gaussian sampling step. Non-Gaussian outlier limits are also known for deformed Wigner matrices \cite{cdf09}, where they can depend on the entry distribution. Here the non-Gaussian component comes from the Haar-induced population spike shifts.

\textbf{Random population covariance matrices.} Sample covariance matrices with random population covariance have also been studied. Silverstein \cite{silverstein95} established almost sure convergence of the empirical spectral distribution in this setting. Kwak, Lee and Park \cite{kwak21} studied extremal eigenvalues for random diagonal population matrices, while Lee and Li \cite{leeli22} proved a CLT for linear spectral statistics in a random diagonal population model. Our focus is on individual sample outliers generated by fixed-rank Haar perturbations of a scalar covariance matrix.

\textbf{Our contribution.} For fixed ranks and pairwise distinct fixed spikes, our results describe the interaction between Haar-induced population spike shifts of order $p^{-1}$ and Gaussian sampling fluctuations of order $n^{-1/2}$. We center each sample outlier at $\phi_{p/n}(\theta_\ell)$, which includes the deterministic sampling bias. In the low-aspect-ratio regime $p/n\to0$, the ratio $p^2/n$ determines the fluctuation regime. When $p^2/n$ tends to zero, a positive finite limit, or infinity, the appropriately normalized fluctuation of each fixed sample outlier has a Haar limit, a convolution of the Haar and Gaussian laws, or a Gaussian limit, respectively. The Haar limit is a linear combination of independent chi-square variables with one degree of freedom in the real case, and independent unit-mean exponential variables in the complex case. The Haar and Gaussian components of the mixed limit are independent. We also recover the Gaussian fluctuation law for each supercritical spike in the proportional regime.

\subsection{Organization of the paper}
\label{ss:org}

The rest of the paper is organized as follows. Section~\ref{s:mod}
states the model, sets the notation, and recalls the results
on Haar matrices and on the classical spiked model that we use as
inputs. It also fixes the labeling of the eigenvalues and defines the
asymptotic regimes. Section~\ref{s:main results} contains the main
results. We first treat the disproportional
regime $p/n\to0$, which is the main object of the paper. In the
disproportional regime the ratio $p^2/n$ separates three cases: the Haar-dominated regime $p^2/n\to0$, the
mixed regime $p^2/n\to\tau\in(0,\infty)$, both treated in
Theorem~\ref{t:main}, and the sample-dominated regime
$p^2/n\to\infty$, treated in Theorem~\ref{t:infy}. We then treat the proportional regime $p/n \to c>0$, where the classical Bai--Yao CLT applies. Next, we state the
results on the two scalings, give a rank-one example, and present applications of these results and a spectral test. Section~\ref{s:prf} proves these
results. All the proofs share one first step, the reduction to a
deterministic spiked model conditionally on the Haar matrix; after
this step each proof is specific to its regime.
Appendix~\ref{s:appfacts} collects the facts used several times
in the paper, and the more technical proofs are given in
Appendices~\ref{appendix: proof of p:pop}
and~\ref{appendix: proof of p:dclt}.

\section{Model, notation and preliminaries}\label{s:mod}

In this section, we recall the model introduced in Section~\ref{ss:model}, we introduce important definitions and conventions, we state important results on Haar matrices and we formalize the foundational theorems of random matrix theory mentioned in Section~\ref{ss:spiked}.

\subsection{The model}\label{ss:setup}

Throughout the paper we work over a field $\mathbb F$, which is either $\mathbb R$ or
$\mathbb C$, fixed once and for all. The choice of $\mathbb F$ ties together three
objects: the entries of the matrices, the group from which the Haar matrix is drawn, and
the law of the samples. We use the following terminology.

\begin{itemize}
  \item In the real case, $\mathbb F=\mathbb R$: the matrices $A_p$, $B_p$ have
        real entries, $O_p$ is Haar distributed on the orthogonal group $O(p)$, and the
        samples are real Gaussian vectors.
  \item In the complex case, $\mathbb F=\mathbb C$: the matrices $A_p$, $B_p$ have
        complex entries, $O_p$ is Haar distributed on the unitary group $U(p)$, and the
        samples are circular complex Gaussian vectors.
\end{itemize}

The groups $O(p)$ and $U(p)$, the Haar measure, the circular complex Gaussian law, and
the constant $\kappa$ ($\kappa=2$ in the real case, $\kappa=1$ in the complex case) are
defined in Section~\ref{s:grp}. We write $M^*$ for the transpose in the real case and
the conjugate transpose in the complex case.

Fix two integers $r,s\ge1$ and two families of strictly positive numbers
\[
  \alpha_1,\ldots,\alpha_r,\qquad \beta_1,\ldots,\beta_s,
\]
which we assume to be all distinct, and fix a number $\sigma^2>0$. For every
$p\ge r+s$, let $A_p$ and $B_p$ be deterministic positive semidefinite $p\times p$
matrices over $\mathbb F$, of rank $r$ and $s$ respectively, whose nonzero eigenvalues
are $\alpha_1,\ldots,\alpha_r$ for $A_p$ and $\beta_1,\ldots,\beta_s$ for $B_p$. Let
$O_p$ be Haar distributed on the orthogonal group $O(p)$ when $\mathbb F=\mathbb R$, and
on the unitary group $U(p)$ when $\mathbb F=\mathbb C$. The population covariance matrix
is
\begin{equation}
  \Sigma_p=A_p+O_pB_pO_p^*+\sigma^2I_p .
  \label{e:sigm2}
\end{equation}

Conditionally on $O_p$, we observe $n$ independent centered Gaussian vectors with
covariance $\Sigma_p$,
\[
  x_{p,1},\ldots,x_{p,n}\mid O_p
  \stackrel{\mathrm{i.i.d.}}{\sim}\mathcal N_{\mathbb F}(0,\Sigma_p),
\]
where $\mathcal N_{\mathbb F}(0,\Sigma_p)$ is the centered real Gaussian law with
covariance $\Sigma_p$ when $\mathbb F=\mathbb R$, and the centered circular complex
Gaussian law with covariance $\Sigma_p$ when $\mathbb F=\mathbb C$. From these vectors
we form the sample covariance matrix
\[
  \widehat\Sigma_p=\frac1n\sum_{k=1}^n x_{p,k}x_{p,k}^* .
\]

We order the eigenvalues of $\Sigma_p$ and of $\widehat\Sigma_p$ decreasingly,
\[
  \lambda_{p,1}\ge\cdots\ge\lambda_{p,p},
  \qquad
  \widehat\lambda_{p,1}\ge\cdots\ge\widehat\lambda_{p,p},
\]
and we write $\theta_1>\cdots>\theta_{r+s}>0$ for the elements of
$\{\alpha_1,\ldots,\alpha_r,\beta_1,\ldots,\beta_s\}$ sorted decreasingly. We call the
$\theta_\ell$ the spikes of the model. As shown in Section~\ref{s:diag}, the matrix
$A_p+O_pB_pO_p^*$ has rank $r+s$ almost surely, so $\Sigma_p$ has $r+s$ eigenvalues
different from $\sigma^2$, called the spike eigenvalues, and the eigenvalue $\sigma^2$
with multiplicity $p-r-s$. The precise correspondence between a population spike
eigenvalue $\lambda_{p,\ell}$ and the associated sample spike eigenvalue
$\widehat\lambda_{p,\ell}$ is fixed in Section~\ref{s:lab}.

\subsection{Groups, Haar measure, and conventions}\label{s:grp}

For $p\ge1$ we write
\[
  O(p)=\{Q\in M_p(\mathbb R):Q^\top Q=I_p\},
  \qquad
  U(p)=\{Q\in M_p(\C):Q^*Q=I_p\},
\]
and we denote by $O(p)$ the unitary (orthogonal) group when $\mathbb F=\mathbb R$, and the
group $U(p)$ when $\mathbb F=\C$. All the matrices $O_p$, $p\ge1$, are defined on the same probability space, so that
almost sure statements about the sequence $(\lambda_{p,\ell})_{p\ge1}$ make sense. Since
our results are convergences in law, they do not depend on the joint law chosen for the
$O_p$, and the reader may take them independent.

$O_p$ and $U_p$ are compact groups. On such a group there is
a unique Borel probability measure which is invariant under left multiplication by a
fixed element of the group; it is called the Haar measure, and compactness implies
that it is also invariant under right multiplication and under inversion. A random
matrix is Haar distributed when its law is this measure. This is the
precise meaning of the informal description given in Section~\ref{ss:model}: no
orientation is more likely than another, and a Haar orthogonal (resp.unitary) matrix is a matrix sampled uniformly on $O(p)$ (resp. $U(p)$).

In the real case, $M^*$ denotes the transpose $M^\top$; in the complex case, it
denotes the conjugate transpose. We use the same notation in both cases.

A standard real Gaussian variable is $\mathcal N(0,1)$. A standard circular complex
Gaussian variable is $z=x+\mathrm{i}y$, where $x$ and $y$ are independent
$\mathcal N(0,1/2)$; hence $\mathbb E|z|^2=1$ and $|z|^2\sim\operatorname{Exp}(1)$.
Recall from Section~\ref{ss:model} that
\begin{equation}
  \kappa=
  \begin{cases}
    2,&\text{in the real case},\\
    1,&\text{in the complex case}.
  \end{cases}
  \label{e:kapp}
\end{equation}
With these conventions, a standard real or complex Gaussian variable $z$
satisfies $\operatorname{Var}(|z|^2-1)=\kappa$. 

We write $\mathcal F_p=\sigma(O_p)$ for the sigma-algebra generated by the Haar matrix. When we condition on $\mathcal F_p$, the
notation refers to the remaining Gaussian randomness.

If $X_p$ and $X$ are random variables, then $X_p\xrightarrow{a.s.}X$ means almost sure convergence, $X_p\xrightarrow{\mathbb P}X$ means
convergence in probability, and $X_p\xrightarrow{\ \mathrm{law}\ }X$ means
convergence in distribution. If $a_p>0$, then $X_p=o_{\mathbb P}(a_p)$ means
$X_p/a_p\xrightarrow{\mathbb P}0$, and $X_p=\mathrm O_{\mathbb P}(a_p)$ means that
$\{X_p/a_p\}_{p\ge1}$ is stochastically bounded; that is, for any $\epsilon > 0$, there exists a constant $M > 0$ such that $\sup_{p \ge 1} \mathbb{P}(\vert{}X_p\vert{}/a_p > M) < \epsilon$. We say that a sequence of events $\Omega_p$ holds
with high probability when $\mathbb P(\Omega_p)\to1$. We also use uniform versions of
this notation: if $I$ is a fixed compact interval, then
$\sup_{t\in I}|X_p(t)|=o_{\mathbb P}(a_p)$ means that the non-negative random
variable on the left, divided by $a_p$, converges to zero in probability.

Finally, we say that a sequence $X_p$ of real random variables diverges in
probability, and we write $|X_p|\xrightarrow{\mathbb P}\infty$, when
$\mathbb P(|X_p|\le M)\to0$ for every constant $M>0$. By Lemma~\ref{a:nolimit}, such a
sequence is not stochastically bounded, and no subsequence of it converges in distribution.

\subsection{Reduction to the diagonal case and facts on $\Sigma_p$}\label{s:diag}

\begin{rem}[Reduction to the diagonal case]\label{r:wlog}
In the rest of the paper, we may and do assume without loss of generality that
\[
    A_p=\diag(\alpha_1,\ldots,\alpha_r,0,\ldots,0),
    \qquad
    B_p=\diag(\beta_1,\ldots,\beta_s,0,\ldots,0).
\]
Indeed, write $A_p=UCU^*$ and $B_p=VDV^*$, where $C$ and $D$ are the two diagonal
matrices above, and where $U,V\in U(p)$ are deterministic. Put
$\widetilde O_p=U^*O_pV$ and $\widetilde x_{p,k}=U^*x_{p,k}$. By the left and right
invariance of the Haar measure, $\widetilde O_p$ is again Haar distributed, and,
conditionally on $O_p$, the vectors $\widetilde x_{p,1},\ldots,\widetilde x_{p,n}$
are i.i.d. with law $\mathcal N_\mathbb F(0,\widetilde\Sigma_p)$, where
\[
  \widetilde\Sigma_p
  :=U^*\Sigma_pU
  =C+\widetilde O_pD\widetilde O_p^*+\sigma^2I_p .
\]
Thus $(\widetilde O_p,\widetilde x_{p,1},\ldots,\widetilde x_{p,n})$ has exactly the
law of the data of the diagonal model. Moreover, for each realization of $O_p$ and of
the samples,
\[
  \frac1n\sum_{k=1}^n\widetilde x_{p,k}\widetilde x_{p,k}^*=U^*\widehat\Sigma_pU ,
\]
and a conjugation by $U$ does not change eigenvalues. Hence the vector
$(\lambda_{p,\ell},\widehat\lambda_{p,\ell})_{\ell\le r+s}$ is the same random
variable in the two models, and its law is the same for every $p$. Only the nonzero
eigenvalues of $A_p$ and $B_p$ enter our results.
\end{rem}

From now on, $A_p$ and $B_p$ are the two diagonal matrices of Remark~\ref{r:wlog}. In
particular, the eigenvectors of $A_p$ are $e_1,\ldots,e_r$, and the eigenvectors of
$O_pB_pO_p^*$ are the first $s$ columns of $O_p$.

Both $A_p$ and $O_pB_pO_p^*$ are positive semidefinite, so $\Sigma_p$ is positive
definite, and $\|\Sigma_p\|_{\mathrm{op}}\le\sigma^2+\max_i\alpha_i+\max_j\beta_j$ for
every $p$. By Proposition~\ref{p:pop}, the spikes $\lambda_{p,\ell}-\sigma^2$ converge in
probability to the distinct positive limits $\theta_\ell$, $\ell=1,\ldots,r+s$; hence,
with probability tending to one, $\Sigma_p$ has $r+s$ eigenvalues different from
$\sigma^2$, which we call the spike eigenvalues, and the eigenvalue $\sigma^2$ with
multiplicity at least $p-r-s$.

\subsection{Two properties of Haar matrices}\label{s:quo}

The following result is the fixed-block Gaussian limit for Haar matrices. In the real
case it is a direct consequence of Jiang \cite[Theorem~1]{jiang06}, which gives
convergence in total variation for blocks of size $r\times s$ with $r,s$ fixed. The
complex case is stated in Collins et al. \cite[Lemma~2.1]{collins24}.

\begin{prop}[Fixed-block limit for Haar matrices]\label{p:haar}
Let $r,s$ be fixed, and let $O_p=(o_{ij}^{(p)})$ be Haar distributed on $O(p)$ or
$U(p)$. Then
\[
  \bigl(\sqrt p\,o_{ij}^{(p)}\bigr)_{1\leq i\leq r,\,1\leq j\leq s}
  \xrightarrow{\ \mathrm{law}\ }
  (g_{ij})_{1\leq i\leq r,\,1\leq j\leq s},
\]
where the entries of $(g_{ij})$ are i.i.d. standard real Gaussian variables in the real case, and i.i.d. standard circular complex Gaussian variables in the complex case.
\end{prop}

The next result gives the Haar fluctuation of the population spike eigenvalues, at
the scale $1/p$. The complex case is a specialization of a result shown by Collins et
al. \cite[Theorem~ 1.2 and 1.3]{collins24}. We prove the result for the real case in Appendix \ref{appendix: proof of p:pop}, using the same fixed-block perturbation computation as in the complex case, and then using the real case of Proposition~\ref{p:haar}. 
\begin{prop}[Population spike fluctuations at scale $1/p$]\label{p:pop} Recall that the $\alpha_i$ and $\beta_j$ are all distinct. For each $\ell\in\{1,\ldots,r+s\}$,
\begin{equation}
    \lambda_{p,\ell}\xrightarrow{\ \mathrm{a.s.}\ } \sigma^2+\theta_\ell,
\end{equation}
\begin{equation}
   p\bigl(\lambda_{p,\ell}-\sigma^2-\theta_\ell\bigr)
  \xrightarrow{\ \mathrm{law}\ } W_\ell=
  \begin{cases}
  \displaystyle
  \sum_{j=1}^s
  \frac{\alpha_i\beta_j}{\alpha_i-\beta_j}|g_{ij}|^2,
  &\theta_\ell=\alpha_i,\\[3mm]
  \displaystyle
  \sum_{i=1}^r
  \frac{\alpha_i\beta_j}{\beta_j-\alpha_i}|g_{ij}|^2,
  &\theta_\ell=\beta_j.
  \end{cases}
  \label{e:well}
\end{equation}
where $|g_{ij}|^2\stackrel{\mathrm{i.i.d.}}{\sim}\chi_1^2$ in the real case, and
$|g_{ij}|^2\stackrel{\mathrm{i.i.d.}}{\sim}\operatorname{Exp}(1)$ in the complex case.
Moreover, the convergence in law holds jointly for $\ell=1,\ldots,r+s$, with all $W_\ell$ defined using the same family $(g_{ij})$.
\end{prop}

\subsection{Known results for the classical spiked model}\label{s:quoted}

The three results of this subsection concern the classical spiked covariance model, in
which the population covariance is deterministic. We state them as external inputs and we use them extensively to prove our results in Section~\ref{s:prf}.

Fix an integer $K\ge1$, and let
\[
  \Lambda_p^{\mathrm{id}}
  =\diag(\sigma^2+\theta_1,\ldots,\sigma^2+\theta_K,\sigma^2,\ldots,\sigma^2),
\]
where $\theta_1>\cdots>\theta_K>0$ are fixed. Let $y_{p,1},\ldots,y_{p,n}$ be i.i.d.
$\mathcal N_{\mathbb F}(0,\Lambda_p^{\mathrm{id}})$, and let
$\widehat\lambda_{p,1}^{\mathrm{id}}\ge\cdots\ge\widehat\lambda_{p,p}^{\mathrm{id}}$ be
the ordered eigenvalues of the associated sample covariance matrix. Recall from
Section~\ref{ss:main} the functions
\begin{equation}
    \label{eq:definition of phi_c and V_c}
    \phi_c(\theta)=(\sigma^2+\theta)\Bigl(1+\frac{c\,\sigma^2}{\theta}\Bigr),
  \qquad
  V_c(\theta)=(\sigma^2+\theta)^2\Bigl(1-\frac{c\,\sigma^4}{\theta^2}\Bigr).
\end{equation}
A spike $\theta_\ell$ is called supercritical when $\theta_\ell>\sigma^2\sqrt c$, and subcritical when $\theta_\ell\leq\sigma^2\sqrt c$.

\begin{thm}[Limit of the spikes, BBP transition \cite{baik06}]\label{t:bbp}
Assume $p/n\to c\in(0,\infty)$. For each $\ell\le K$: if $\theta_\ell>\sigma^2\sqrt c$,
then $\widehat\lambda_{p,\ell}^{\mathrm{id}}\to\phi_c(\theta_\ell)$ almost surely, and
$\phi_c(\theta_\ell)>\sigma^2(1+\sqrt c)^2$; if $0<\theta_\ell\le\sigma^2\sqrt c$, then
$\widehat\lambda_{p,\ell}^{\mathrm{id}}\to \sigma^2(1+\sqrt c)^2$ almost surely.
\end{thm}

\begin{thm}[Gaussian fluctuations of a supercritical spike \cite{bai08,bai12}]\label{t:baiyao}
Assume $p/n\to c\in(0,\infty)$. Fix $\ell\le K$ and assume $\theta_\ell$ is
supercritical. Then
\[
  \sqrt n\,\bigl(\widehat\lambda_{p,\ell}^{\mathrm{id}}-\phi_{p/n}(\theta_\ell)\bigr)
  \xrightarrow{\ \mathrm{law}\ }
  \mathcal N\bigl(0,\kappa V_c(\theta_\ell)\bigr),
\]
with $\kappa$ given by \eqref{e:kapp}. The fluctuations of distinct supercritical spikes
are asymptotically independent.
\end{thm}

Let $\mathrm{TW}$ denote the Tracy--Widom law of index $2$ when $\mathbb F=\mathbb C$ and
of index $1$ when $\mathbb F=\mathbb R$, and define the deterministic centering and
scaling
\begin{equation}
  m_{p,n}^+=\frac{\sigma^2}{n}\bigl(\sqrt{n-1}+\sqrt p\bigr)^2,
  \qquad
  s_{p,n}^+=\frac{\sigma^2}{n}\bigl(\sqrt{n-1}+\sqrt p\bigr)
  \Bigl(\frac1{\sqrt{n-1}}+\frac1{\sqrt p}\Bigr)^{1/3}>0 .
  \label{e:john}
\end{equation}

\begin{thm}[Tracy--Widom fluctuations of the largest non-spiked
eigenvalue \cite{johnstone01,johansson00,baik06,bkyy16}]\label{t:twedge}
Assume $p/n\to c\in(0,\infty)$ and that every spike is supercritical, that is,
$\theta_K>\sigma^2\sqrt c$. Then
\[
  \text{(a)}\ \ \widehat\lambda_{p,K+1}^{\mathrm{id}}\to \sigma^2(1+\sqrt c)^2
  \ \ \text{almost surely},
  \qquad
  \text{(b)}\ \
  \frac{\widehat\lambda_{p,K+1}^{\mathrm{id}}-m_{p,n}^+}{s_{p,n}^+}
  \xrightarrow{\ \mathrm{law}\ }\mathrm{TW}.
\]
The centering $m_{p,n}^+$ and scaling $s_{p,n}^+$ depend only on $(n,p,\sigma^2)$, not on the spikes.
\end{thm}
Part (a) follows from the limiting locations of the sample eigenvalues in the spiked model \cite{baik06}. For part (b), eigenvalue sticking away from the BBP threshold transfers the edge fluctuations of the unspiked Wishart model to the largest non-outlier eigenvalue \cite[Theorem~2.7 and Remark~8.7]{bkyy16}. Since all spikes are fixed and supercritical, this gives the Tracy--Widom limit with the centering and scaling in \eqref{e:john}.

\subsection{Labeling of the eigenvalues and asymptotic regimes}\label{s:lab}

The correspondence between sample and population eigenvalues announced at the end of Section~\ref{ss:setup} is made precise as
follows. By Proposition~\ref{p:pop}, the population spike eigenvalues satisfy
$\lambda_{p,\ell}\xrightarrow{\mathbb P}\sigma^2+\theta_\ell$ for each $\ell\le r+s$, and
the limits $\sigma^2+\theta_1>\cdots>\sigma^2+\theta_{r+s}$ are distinct. Hence, on an
event of probability tending to one, the population spikes are separated, and each of
them is close to exactly one eigenvalue of $\widehat\Sigma_p$. We write
$\widehat\lambda_{p,\ell}$ for the eigenvalue of $\widehat\Sigma_p$ associated in this way
with $\theta_\ell$, and it is the object whose fluctuations we study. 

We study four asymptotic regimes. The main object of the paper is the disproportional
regime $p/n\to0$, in which the ratio $p^2/n$ separates three cases:
\begin{align}
\text{the Haar-dominated regime}\quad & \frac{p^2}{n}\to 0, \\
  \text{the mixed regime}\qquad
  &\frac{p^2}{n}\to\tau\in(0,\infty),\label{e:regm}\\
  \text{the sample-dominated regime}\qquad
  &\frac pn\to0,\quad\frac{p^2}{n}\to\infty.\label{e:infr}
\end{align}
Theorem~\ref{t:main} covers the first two of them at once, that is the whole range
$\tau\in[0,\infty)$, with a single proof. In addition, we also treat the proportional regime
\begin{equation}
  \frac pn\to c\in(0,\infty).\label{e:prop}
\end{equation}
The role of the ratio $p^2/n$, and the meaning of the two scales $p$ and $\sqrt n$, were
explained in Section~\ref{ss:why}.

All the objects needed to state our results are now defined: the model, the two sources of
randomness, the classical results from RMT, and the four regimes. We now present our results.

\section{Main results}
\label{s:main results}

In this section we state the fluctuation results announced in Section~\ref{ss:main} as well as other relevant new results. The disproportional regime $p/n\to0$ is the
main object of the paper, so we treat it first. Theorem~\ref{t:main} covers at once the
Haar-dominated regime $p^2/n\to0$ and the mixed regime $p^2/n\to\tau\in(0,\infty)$, both
at the scale $p$; Theorem~\ref{t:infy} covers the sample-dominated regime
$p^2/n\to\infty$, at the scale $\sqrt n$. We then treat the proportional regime
$p/n\to c>0$ in Theorem~\ref{t:prop}, where the classical properties of the spiked model are recovered. Section~\ref{ss:scal} explains why each regime uses the scaling it uses, and
completes Table~\ref{t:four}. All the statements are given at a fixed spike
$\theta_\ell$, and none of the limiting laws depends on the realization of $O_p$.

\subsection{The $p^2/n \to \tau \in [0,\infty[$ regime theorem}

We start with the regime $p^2/n\to\tau\in[0,\infty)$. It contains both the
Haar-dominated case $\tau=0$ and the mixed case $\tau\in(0,\infty)$, and the natural
scaling in both is $p$. We give the proof of Theorem \ref{t:main} in Section \ref{proof of theorem main}.

\begin{thm}[The $p^2/n \to \tau \in [0,\infty[$ regime theorem]\label{t:main}
Consider the model of Section~\ref{s:mod}. Assume that we are in the regime where
\[
  \frac{p^2}{n} \to\tau\in[0,\infty).
\]
For each fixed $\ell\in\{1,\ldots,r+s\}$,
\begin{equation}
  p\bigl(\widehat\lambda_{p,\ell}-\sigma^2-\theta_\ell\bigr)
  \xrightarrow{\ \mathrm{law}\ }
  W_\ell
  +\sqrt{\kappa\tau V_0(\theta_\ell)}\,G_\ell
  +\tau\frac{\sigma^2(\sigma^2+\theta_\ell)}{\theta_\ell},
  \label{e:main}
\end{equation}
or equivalently, 
\begin{equation}
  p\bigl(\widehat\lambda_{p,\ell}-\phi_{p/n}(\theta_\ell)\bigr)
  \xrightarrow{\ \mathrm{law}\ }
  W_\ell+\sqrt{\kappa\tau V_0(\theta_\ell)}G_\ell.
  \label{e:mainc}
\end{equation}
where $G_\ell\sim \mathcal N (0,1)$ is independent of $W_\ell$, $W_\ell$ is given by \eqref{e:well}, $\kappa$ is given by \eqref{e:kapp} and $\phi_{p/n}(\cdot), V_0(\cdot)$ are given by \eqref{eq:definition of phi_c and V_c}. 
\end{thm}

\begin{rem}[Theorem \ref{t:main} covers both the Haar-regime and the mixed regime.]
    Theorem~\ref{t:main} covers the whole range $\tau\in[0,\infty)$ with a single proof, i.e. it gives the limit both in the Haar-dominated regime and the mixed regime as defined in Section \ref{s:lab}. However, it produces a different limiting law in both regime, which is why we separate those regimes in Table \ref{t:four}. In the Haar regime $\tau=0$, the Gaussian term $\sqrt{\kappa\tau}(\sigma^2+\theta_\ell)G_\ell$
and the deterministic shift $\tau\sigma^2(\sigma^2+\theta_\ell)/\theta_\ell$ both vanish,
and \eqref{e:mainc} reduces to
\[
 p\bigl(\widehat\lambda_{p,\ell}-\phi_{p/n}(\theta_\ell)\bigr)
  \xrightarrow{\ \mathrm{law}\ }
  W_\ell.
\]
The limit is then $W_\ell$ alone, that is, a linear combination of independent chi-square variables with one degree of freedom in the real case,
and exponential variables in the complex case. It carries no Gaussian component.
\end{rem}

\begin{rem}\label{r:randomlimit}
The limit in \eqref{e:main} can be read as
\[\begin{split}
  p\bigl(\widehat\lambda_{p,\ell}-\sigma^2-\theta_\ell\bigr)
  \xrightarrow{\ \mathrm{law}\ }{}&
  \underbrace{W_\ell}_{\text{Haar fluctuation}}
  +\underbrace{\sqrt{\kappa\tau}\,(\sigma^2+\theta_\ell)G_\ell}_{\text{sample fluctuation}}
  +\underbrace{\tau\frac{\sigma^2(\sigma^2+\theta_\ell)}{\theta_\ell}}_{\text{bulk shift}}.
\end{split}\]
The variable $W_\ell$ comes from the Haar matrix. The Gaussian variable $G_\ell$ comes from the sample vectors. These two sources are independent before taking the limit, and this is why they are independent in the limit. The bulk shift is deterministic: it is not a third source of noise but
a bias. It is the same finite-sample bias as in the classical spiked
model: the sample spike does not estimate $\sigma^2+\theta_\ell$ but
$\phi_{p/n}(\theta_\ell)=\sigma^2+\theta_\ell
+\frac{p}{n}\frac{\sigma^2(\sigma^2+\theta_\ell)}{\theta_\ell}$, and
the difference $\phi_{p/n}(\theta_\ell)-\sigma^2-\theta_\ell$ is of
order $p/n$. After scaling by $p$, it becomes
$\frac{p^2}{n}\frac{\sigma^2(\sigma^2+\theta_\ell)}{\theta_\ell}$,
which converges to
$\tau\,\sigma^2(\sigma^2+\theta_\ell)/\theta_\ell$.
\end{rem}

\subsection{The sample-dominated theorem}

The formula of Theorem~\ref{t:main} suggests what happens when $\tau=\infty$. Heuristically, when replacing
$\tau$ by $p^2/n$ in \eqref{e:main}, one expects
$p(\widehat\lambda_{p,\ell}-\phi_{p/n}(\theta_\ell))$ to behave like
\[
  W_\ell
  +\sqrt{\kappa\frac{p^2}{n}}\,(\sigma^2+\theta_\ell)G_\ell
 .
\]
When $p^2/n\to\infty$, this quantity is not tight: the standard deviation of the Gaussian
term grows like $\sqrt{p^2/n}$ while $W_\ell$ remains bounded. The right
normalization is therefore obtained by dividing by $\sqrt{p^2/n}$, which amounts to scaling with $\sqrt n$ instead of $p$. This is what the next theorem does. We give its proof in Section \ref{proof of theorem infinity}.

\begin{thm}[Sample-dominated  regime theorem]\label{t:infy}
Consider the model of Section~\ref{s:mod}. Assume that we are in the sample-dominated regime, where
\[
  \frac pn\to0,
  \qquad
  \frac{p^2}{n} \to\infty.
\]
Then, for each fixed $\ell\in\{1,\ldots,r+s\}$,
\begin{equation}
  \sqrt n
  \left(
    \widehat\lambda_{p,\ell}-\phi_{p/n}(\theta_\ell)
  \right)
  \xrightarrow{\ \mathrm{law}\ }
  \sqrt{\kappa V_0(\theta_\ell)}G_\ell.
  \label{e:infy2}
\end{equation}
where $G_\ell\sim \mathcal N (0,1)$, $\kappa$ is given by \eqref{e:kapp} and $\phi_{p/n}(\cdot), V_0(\cdot)$ are given by \eqref{eq:definition of phi_c and V_c}.  
\end{thm}
As shown later in Proposition~\ref{p:diverge}, the scaling by $p$ has no limit
here, so $\sqrt n$ is the only correct choice. In the sample-dominated regime and with the $\sqrt n$  scaling, only the Gaussian fluctuation remains. The Haar variable $W_\ell$ has
disappeared from the limit, because it is of order $1/p$, which is negligible compared to the sample fluctuations of order $1/\sqrt n$ when $p^2/n\to\infty$. 

Notice that in the sample-dominated regime, the centering is necessarily  $\phi_{p/n}(\theta_\ell)
=\sigma^2+\theta_\ell+\frac pn\frac{\sigma^2(\sigma^2+\theta_\ell)}{\theta_\ell}$ and it can't simply be $\sigma^2+\theta_\ell$ as in Theorem \ref{t:main}. The
deterministic bias $\frac pn \frac{\sigma^2(\sigma^2+\theta_\ell)}{\theta_\ell}$ can't be left in the limit as it is of order $p/n$; after scaling by $\sqrt n$, it diverges since $\sqrt n\cdot p/n=p/\sqrt n\to\infty$. This deterministic bias must therefore be subtracted to be removed from the limit, exactly as in
the classical Bai--Yao theorem.

\subsection{The proportional regime theorem}

We now treat the proportional regime $p/n\to c>0$. Although point (ii) of Theorem \ref{t:prop} is similar to the content of Theorem \ref{t:infy}, the proportional regime differs from the sample dominated in three ways. First, the proof techniques are different and rely on different tools. Second, the CLT presented in Theorem \ref{t:prop} only holds for supercritical spikes $\theta_\ell>\sigma^2\sqrt c$; this condition was automatic in the disproportional case. Third, in the proportional regime the bulk of the spectrum spreads
over the interval $[\sigma^2(1-\sqrt c)^2,\,\sigma^2(1+\sqrt c)^2]$:
the empirical spectral distribution converges to the Marchenko--Pastur
law, and the largest non-spiked eigenvalue has Tracy--Widom
fluctuations at the scale $n^{-2/3}$, as in the classical
deterministic spiked model. The following theorem collects these facts and is proven in Section \ref{proof of theorem prop}.

\begin{thm}[Limits and fluctuations in the proportional regime]\label{t:prop}
Consider the model of Section~\ref{s:mod}. Assume that we are in the proportional regime, where 
$$p/n\to c\in(0,\infty).$$ 
Then:
\begin{enumerate}
\item[\textnormal{(i)}] The empirical spectral distribution of $\widehat\Sigma_p$
converges weakly, almost surely, to the Marchenko--Pastur law $\mu_{c,\sigma^2}$,
whose support has upper edge $\sigma^2(1+\sqrt c)^2$.
\item[\textnormal{(ii)}] For each $\ell \in \{1,\dots,r+s\}$ with $\theta_\ell>\sigma^2\sqrt c$,
\[
  \sqrt n\,\bigl(\widehat\lambda_{p,\ell}-\phi_{p/n}(\theta_\ell)\bigr)
  \xrightarrow{\ \mathrm{law}\ }  \sqrt{\kappa V_c(\theta_\ell)}G_\ell, \qquad \text{and} \qquad \widehat\lambda_{p,\ell}\xrightarrow{\mathbb P}\phi_c(\theta_\ell).
\]
\item[\textnormal{(iii)}] For each $\ell \in \{1,\dots,r+s\}$ with $0<\theta_\ell\leq\sigma^2\sqrt c$,
$\widehat\lambda_{p,\ell}\xrightarrow{\mathbb P}\sigma^2(1+\sqrt c)^2$.
\item[\textnormal{(iv)}] If moreover every spike is supercritical, that is
$\theta_{r+s}>\sigma^2\sqrt c$, then the largest bulk eigenvalue satisfies
\[
  \frac{\widehat\lambda_{p,r+s+1}-m_{p,n}^+}{s_{p,n}^+}
  \xrightarrow{\ \mathrm{law}\ }\mathrm{TW},
  \qquad\text{and}\qquad
  \widehat\lambda_{p,r+s+1}
  =\sigma^2\bigl(1+\sqrt{p/n}\bigr)^2
  +\mathrm O_{\mathbb P}\bigl(n^{-2/3}\bigr).
\]
\end{enumerate}
where $G_\ell\sim \mathcal N (0,1)$, $\kappa$ is given by \eqref{e:kapp}, $\phi_{p/n}(\cdot), V_c(\cdot)$ are given by \eqref{eq:definition of phi_c and V_c} and $m_{p,n}^+$, $s_{p,n}^+$ are given by \eqref{e:john}. 
None of these limiting laws depends on $O_p$.
\end{thm}

The four parts describe the whole spectrum. Part (i) is the Marchenko--Pastur limit of
the bulk; part (ii) gives the Gaussian fluctuations of the supercritical spikes, with the
classical variance $V_c$; part (iii) shows that a subcritical spike does not produce an outlier
in the sample spectrum: the corresponding sample eigenvalue converges
to the upper edge $\sigma^2(1+\sqrt c)^2$ of the bulk, regardless of
the value of $\theta_\ell$; part (iv) locates the largest
non-spiked eigenvalue and gives its Tracy--Widom fluctuations. 

The final sentence stating that
none of these laws depends on $O_p$, is the precise form of the statement that our model behaves like a classical spiked covariance model in the proportional regime: the Haar matrix is drawn once, but the Haar fluctuations are negligible compared to the sample fluctuations. Hence, the Haar matrix effect on the spikes is of lower order and does not enter any of the four limits.

\subsection{Two possible scalings}
\label{ss:scal}

Theorems~\ref{t:main},~\ref{t:infy} and~\ref{t:prop} use different
scalings: the first uses $p$, the second and the third use $\sqrt n$.
The next two propositions show that, outside the mixed regime, these choices are the only ones that give a limit: under the other scaling,
the rescaled fluctuation diverges in probability, so that no
subsequence converges in distribution. In the mixed regime
$p^2/n\to\tau\in(0,\infty)$, both scalings give a limit, and they
differ by the deterministic factor $\sqrt n/p\to\tau^{-1/2}$.

\begin{prop}[The other scaling has no limit]\label{p:diverge}
Fix $\ell\in\{1,\ldots,r+s\}$.
\begin{enumerate}
  \item[\textnormal{(i)}] If $p^2/n\to0$, then
  \[
    \sqrt n\,\bigl|\widehat\lambda_{p,\ell}-\phi_{p/n}(\theta_\ell)\bigr|
    \xrightarrow{\mathbb P}\infty .
  \]
  \item[\textnormal{(ii)}] If $p/n\to0$ and $p^2/n\to\infty$, then
  \[
    p\,\bigl|\widehat\lambda_{p,\ell}-\phi_{p/n}(\theta_\ell)\bigr|
    \xrightarrow{\mathbb P}\infty .
  \]
  \item[\textnormal{(iii)}] If $p/n\to c>0$ and $\theta_\ell>\sigma^2\sqrt c$, then
  \[
    p\,\bigl|\widehat\lambda_{p,\ell}-\phi_{p/n}(\theta_\ell)\bigr|
    \xrightarrow{\mathbb P}\infty .
  \]
\end{enumerate}
In each of the three cases the scaled sequence is not tight, and no subsequence of it
converges in distribution to a real random variable.
\end{prop}

In the mixed regime the scalings by $p$ and $\sqrt n$ are equivalent and both yield a finite limit, since they differ by the bounded
factor $\sqrt n/p\to\tau^{-1/2}$. For completeness we record the limit under the scaling
by $\sqrt n$, which is the one used in the third row of Table~\ref{t:four}.

\begin{prop}[In the mixed regime both scalings work]\label{p:mixed}
Assume $p^2/n\to\tau\in(0,\infty)$ and fix $\ell\in\{1,\ldots,r+s\}$. Then
\[
  \sqrt n\,\bigl(\widehat\lambda_{p,\ell}-\phi_{p/n}(\theta_\ell)\bigr)
  \xrightarrow{\ \mathrm{law}\ }
  \tau^{-1/2}W_\ell+\sqrt{\kappa V_0(\theta_\ell)}\,G_\ell ,
\]
where $G_\ell\sim\mathcal N(0,1)$ is independent of $W_\ell$.
\end{prop}

The proofs of Propositions \ref{p:diverge} and \ref{p:mixed} are given in Section \ref{s:scalprf}.

\begin{rem}[Why we use the scaling by $p$ when $\tau$ is finite]\label{r:conv}
When $\tau\in(0,\infty)$, the two scaled sequences differ by the deterministic factor
$\sqrt n/p$, which converges to $\tau^{-1/2}\in(0,\infty)$, and Proposition~\ref{p:mixed}
shows that both converge. In this regime the choice is only a matter of presentation.
We use $p$ because it is the only one of the two which gives a limit on the whole range
$\tau\in[0,\infty)$, since the scaling by $\sqrt n$ diverges when $\tau=0$, by
Proposition~\ref{p:diverge}.
\end{rem}

Before turning to the application, we illustrate the two theorems on the smallest instance
of the model, where every quantity can be written in closed form.

\subsection{An example: two rank-one perturbations}

This example is the simplest instance of the model, with $r=s=1$, and it can be computed in closed form. We first look at the
population matrix, then at the sample eigenvalues.

\medskip
\noindent\textbf{The population spikes as a function of the angle between eigenvectors.}
We first drop the Haar assumption and compute the eigenvalues of any rank-two perturbation (unrelated to a Haar matrix). Let $a,b\in\mathbb F^p$ be two unit
vectors, and let $t=|a^*b|^2\in[0,1]$ be the squared overlap of the
two directions, so that $t=0$ when $a$ and $b$ are orthogonal, and $t=1$ when $a$ and
$b$ are colinear. The matrix $\alpha\,aa^*+\beta\,bb^*$ with
$\alpha>\beta>0$ has rank at most two, and its two nonzero eigenvalues
are
\begin{equation}
  \mu_\pm(t)=\frac{\alpha+\beta}{2}
  \pm\frac12\sqrt{(\alpha-\beta)^2+4\alpha\beta\,t}.
  \label{e:mupm}
\end{equation}
When $t=0$ the eigenvalues are $\alpha$ and $\beta$; when $t=1$ they
are $\alpha+\beta$ and $0$ because the eigenvectors are colinear. The overlap $t$ therefore measures how far the spikes
are pushed away from $\alpha$ and $\beta$.

\medskip
\noindent\textbf{The Haar model corresponds to $t\to 0$.} In our model $A_p = \diag(\alpha,0,\dots,0)$, $B_p =\diag(\beta,0,\dots,0)$, $a=e_1$ and $b=O_pe_1$,
so the overlap is $t_p=|e_1^*O_pe_1|^2=|o_{11}^{(p)}|^2$. By Proposition~\ref{p:haar},
$\sqrt p\,o_{11}^{(p)}\xrightarrow{\mathrm{law}}g$, a standard Gaussian variable, so
$$t_p =|e_1^*O_pe_1|^2=O_\mathbb P\left(\frac 1p\right) \to 0.$$
Hence $e_1$ and $O_pe_1$ are asymptotically orthogonal. This is why the population spikes converge to $\alpha$ and
$\beta$, and not to a mixture of them. A Taylor expansion of \eqref{e:mupm} at $t=0$ gives
$$\mu_+(t)=\alpha+\frac{\alpha\beta}{\alpha-\beta}t+\mathrm O(t^2), \qquad\mu_-(t)=\beta-\frac{\alpha\beta}{\alpha-\beta}t+\mathrm O(t^2),$$
so, with
$p\,t_p\xrightarrow{\mathrm{law}}|g|^2$ and the
continuous mapping theorem, and since $\lambda_{p,1}-\sigma^2$
and $\lambda_{p,2}-\sigma^2$ are the two nonzero eigenvalues
$\mu_+(t_p)$ and $\mu_-(t_p)$ of $A_p+O_pB_pO_p^*$,
\[
  p(\lambda_{p,1}-\sigma^2-\alpha)\xrightarrow{\ \mathrm{law}\ }
  \frac{\alpha\beta}{\alpha-\beta}|g|^2=W_1,
  \qquad
  p(\lambda_{p,2}-\sigma^2-\beta)\xrightarrow{\ \mathrm{law}\ }
  -\frac{\alpha\beta}{\alpha-\beta}|g|^2=W_2.
\]
The $\mathrm O(t_p^2)$ remainders are of order $1/p^2$, so they
vanish after scaling by $p$.
This recovers Proposition~\ref{p:pop} in the rank-one case, and shows that
$W_1=-W_2$. 

\begin{ex}[Sample eigenvalues in this example]\label{e:rank}
Keeping the notation above and adding the sample
step, Theorem~\ref{t:main} gives, for the largest sample eigenvalue in the mixed regime
$p^2/n\to\tau\in[0,\infty)$,
\begin{equation}
\begin{split}
  p\bigl(\widehat\lambda_{p,1}-\sigma^2-\alpha\bigr)
  \xrightarrow{\ \mathrm{law}\ }{}&
  \frac{\alpha\beta}{\alpha-\beta}|g|^2
  +\sqrt{\kappa\tau}\,(\sigma^2+\alpha)G
  +\tau\frac{\sigma^2(\sigma^2+\alpha)}{\alpha},
\end{split}
  \label{e:rana}
\end{equation}
where $G$ is a standard real Gaussian variable independent of $g$, with $|g|^2\sim\chi_1^2$
in the real case and $|g|^2\sim\operatorname{Exp}(1)$ in the complex case. The limit is a
sum of two independent variables: the Haar variable $\frac{\alpha\beta}{\alpha-\beta}|g|^2$,
which is positive, asymmetric and non-Gaussian, and the sample variable
$\sqrt{\kappa\tau}\,(\sigma^2+\alpha)G$, which is centered Gaussian.  When $\tau=0$, only the term
$W_1=\frac{\alpha\beta}{\alpha-\beta}|g|^2$ survives; it is a
positive random variable since $\alpha > \beta$, so
$p(\widehat\lambda_{p,1}-\sigma^2-\alpha)$ converges to a
distribution supported on $[0,\infty)$, and the largest sample spike
lies asymptotically above $\sigma^2+\alpha$. As $\tau$ grows, the
Gaussian term is centered with standard deviation of order
$\sqrt\tau$, so it eventually dominates $W_1$, which does not
depend on $\tau$; the limiting law then becomes closer and closer to a
Gaussian, and its asymmetry disappears.
\end{ex}

\subsection{Numerical experiments}
\label{ss:num}

The rank-one example above is also the simplest setting in which our results can be checked numerically \footnote{The code to reproduce the figures of this section is available at \url{https://github.com/jeanlouvaleau/Internship-Random-Matrix-for-ML-/tree/main}.}. We keep the parameters of this example, namely $r=s=1$, and we set
$\alpha=10$, $\beta=5$ and $\sigma^2=1$, so that the tracked spike is
$\theta_1=\alpha=10$ and the Haar fluctuation is the explicit variable
\[
  W_1=\frac{\alpha\beta}{\alpha-\beta}\,|g|^2=10\,|g|^2 ,
\]
with $|g|^2\sim\chi^2_1$ in the real case and $|g|^2\sim\operatorname{Exp}(1)$ in the
complex case. For each pair $(p,n)$ we do $2000$ independent simulations of the model:
we sample a Haar matrix $O_p$, then $n$ Gaussian vectors with covariance $\Sigma_p$ conditionally on
$O_p$, and we compute the largest eigenvalue $\widehat\lambda_{p,1}$ of the sample
covariance matrix.

Let us emphasize that our results describe limits along sequences $(p,n)$, while a simulation is run at a single
pair $(p,n)$, at which $p^2/n$ is a finite number. We therefore do not compare the
empirical fluctuation to a limiting law, but to the law obtained from Theorem~\ref{t:main}
by replacing $\tau$ by $p^2/n$ and $V_c$ by $V_{p/n}$,
\begin{equation}
  L_{p,n}
  :=W_\ell+\sqrt{\kappa\,\frac{p^2}{n}\,V_{p/n}(\theta_\ell)}\;G_\ell ,
  \label{e:reflaw}
\end{equation}
where $G_\ell\sim\mathcal N(0,1)$ is independent of $W_\ell$. This is the reference at the
scale $p$; at the scale $\sqrt n$ the two statistics differ by the factor $p/\sqrt n$, so
the reference is
\begin{equation}
  \frac{\sqrt n}{p}\,L_{p,n}
  =\sqrt{\frac{n}{p^2}}\;W_\ell+\sqrt{\kappa V_{p/n}(\theta_\ell)}\;G_\ell .
  \label{e:reflawsqrt}
\end{equation}
Formula \eqref{e:reflaw} is a single expression covering the four regimes, so the same
reference can be used in all the configurations we simulate. Along a sequence with
$p^2/n\to\tau\in[0,\infty)$ it converges to the limit of Theorem~\ref{t:main}. Along a
sequence with $p^2/n\to\infty$ it does not converge, but \eqref{e:reflawsqrt} does: the
coefficient $\sqrt{n/p^2}$ tends to zero and only the Gaussian
$\sqrt{\kappa V_c(\theta_\ell)}\,G_\ell$ of Theorems~\ref{t:infy} and~\ref{t:prop}(ii)
remains in the limit law.

We report two kinds of experiments: histograms of the fluctuation against the density of
the reference law, and a measure of the distance between the two as $p^2/n$ varies.

\begin{figure}[H]
    \centering
    \includegraphics[width=1\linewidth]{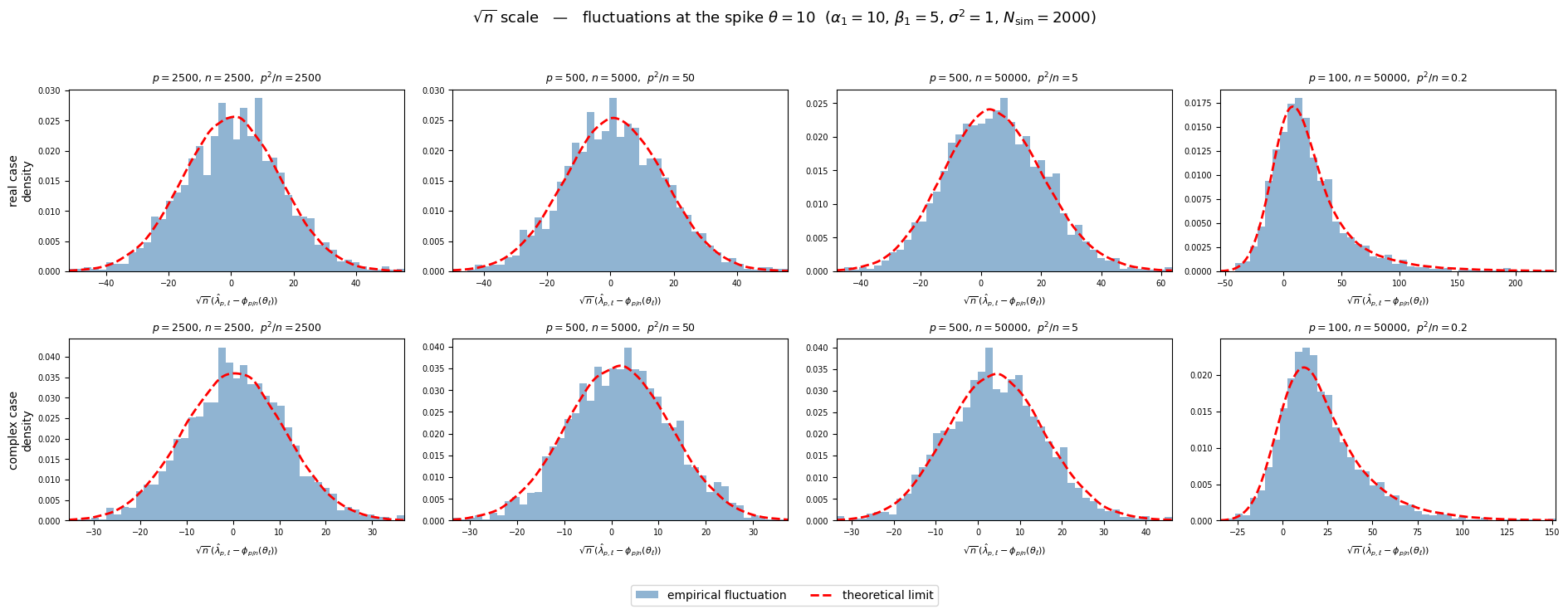}
    \caption{Fluctuations of $\sqrt n\,(\widehat\lambda_{p,1}-\phi_{p/n}(\theta_1))$ at the spike
$\theta_1=10$, over $2000$ simulations of the model, with $r=s=1$, $\alpha=10$, $\beta=5$
and $\sigma^2=1$. Top row: real case. Bottom row: complex case. The four columns are the
four pairs $(p,n)$, with $p^2/n$ decreasing from left to right. The red dashed curve is the
reference law \eqref{e:reflawsqrt}, that is
$\sqrt{n/p^2}\,W_1+\sqrt{\kappa V_{p/n}(\theta_1)}\,G_1$.}
    \label{fig : sqrt n}
\end{figure}

\begin{figure}[H]
    \centering
    \includegraphics[width=1\linewidth]{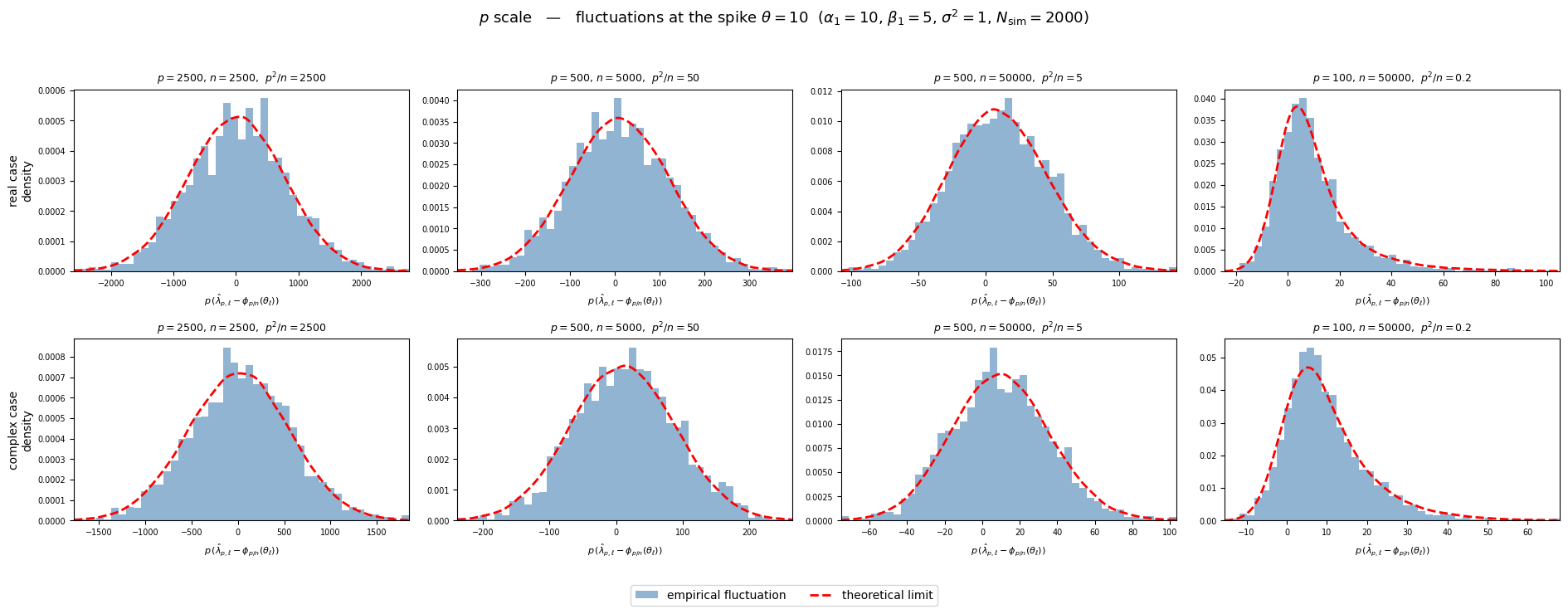}
    \caption{Fluctuations of $p\,(\widehat\lambda_{p,1}-\phi_{p/n}(\theta_1))$, for the same
simulations and the same four pairs $(p,n)$ as in Figure~\ref{fig : sqrt n}. The red dashed
curve is the reference law \eqref{e:reflaw}, that is
$W_1+\sqrt{\kappa (p^2/n) V_{p/n}(\theta_1)}\,G_1$. 
Note the different horizontal scales across columns and between the two rows.}
    \label{fig : sqrt p}
\end{figure}

\medskip
\noindent\textbf{Histograms.}
Figures~\ref{fig : sqrt n} and~\ref{fig : sqrt p} plot the fluctuation at the scale $\sqrt n$ and at
the scale $p$, over $2000$ simulations, in the real and in the complex case, for four configurations of pairs $(p,n)$ chosen so that
$p^2/n$ decreases over four orders of magnitude:
\[
  (2500,2500),\qquad (500,5000),\qquad (500,50000),\qquad (100,50000),
\]
for which $p^2/n$ equals $2500$, $50$, $5$ and $0.2$. The first pair is in the proportional
regime $c=1$, the second and the third have $p/n\to0$, and the last one is the closest to
the Haar-dominated regime that we can reach while keeping $p$ large enough.

In the eight panels the histogram follows the reference law closely, in
both fields and over the whole range of $p^2/n$, which is what the theorems predict.

The two figures give complementary pictures. At the scale $\sqrt n$, the Gaussian part is
essentially fixed across the four configurations, since $V_{p/n}$ varies little, and the
Haar part grows as $p^2/n$ decreases: the histogram widens and becomes skewed. The
skewness is visible only in the last column, where the coefficient $\sqrt{n/p^2}$ in front
of $W_1$ is $2.24$; it is stronger in the real case, where $|g|^2$ is $\chi^2_1$, than in
the complex case, where it is $\operatorname{Exp}(1)$. At the scale $p$, it is
the Haar part $W_1$ that stays fixed and the Gaussian part that grows as $p^2/n$ increases. In the first two columns the law is then almost Gaussian, since $W_1$ is
negligible next to it; the asymmetry only appears in the last column.

Only the second of these two figures shows the divergence of the fluctuations when the scaling is wrong. At the scale $p$, the Gaussian term
of \eqref{e:reflaw} has standard deviation $\sqrt{\kappa(p^2/n)V_{p/n}(\theta_\ell)}$, which
is approximately $774$, $110$, and $7$ in the real case
for the first, second, and last configurations, respectively: this is why the fluctuations take huge values at scale $p$ when $p^2/n \to \infty$. These standard deviations grow to $\infty$ as $p^2/n \to \infty$, as stated in Proposition~\ref{p:diverge}(ii) and~(iii). At the scale $\sqrt n$, the coefficient in front of $W_\ell$ is $\sqrt{n/p^2}$, which
reaches only $2.24$ in the last configuration, so the fluctuations do not take huge values and we don't observe a divergence. To observe the divergence of the fluctuations with a scaling in $\sqrt n$, one would need $p^2/n$ smaller by several orders of magnitude, hence a value of $p$ which is too small for the fluctuations to converge to their limit.

\begin{figure}[H]
    \centering
    \includegraphics[width=1\linewidth]{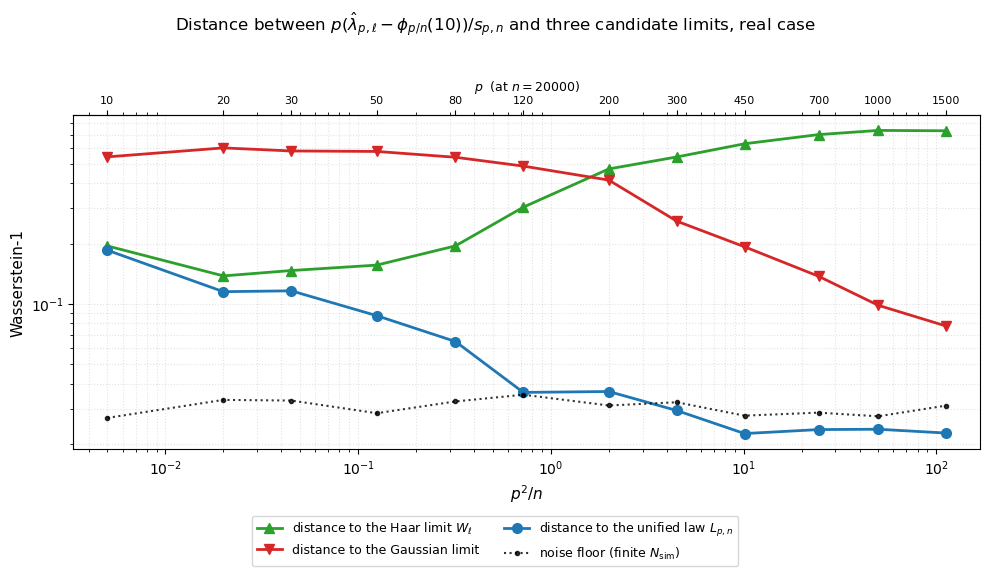}
\caption{Wasserstein-1 distance between the fluctuation
$p\,(\widehat\lambda_{p,1}-\phi_{p/n}(\theta_1))/s_{p,n}$ and three candidate limits, in the
real case, with $n=20000$ fixed and $p$ from $10$ to $1500$, so that $p^2/n$ ranges from
$0.005$ to $112.5$. Both the fluctuation and each candidate are divided by
$s_{p,n}$, the standard deviation of $L_{p,n}$, so that the distances are comparable across
the grid. The dotted line is the Monte-Carlo noise level, computed from two independent
samples of $L_{p,n}$ at the same sample sizes. The upper axis gives the corresponding
values of $p$.}
\label{f:speed}
\end{figure}

\medskip
\noindent\textbf{Distance of the fluctuation to the limit.}
Figure~\ref{f:speed} measures how close the fluctuation is to its limit, on a grid of values
of $p$ at fixed $n = 20 000$, so that $p^2/n$ covers the different regimes. We take $p$ from $10$ to $1500$, so that $p^2/n$ ranges from $0.005$ to
$112.5$, and we work in the real case. To do this, we plot the Wasserstein-1 distance between the normalised
\footnote{The Wasserstein-1 distance has the same unit as the variables being compared, and the
standard deviation of the reference law $L_{p,n}$ ranges from $14$ to $166$ over our grid.
A raw distance would therefore grow from left to right for reasons of scale alone. We
divide both the fluctuation and the reference law by that standard deviation, so that a
given distance means the same thing at every point of the grid.} fluctuation and three limiting laws: $W_\ell$, the Gaussian, and $L_{p,n}$. These distances are computed over $2000$ samples.

In the plots below, the dotted line is the Monte-Carlo noise level: even when the reference law is exactly
right, the distance is not zero, because it is computed between two finite samples,
whose empirical distributions differ by pure sampling error. We measure this non-zero noise level by drawing
two independent samples of $L_{p,n}$, of the same sizes as the ones used for the curves, and
computing the same distance between them; the operation is repeated and averaged, and
it is done at every point of the grid, since the level depends on $(p,n)$ through the shape
of $L_{p,n}$. When a curve reaches this level, the fluctuation and the candidate law are as close to each
other as two samples of the same law would be, so the simulation cannot tell them apart. A curve staying above it reports a discrepancy that
sampling error alone does not explain.

The figure separates the three laws clearly. The distance to $L_{p,n}$
decreases and reaches the noise level from $p^2/n$ of order $2$: past that point the simulation cannot distinguish the
fluctuation from the reference law. The distance to $W_\ell$ alone is small only on the
left and grows steadily, while the distance to the Gaussian alone does the opposite. The
two curves cross around $p^2/n=2$, which is the region where neither of the two terms can
be dropped. Two things are worth noting. On the left of the grid, the distance to
$L_{p,n}$ does not reach the noise level: there, $p$ is between $10$ and $80$, which is too
small for Proposition~\ref{p:pop} to hold. On the right, the distance to the Gaussian alone decreases but is
still above the noise level at $p^2/n=112.5$: the Haar term is of relative size
$\sqrt{n/p^2}$, which is still about $0.09$, so a much larger $p^2/n$ would be needed
before the Gaussian alone becomes accurate.

\subsection{Application and a spectral test}
\label{ss:test}

In Section~\ref{ss:lmm} we introduced the linear mixed model with a continuous batch
effect and explained why the Haar assumption makes the separation of the signal and the
batch possible. Here we develop one application in detail, industrial quality control, we
state precisely when the model is reasonable, and we build a test for one variance
component. We only present the test in the regime $p^2/n\to\tau\in[0,\infty)$, which is
where our test is new; when $p/n \to c>0$ our model behaves similarly to a classical spiked covariance model and some tests already exist: we discuss these at the end of the section.

\subsubsection{An application: industrial quality control}
\label{sss:qc}

A production line in a factory makes items one after another, and after each item is produced, we measure $p$ physical characteristics of it (for instance a length, a mass, an
electrical resistance...). Each item gives one measurement vector $x_{p,k}\in\mathbb R^p$,
whose $\ell$-th coordinate is the value of characteristic $\ell$ for item $k$, and over a
production session we collect $n$ such vectors. A modern line produces and measures a very
large number of items, while the number $p$ of measured characteristics is moderate, of
order at most few hundreds. So $n$ is large and $p$ is moderate, hence $p/n\to0$ and
$p^2/n\to\tau<\infty$.

The three parts of the model $A_p +O_pB_pO_p^*+\sigma^2 I_p$ have a concrete interpretation. The signal covariance $A_p$ defines
the genuine product variation, that is the real differences in terms of quality between items produced in the factory; there are $r$ main modes of variation, each a fixed direction $a_i\in\mathbb R^p$
with strength $\alpha_i$. The batch covariance $O_pB_pO_p^*$ defines the continuous batch effect: when an item is measured, the instrument adds a distortion that moves several
characteristics together along $s$ directions $O_pe_1,\ldots,O_pe_s$, with strengths
$\beta_1,\ldots,\beta_s$. Each item passes through the station on its own and is measured
once, so its batch effect is drawn independently from the others, which is the continuous
batch effect. The noise $\sigma^2I_p$ is the small independent error made on each single
characteristic.

The Haar assumption is reasonable here. The batch effect directions $O_pe_j$ come from the
internal state of the measurement device, such as a thermal drift of one of its
parts, a small misalignment of the item on the table, or a change of focus of a camera.
We have no reason to expect the distortion to act along one combination of characteristics
rather than another, so treating all directions as equally likely is the natural choice.
It is also reasonable to assume the batch effect directions $O_pe_j$ to be unrelated to the product-variation
directions $a_i$: the first come from the measurement device and the second from the
production process, so the two have no reason to be aligned. The assumption fails when the
measurement is done in a way that follows the quality of the product itself, so that the distortion and
the genuine variation become correlated; in that case the spikes of $\Sigma_p$ are no
longer close to the union of the $\alpha_i$ and the $\beta_j$ but a combination of both depending on the correlation between product quality and batch effect.

One consequence of the model matters for the test. By Theorem~\ref{t:main}, the $r+s$
largest sample eigenvalues let us recover the sorted set of strengths
$\{\alpha_i\}\cup\{\beta_j\}$, but the model does not say which strength is a product mode
and which is a batch effect: differentiating the two requires external information. We come back
to this after stating the test.

\subsubsection{The hypothesis we test}
\label{sss:hyp}

Fix one spike $\theta_\ell$ among the sorted signal and batch strengths
$\{\alpha_i\}_{i=1}^r\cup\{\beta_j\}_{j=1}^s$. In practice, one wants to know whether this
spike is strong enough to be considered significant, that is above a threshold $\gamma>0$
chosen by the user. We test
\[
  H_0:\theta_\ell=\gamma
  \qquad\text{against}\qquad
  H_1:\theta_\ell>\gamma .
\]
We calibrate the test at the boundary $\theta_\ell=\gamma$.
For every fixed parameter configuration satisfying the assumptions
with $\theta_\ell<\gamma$, the same argument gives
$T_\ell\to-\infty$ in probability and a bounded critical value,
so the rejection probability tends to zero.
These are pointwise statements; uniform size over parameter
sequences approaching the boundary is not claimed.

\begin{rem}\label{rem:signalspike}
The limit law $W_\ell$ of Proposition~\ref{p:pop} is not the same when $\theta_\ell$ is a
signal spike or a batch spike, so we must fix which of the two we test. In what follows, we
assume for simplicity that $\theta_\ell$ is a signal spike, that is
$\theta_\ell=\alpha_i$ for some $i$. The construction is the same when $\theta_\ell$ is a
batch spike, the only thing that differs is that $W_\ell$ has to be replaced by the corresponding sum over the signal strengths (instead of batch strengths).
\end{rem}

\subsubsection{The test in the regime $p^2/n\to\tau\in[0,\infty)$}
\label{sss:test-tau}

We recall the shape of the fluctuations, since this is what defines the test statistic. 
Write $c=p/n$ and recall the functions defined in \eqref{eq:definition of phi_c and V_c}, $\phi_c(y)=(\sigma^2+y)(1+c\sigma^2/y)$ and $V_c(y)=(\sigma^2+y)^2(1-c\sigma^4/y^2)$.
By Theorem~\ref{t:main}, under $H_0$ and when
$p^2/n\to\tau\in[0,\infty)$,
\begin{equation}
  p\bigl(\widehat\lambda_{p,\ell}-\phi_{p/n}(\gamma)\bigr)
  \xrightarrow{\ \mathrm{law}\ }
  L(\gamma):=
  \sum_{j=1}^s\frac{\gamma\beta_j}{\gamma-\beta_j}Z_j
  +\sqrt{\kappa\tau}\,(\sigma^2+\gamma)G,
  \label{eq:test-law}
\end{equation}
where the $Z_j$ are independent $\chi^2_1$ variables in the real case
or independent $\mathrm{Exp}(1)$ variables in the complex case,
and $G\sim\mathcal N(0,1)$ is independent of them.

To calibrate the test, we must estimate the noise level $\sigma^2$,
the number of spikes, and the batch strengths $\beta_1,\ldots,\beta_s$.
We first estimate the noise level without knowing the number of spikes,
then count the spikes and refine the noise estimate using the bulk.
The centered spike eigenvalues are used to estimate the batch strengths.
These estimates determine the null centering $\phi_{p/n}(\gamma)$
and the reference law. The Gaussian coefficient uses the observed
ratio $p^2/n$ directly, so no separate estimate of its limit $\tau$
is needed. The consistency of these estimates is proved in
Section~\ref{proofs of the properties of the test}.

Once $\hat\sigma^2$, the $\hat\beta_j$ and $\widehat{r+s}$ are computed, the estimated
reference law
\[
  \widehat L=\sum_{j=1}^{\widehat{r+s}-r}\frac{\gamma\,\hat\beta_j}{\gamma-\hat\beta_j}\,Z_j
  +\sqrt{\kappa\tfrac{p^2}{n}}\,(\hat\sigma^2+\gamma)\,G
\]
is an explicit function of the standard variables $Z_1,\ldots,Z_{\hat s},G$, where the
$Z_j$ are i.i.d.\ $\chi^2_1$ if $\mathbb F=\mathbb R$ (resp.\ $\operatorname{Exp}(1)$ if
$\mathbb F=\mathbb C$) and $G\sim\mathcal N(0,1)$ is independent of them. All of these can be
simulated by a standard random number generator, so we can draw from $\widehat L$ and
approximate its $(1-a)$-quantile by the empirical quantile of the draws. Concretely, we
draw $M$ independent copies $c_1,\ldots,c_M$ of $\widehat L$ and sort these copies in increasing order $c_{(1)}\le\cdots\le c_{(M)}$. The empirical $(1-a)$-quantile is then
\[
  \hat q_{1-a}=c_{(\lceil M(1-a)\rceil)},
\]
We then reject $H_0$ when
$T_\ell>\hat q_{1-a}$. As $M\to\infty$, $\hat q_{1-a}$ converges to the $(1-a)$-quantile of
$\widehat L$ and we show in Section~\ref{proofs of the properties of the test} that the latter converges to the quantile of the true reference law $L(\gamma)$.

The reference law $W_\ell$ is a sum over the batch strengths $\beta_1,\ldots,\beta_s$, so to
sample from it we must know which of the $\widehat{r+s}$ estimated spikes are signals and
which are batch effects. We encode this by the ranks of the signal strengths: for each
$i\in\{1,\ldots,r\}$, let $k_i$ be the position of $\alpha_i$ in the list of the $r+s$
spikes $\{\alpha_i\}_{i=1}^r \cup \{\beta_j\}_{j=1}^s$ sorted in decreasing order, so
that $\theta_{k_i}=\alpha_i$. Knowing $r$ and $\{k_1,\ldots,k_r\}$ is
what lets us remove the $r$ signal-related centered eigenvalues and keep the $\hat s$
remaining ones as estimates $\hat\beta_1,\ldots,\hat\beta_{\hat s}$ of the batch strengths.

For completeness, stop and do not reject if
$\widehat{r+s}\notin\{r+1,\ldots,p-1\}$, if a required signal rank
is unavailable or duplicated, or if any retained batch estimate
belongs to $\{0,\gamma\}$.
The proof shows that these exceptional events have probability
tending to zero, so this convention does not affect the asymptotic
conclusions.

\begin{algorithm}[H]
\caption{Spectral test in the regime $p^2/n\to\tau\in[0,\infty)$}
\label{alg:spectral_test}
\begin{algorithmic}[1]
\State \textbf{Input:} eigenvalues $\widehat\lambda_{p,1}\ge\cdots\ge\widehat\lambda_{p,p}$; threshold $\gamma$; tested rank $\ell$; level $a$; margin $\epsilon_p$; number of draws $M$; number of signals $r$ and signal ranks $\{k_1,\dots,k_r\}$.
\State $\tilde\sigma^2\gets\frac1p\sum_{k=1}^p\widehat\lambda_{p,k}$ \Comment{preliminary noise estimate}
\State $\widehat{r+s}\gets\#\{k:\widehat\lambda_{p,k}>\tilde\sigma^2(1+\sqrt{p/n})^2+\epsilon_p\}$ \Comment{number of spikes}
\State $\hat\sigma^2\gets\frac1{p-\widehat{r+s}}\sum_{k>\widehat{r+s}}\widehat\lambda_{p,k}$ \Comment{improved noise estimate}
\State remove the $r$ signal ranks; let $\hat\beta_1,\dots,\hat\beta_{\hat s}$ be the remaining centered eigenvalues, $\hat s=\widehat{r+s}-r$
\State $T_\ell\gets p\left(\widehat\lambda_{p,\ell}-(\hat\sigma^2+\gamma)(1+\tfrac pn\tfrac{\hat\sigma^2}{\gamma})\right)$
\State draw $c_1,\dots,c_M$ i.i.d.\ from $\widehat L=\sum_{j=1}^{\hat s}\tfrac{\gamma\hat\beta_j}{\gamma-\hat\beta_j}Z_j+\sqrt{\kappa\tfrac{p^2}{n}}(\hat\sigma^2+\gamma)G$
\State $\hat q_{1-a}\gets c_{(\lceil M(1-a)\rceil)}$ \Comment{empirical $(1-a)$-quantile}
\State \textbf{reject} $H_0$ iff $T_\ell>\hat q_{1-a}$
\end{algorithmic}
\end{algorithm}

Notice that the same procedure
covers both $\tau=0$ and $\tau>0$. When $\tau>0$ (the mixed regime), the reference law has both a Haar part and a
Gaussian part. When $\tau=0$ (the Haar-dominated regime), the plugged-in Gaussian scale
$\sqrt{\kappa p^2/n}(\hat\sigma^2+\gamma)$ tends to zero, so $\widehat L$ reduces to the
weighted sum $\sum_j\frac{\gamma\hat\beta_j}{\gamma-\hat\beta_j}Z_j$ alone, and
the test is calibrated against $W_\ell(\gamma)$ with no Gaussian part. 

We collect the assumptions under which the test has the stated
properties.

\begin{Assumption}\label{ass:common}
\begin{itemize}
\item[(A1)] The data $x_{p,1},\ldots,x_{p,n}$ follow the model of Section~\ref{s:mod}, that
is $\Sigma_p=A_p+O_pB_pO_p^*+\sigma^2I_p$ with $O_p$ Haar, $A_p$, $B_p$ diagonal as in
Remark~\ref{r:wlog}.
\item[(A2)] The integers $r,s\ge1$ are fixed, and the $\alpha_i$ and $\beta_j$ are distinct
and nonzero. Every spike other than the tested one is different from $\gamma$.
\item[(A3)] The margin satisfies $\epsilon_p>0$, $\epsilon_p\to0$, and $\sqrt{p/n}=\mathrm o(\epsilon_p)$.
\item[(A4)] The dimensions grow with $p^2/n\to\tau\in[0,\infty)$, so in particular
$p/n\to0$.
\item[(A5)] The tested spike $\theta_\ell$ is a signal spike.\footnote{Assumed for
simplicity; the test also holds when $\theta_\ell$ is a batch spike, on replacing (A6) and
the law $W_\ell$ accordingly.}
\item[(A6)] The number of signals $r$ and the ranks $\{k_1,\ldots,k_r\}$ of
$\alpha_1,\ldots,\alpha_r$ in the sorted list of the $r+s$ spikes are known.
\item[(A7)] The number of Monte-Carlo draws is an integer sequence $M=M_p\to\infty$. Conditionally on the data, the draws are independent with law $\widehat L$, generated using auxiliary randomness independent of the data.
\end{itemize}
\end{Assumption}

\begin{thm}[Level and consistency]\label{thm:test}
Under assumptions (A1)--(A7), in the regime $p^2/n\to\tau\in[0,\infty)$, the
test in Algorithm~\ref{alg:spectral_test} is asymptotically of level $a$ and consistent: under $H_0$,
$\mathbb P(\text{reject})\to a$; under $H_1$, $\mathbb P(\text{reject})\to1$.
\end{thm}

\noindent The proof is given in Section~\ref{subsubsection: Proof of consistency and level of the test}.

\subsubsection{The regime $p^2/n\to\infty$, and a comparison}
\label{sss:compare}
In the sample-dominated regime of Theorem~\ref{t:infy}, or in the
proportional regime of Theorem~\ref{t:prop}, the Haar contribution is
negligible at the sampling fluctuation scale.
Under $H_0$, the limiting reference law is
$\mathcal N(0,\kappa V_c(\gamma))$, provided that
$\gamma>\sigma^2\sqrt c$ when $c>0$.
It depends only on $\gamma$, $\sigma^2$ and $c=\lim p/n$.
With a noise estimator satisfying
$\sqrt n(\hat\sigma^2-\sigma^2)\xrightarrow{\mathbb P}0$,
the statistic
$\sqrt n\bigl(\widehat\lambda_{p,\ell}
-(\hat\sigma^2+\gamma)(1+\tfrac pn\tfrac{\hat\sigma^2}{\gamma})\bigr)$
can be compared with a Gaussian quantile of estimated variance
$\kappa\hat V(\gamma)$, where
$\hat V(\gamma)=(\hat\sigma^2+\gamma)^2
(1-\tfrac pn\tfrac{(\hat\sigma^2)^2}{\gamma^2})$.
No Monte-Carlo step is needed for this Gaussian reference law.
Related spectral procedures have been developed for the problem of
detecting or counting signals, in the proportional subcase $p/n\to c>0$:
\cite{nadakuditi2008sample} estimate the number of signals from the sample eigenvalues and
describe the detection threshold $\sigma^2\sqrt{p/n}$; \cite{kritchman2009nonparametric}
build a sequence of hypothesis tests in which each eigenvalue is compared to a
Tracy--Widom quantile; and \cite{bianchi2011performance} compute the level and the power of
a generalized likelihood ratio test for single-source detection. 

It is interesting to compare the difficulty of the two regimes. In the classical spiked model setting without Haar fluctuation, the proportional regime $p/n\to c>0$ is the most difficult setup: weak non-supercritical
spikes can't be detected, and the sample covariance is not consistent in
operator norm. The regime $p/n\to0$ is easier in the classical spiked model because the sample
covariance is consistent and every fixed spike is detectable.

In our model, calibration in the regime $p^2/n\to\tau<\infty$
requires information about the batch strengths because the Haar
contribution is not negligible. After removing the deterministic
sampling bias, the Haar fluctuation dominates when $\tau=0$,
whereas the Haar and sampling fluctuations have the same order
when $\tau>0$. The reference law therefore requires estimates of
the batch strengths and information identifying the batch spikes.
In the sample-dominated and proportional regimes considered above,
the Gaussian reference law for a separated sample outlier does not
involve these batch parameters.

\subsubsection{Relaxing assumption (A6)}
\label{sss:relax assumption A6}

Assumption (A6) is the strongest of the list: it requires external knowledge of the signal
ranks, that is, to know which of the $r+s$ spikes are signals and which are batch effects. We
present three ways to weaken this assumption. The first is a special case where (A6) holds for free; the
second replaces it by prior knowledge of the scale of the signals. Throughout we write, as in
Algorithm~\ref{alg:spectral_test}, $\hat\mu_\ell=\widehat\lambda_{p,\ell}-\hat\sigma^2$ for
the centered eigenvalue at rank $\ell$, $\hat\phi(\gamma)=(\hat\sigma^2+\gamma)(1+\tfrac
pn\tfrac{\hat\sigma^2}{\gamma})$ for the estimated centering,
$\hat v=\sqrt{\kappa\,p^2/n}\,(\hat\sigma^2+\gamma)$ for the estimated Gaussian scale, and
$v_0=\sqrt{\kappa\tau\,V_0(\gamma)}=\sqrt{\kappa\tau}\,(\sigma^2+\gamma)$ for its limit.

\paragraph{A single signal dominating the batch effects.}
Suppose $r=1$ and $\alpha_1>\beta_j$ for every $j\in\{1,\dots,s\}$, that is, one signal
stronger than every batch effect. Then $\alpha_1$ is the largest spike, so its rank is
$k_1=1$: the tested spike is $\theta_1=\alpha_1$, and the batch effects are exactly the
remaining spikes $\theta_2,\dots,\theta_{1+s}$. The signal ranks are therefore known from
the ordering of the spikes alone, with no external information, and (A6) holds
automatically. In practice this covers any setting where a single real effect clearly
dominates every batch effect. The test itself is unchanged: the reference law still
depends on the batch strengths $\beta_1,\dots,\beta_s$ through
$W_1(\gamma)=\sum_{j=1}^s\tfrac{\gamma\beta_j}{\gamma-\beta_j}|g_{1j}|^2$, which are
estimated from the remaining spikes, so the Monte-Carlo calibration is still needed. 

\paragraph{Relaxing the ranks with prior knowledge of the scale of the signals.}
Suppose now that for each signal we have a prior approximate value $\alpha_i^\circ$ of its
strength $\alpha_i$, available for instance from an earlier study or from the physical scale
of the effect. Let $\mathcal S:=\{\theta_1,\dots,\theta_{r+s}\}$ be the set of spikes, and
recall that $k_i$ denotes the rank of $\alpha_i$, so that $\theta_{k_i}=\alpha_i$ and
$\mathcal K_{\mathrm{sig}}:=\{k_1,\dots,k_r\}$ is the set of signal ranks. The idea is to
assign to each prior value $\alpha_i^\circ$ the closest centered eigenvalue.

\begin{prop}[Calibrated signal values]
\label{prop:calibrated}
Assume (A1)--(A5), that the number of signals $r$ is known (the first half of (A6)), (A7),
and:
\begin{itemize}
\item[(C1)] for each $i\in\{1,\dots,r\}$, the prior value $\alpha_i^\circ$ has $\alpha_i$ as
its unique closest spike: $|\alpha_i^\circ-\alpha_i|<|\alpha_i^\circ-\theta|$ for every spike
$\theta\in\mathcal S$ with $\theta\ne\alpha_i$.
\end{itemize}
Estimate each signal rank by
$\hat k_i:=\arg\min_{\ell\in\{1,\dots,\widehat{r+s}\}}|\hat\mu_\ell-\alpha_i^\circ|$. Then
$\mathbb P\bigl(\{\hat k_1,\dots,\hat k_r\}=\mathcal K_{\mathrm{sig}}\bigr)\to1$, and the
rank-$\ell$ test of Algorithm~\ref{alg:spectral_test} run with these estimated ranks is
asymptotically of level $a$ under $H_0$ and consistent under $H_1$.
\end{prop}

The separation condition (C1) is what makes the nearest-spike assignment reliable; it holds,
for example, whenever the prior values are closer to the truth than half the smallest gap
between two spikes. A cleaner statement is available in the common case where the signals
are known to lie above the batch effects on a fixed scale: then a single threshold separates
the two groups, and one no longer needs a prior value for each signal.

\begin{prop}[Test under a known scale separation]
\label{prop:scale-separation}
Assume (A1)--(A5), that $r$ is known, (A7), and:
\begin{itemize}
\item[(S1)] a threshold $\eta>0$ is known with $\alpha_i>\eta>\beta_j$ for all $i,j$.
\end{itemize}
Estimate the set of signal ranks by
$\widehat{\mathcal K}_{\mathrm{sig}}:=\{\ell\in\{1,\dots,\widehat{r+s}\}:\hat\mu_\ell>\eta\}$,
and run the rank-$\ell$ test for a signal rank $\ell\in\{1,\dots,r\}$ with batch ranks
$\{1,\dots,\widehat{r+s}\}\setminus\widehat{\mathcal K}_{\mathrm{sig}}$. Then
$\mathcal K_{\mathrm{sig}}=\{1,\dots,r\}$,
$\mathbb P(\widehat{\mathcal K}_{\mathrm{sig}}=\mathcal K_{\mathrm{sig}})\to1$, and the test
is asymptotically of level $a$ under $H_0$ and consistent under $H_1$.
\end{prop}

The proofs of both propositions are given in
Section~\ref{subsubsection:Proofs of the properties of the test when relaxing assumption (A6}.

\newpage
\section{Proofs}\label{s:prf}

All the proofs below start from the same common step. Conditionally on the Haar matrix $O_p$, the population covariance $\Sigma_p$ is deterministic, so
the eigenvalues of $\widehat\Sigma_p$ are those of a classical spiked
model. This reduction is stated in Section~\ref{s:cond}. The remaining
difficulty is that the conditional spikes $\lambda_{p,\ell}-\sigma^2$
are not fixed: they depend on $p$ and converge to $\theta_\ell$ only as
$p\to\infty$, while the classical spiked covariance model results are stated for fixed spikes.
Each regime handles this difficulty at its own scale and using specific tools:
\begin{itemize}
    \item In the regime $p^2/n\to\tau$ (Section~\ref{proof of theorem main}), the
scale is $p$, and we expand the sample eigenvalue to second order using
the perturbation Lemma~\ref{l:pert}.
\item In the sample-dominated regime
(Section~\ref{proof of theorem infinity}), the scale is $\sqrt n$, at
which a fixed-order expansion is no longer accurate, and we use instead
a central limit theorem (Proposition~\ref{p:dclt}) relatively similar to the Bai--Yao CLT; but applying specificaly to the case $p/n \to 0$. 
\item In the proportional regime
(Section~\ref{proof of theorem prop}), the classical results of the deterministic spiked model are transferred to our model with non-fixed spikes using Lemma~\ref{l:spikecont}.

\end{itemize}
The results on the two scalings are proved in
Section~\ref{s:scalprf} and the properties of the spectral test are proven in Section~\ref{proofs of the properties of the test}.

\subsection{Conditional reduction to a deterministic spiked model}\label{s:cond}

This reduction is the structural core of the argument. The eigenvectors of the
population covariance are random, but conditionally on the Haar matrix the eigenvalue
distribution of the sample covariance depends only on the population eigenvalues.

\begin{lem}[Conditional diagonalization]\label{l:cond}
Let $\Sigma_p=U_p\Lambda_pU_p^*$ be an $\mathcal F_p$-measurable spectral
decomposition, with
\[
  \Lambda_p=\operatorname{diag}(\lambda_{p,1},\ldots,\lambda_{p,p}).
\]
If $y_{p,k}=U_p^*x_{p,k}$, then conditionally on $\mathcal F_p$, the vectors
$y_{p,1},\ldots,y_{p,n}$ are independent centered Gaussian vectors with covariance
$\Lambda_p$, and
\[
  U_p^*\widehat\Sigma_pU_p
  =\frac1n\sum_{k=1}^n y_{p,k}y_{p,k}^*.
\]
In particular, the conditional eigenvalue distribution of $\widehat\Sigma_p$ depends
only on $\Lambda_p$, that is, only on the population eigenvalues.
\end{lem}

\begin{proof}
After conditioning on $\mathcal F_p$, the matrix $U_p$ is deterministic. Orthogonal or
unitary invariance of the Gaussian distribution implies that $U_p^*x_{p,k}$ is
centered Gaussian with covariance
\[
  U_p^*\Sigma_pU_p=\Lambda_p.
\]
Independence is preserved. The identity for $U_p^*\widehat\Sigma_pU_p$ follows by
direct substitution. Orthogonal or unitary conjugation does not change eigenvalues.
\end{proof}

\begin{cor}[Representation by a standard Gaussian matrix]\label{c:repr}
When discussing the conditional eigenvalue distribution, we may replace
$\widehat\Sigma_p$ by
\begin{equation}
  \widetilde\Sigma_p
  =\Lambda_p^{1/2}\left(\frac1nZ_pZ_p^*\right)\Lambda_p^{1/2},
  \label{e:repr}
\end{equation}
where $Z_p=(z_{qk})_{1\leq q\leq p,\,1\leq k\leq n}$ has independent standard real
Gaussian entries in the real case and independent standard circular complex Gaussian
entries in the complex case. The matrix $Z_p$ is independent of $\mathcal F_p$.
\end{cor}

\begin{proof}
The conditional law depends only on $\Lambda_p$, and the standard Gaussian variables
used in the representation can be chosen independently of $\mathcal F_p$.
\end{proof}

\begin{rem}[Randomness left after conditioning]\label{r:gau}
After conditioning on $\mathcal F_p$, the Haar matrix $O_p$, the covariance matrix
$\Sigma_p$, and all eigenvalues and eigenvectors of $\Sigma_p$ are fixed
deterministic objects. The only randomness left in the conditional model is the
Gaussian randomness used to generate the independent sample vectors
$x_{p,1},\ldots,x_{p,n}$. There is no remaining Haar randomness inside the
conditional law. 
\end{rem}

\subsection{Proof of Theorem \ref{t:main}} \label{proof of theorem main}

The idea is to expand the sample eigenvalue to second order in the sample error matrix,
conditionally on the Haar matrix, and then to add back the Haar fluctuation. The
first-order term gives the Gaussian piece; the second-order term splits into a negligible spike part and the deterministic bulk shift; the third-order remainder vanishes after scaling by $p$; and the Haar fluctuation is added at the end using independence. 

\subsubsection{Second order perturbation lemma}
\label{subsec: second order perturbation lemma}

\begin{lem}[Second-order eigenvalue perturbation]\label{l:pert}
Let $M$ be Hermitian, let $\lambda$ be a simple eigenvalue of $M$, and let $v$ be a
corresponding unit eigenvector. Let $\{v_q\}$ be an orthonormal eigenbasis of $M$,
with $Mv_q=\lambda_qv_q$, and let $\ell$ be the index such that $v_\ell=v$. Set
\[
  \delta=\min_{\mu\in\operatorname{Spec}(M)\setminus\{\lambda\}}|\lambda-\mu|>0.
\]
If $E$ is Hermitian and $\|E\|_{\mathrm{op}}\leq\delta/4$, then the eigenvalue
$\widehat\lambda$ of $M+E$ associated with $\lambda$ satisfies
\begin{equation}
  \widehat\lambda-\lambda
  =v^*Ev
  +\sum_{q\neq\ell}
  \frac{|v_q^*Ev|^2}{\lambda-\lambda_q}
  +R,
  \label{e:pert}
\end{equation}
and there exists a universal constant $C>0$ such that
\begin{equation}
  |R|\leq C\delta^{-2}\|E\|_{\mathrm{op}}^3.
  \label{e:remd}
\end{equation}
\end{lem}

The proof of Lemma \ref{l:pert} is given in Appendix \ref{appendix: second order perturbation lemma}.

\begin{rem}[Why a first-order expansion is not enough]
Conditionally on the Haar matrix, the sample eigenvalue admits an
expansion around the population eigenvalue $\lambda_{p,\ell}$,
\[
  \widehat\lambda_{p,\ell}-\lambda_{p,\ell}
  =\text{(first order)}+\text{(second order)}+\text{(remainder)},
\]
with the three terms of order $n^{-1/2}$, $p/n$ and $(p/n)^{3/2}$. After
multiplication by $p$, the first two terms are of order
$\sqrt{p^2/n}$ and $p^2/n$, hence both of order one when $p^2/n\to\tau$,
while the remainder is of order $(p^2/n)\sqrt{p/n}$ and tends to zero. A
first-order expansion would keep only the first term and discard the
second, but the second term contributes a deterministic bias of order
$\tau$ that must be retained. This is why we expand to second order. In
the proof, the second-order term is a sum over the other eigenvalue
directions, which splits into a finite sum over the remaining spikes,
negligible after scaling, and a sum over the $p-r-s$ bulk directions,
which gives the deterministic shift.
\end{rem}

\subsubsection{Proof of Theorem \ref{t:main}}

We fix one index $\ell\in\{1,\ldots,r+s\}$ throughout the proof. All the perturbation identities below are stated on the good event
$\Omega_p$ of Step~1, where the sample perturbation is small enough for
Lemma~\ref{l:pert} to hold. Since $\mathbb P(\Omega_p)\to1$, what
happens on $\Omega_p^c$ does not affect any limit.

\noindent \textbf{Step 1: the good event.}

Let $K = r+s$ and define
\[
  d_\ell
  =\frac14\min\left\{
      \theta_\ell,
      \min_{1\le q\le K,\,q\neq\ell}|\theta_\ell-\theta_q|
    \right\}>0.
\]
By Proposition~\ref{p:pop}, for each $q\le K$,
\[
  \lambda_{p,q}-\sigma^2\xrightarrow{\mathbb P}\theta_q.
\]
For $q>K$, we have $\lambda_{p,q}=\sigma^2$ exactly.  Hence the event
\[
  \Omega_{p}^{\mathrm{gap}}
  =\left\{
    \min_{q\neq\ell}|\lambda_{p,\ell}-\lambda_{p,q}|\ge 2d_\ell,
    \quad
    \lambda_{p,\ell}-\sigma^2\ge 2d_\ell
  \right\}
\]
satisfies
\begin{equation}
  \mathbb P(\Omega_{p}^{\mathrm{gap}})\to1.
  \label{e:gap}
\end{equation}
The population eigenvalues are also uniformly bounded, because
\[
  \|\Sigma_p\|_{\mathrm{op}}
  \le \sigma^2+\max_i\alpha_i+\max_j\beta_j.
\]

After conditioning on the Haar matrix, let $E_p$ be the sample error matrix in the population eigenbasis, defined in Step~2 below. Corollary~\ref{a:wishartop}, applied conditionally on $\mathcal F_p$ as explained in
Remark~\ref{a:wishartrand}, together with $p/n\to0$, gives
\begin{equation}
  \|E_p\|_{\mathrm{op}}
  =\mathrm O_{\mathbb P}\left(\sqrt{\frac pn}\right)
  =o_{\mathbb P}(1).
  \label{e:enrm}
\end{equation}
Set
\begin{equation}
  \Omega_p
  =\Omega_p^{\mathrm{gap}}\cap
  \{\|E_p\|_{\mathrm{op}}\le d_\ell/2\}.
  \label{e:good}
\end{equation}
Then $\mathbb P(\Omega_p)\to1$.  On $\Omega_p$, the eigenvalue of the sample covariance associated with $\lambda_{p,\ell}$ is unique.  It is the eigenvalue denoted by $\widehat\lambda_{p,\ell}$ in Section~\ref{s:mod}.

All perturbation identities below are first proved on $\Omega_p$.  Since $\mathbb P(\Omega_p^c)\to0$, replacing a random variable by any harmless value on $\Omega_p^c$ does not change its limit in probability or in distribution.

\noindent \textbf{Step 2: the conditional diagonal model.}

By Lemma~\ref{l:cond} and Corollary~\ref{c:repr}, conditionally on $\mathcal F_p=\sigma(O_p)$ we may work with
\[
  \widetilde\Sigma_p
  =\Lambda_p^{1/2}\left(\frac1nZ_pZ_p^*\right)\Lambda_p^{1/2},
  \qquad
  \Lambda_p=\operatorname{diag}(\lambda_{p,1},\ldots,\lambda_{p,p}),
\]
where $Z_p=(z_{qk})$ has independent standard real or complex Gaussian entries and is independent of $\mathcal F_p$.  This is not an extra assumption; it is a representation of the conditional law of the original sample covariance eigenvalues.

Put
\[
  E_p=\widetilde\Sigma_p-\Lambda_p.
\]
If $z_q=(z_{q1},\ldots,z_{qn})$ is the $q$-th row of $Z_p$, then the entries of $E_p$ are
\begin{equation}
  (E_p)_{ab}
  =\sqrt{\lambda_{p,a}\lambda_{p,b}}
  \left(
    \frac1n\sum_{k=1}^n z_{ak}\overline{z_{bk}}-\delta_{ab}
  \right).
  \label{e:entr}
\end{equation}
In the real case, the complex conjugate has no effect.

\noindent \textbf{Step 3: second-order perturbation.}

On $\Omega_p$, Lemma~\ref{l:pert} applies to $M=\Lambda_p$, $E=E_p$, and $\lambda=\lambda_{p,\ell}$.  Therefore
\begin{equation}
  \widehat\lambda_{p,\ell}-\lambda_{p,\ell}
  =(E_p)_{\ell\ell}
  +\sum_{q\neq\ell}
  \frac{|(E_p)_{q\ell}|^2}{\lambda_{p,\ell}-\lambda_{p,q}}
  +R_{p,\ell},
  \label{e:expn}
\end{equation}
where
\begin{equation}
  |R_{p,\ell}|
  \le C d_\ell^{-2}\|E_p\|_{\mathrm{op}}^3
  \qquad\text{on }\Omega_p.
  \label{e:rbd}
\end{equation}
We split the second-order sum as
\[
  \sum_{q\neq\ell}
  =
  \sum_{\substack{1\le q\le K\\q\neq\ell}}
  +
  \sum_{q>K}.
\]
The first part contains only finitely many spike directions.  The second part contains the $p-K$ bulk directions.

\noindent \textbf{Step 4: the first-order term.}

Taking $a=b=\ell$ in \eqref{e:entr} gives
\begin{equation}
  (E_p)_{\ell\ell}
  =\frac{\lambda_{p,\ell}}{n}
  \sum_{k=1}^n\bigl(|z_{\ell k}|^2-1\bigr).
  \label{e:linr}
\end{equation}
Define
\[
  S_{p,\ell}
  =\frac1{\sqrt n}
  \sum_{k=1}^n\bigl(|z_{\ell k}|^2-1\bigr).
\]
The summands have mean zero and variance $\kappa$.  The classical one-dimensional central limit theorem gives
\begin{equation}
  S_{p,\ell}\xrightarrow{\ \mathrm{law}\ } \sqrt\kappa\,G_\ell,
  \qquad G_\ell\sim \mathcal N (0,1).
  \label{e:cltt}
\end{equation}
Since $Z_p$ is independent of $\mathcal F_p$, the random variable $S_{p,\ell}$ is independent of all Haar-measurable quantities.

Moreover,
\begin{equation}
  p(E_p)_{\ell\ell}
  =\frac p{\sqrt n}\lambda_{p,\ell}S_{p,\ell}.
  \label{e:ling}
\end{equation}
The coefficient satisfies
\begin{equation}
  \frac p{\sqrt n}\lambda_{p,\ell}
  \xrightarrow{\mathbb P}
  \sqrt\tau(\sigma^2+\theta_\ell).
  \label{e:lincoef}
\end{equation}

\noindent \textbf{Step 5: the other spike directions are negligible.}

Let
\[
  \Pi_{p,\ell}^{\mathrm{sp}}
  =p\sum_{\substack{1\le q\le K\\q\neq\ell}}
  \frac{|(E_p)_{q\ell}|^2}{\lambda_{p,\ell}-\lambda_{p,q}}.
\]
On $\Omega_p$, all denominators in this finite sum have absolute value at least $2d_\ell$.
Moreover $\Omega_p\subset\Omega_p^{\mathrm{gap}}$, and $\Omega_p^{\mathrm{gap}}$ is
$\mathcal F_p$-measurable, since it is defined by the population eigenvalues only. Hence
\[
  \mathbf 1_{\Omega_p}\bigl|\Pi_{p,\ell}^{\mathrm{sp}}\bigr|
  \le \mathbf 1_{\Omega_p^{\mathrm{gap}}}\,\frac{p}{2d_\ell}
  \sum_{\substack{1\le q\le K\\ q\neq\ell}}|(E_p)_{q\ell}|^2 ,
\]
where we used that the denominators are bounded below by $2d_\ell$ on
$\Omega_p^{\mathrm{gap}}$. Taking conditional expectations, the indicator is now
$\mathcal F_p$-measurable and can be taken out, so
\[
  \mathbb E\left[
    \mathbf 1_{\Omega_p}|\Pi_{p,\ell}^{\mathrm{sp}}|
    \mid \mathcal F_p
  \right]
  \le \frac{p}{2d_\ell}\sum_{\substack{1\le q\le K\\q\neq\ell}}
  \mathbb E\left[ |(E_p)_{q\ell}|^2\mid\mathcal F_p\right]
  \le C\frac pn,
\]
where the last inequality uses
\[
  \mathbb E\left[ |(E_p)_{q\ell}|^2\mid\mathcal F_p\right]
  =\frac{\lambda_{p,q}\lambda_{p,\ell}}{n}
\]
together with the boundedness of the population eigenvalues.

Since $p/n\to0$, Markov's inequality gives
\begin{equation}
  \Pi_{p,\ell}^{\mathrm{sp}}=o_{\mathbb P}(1).
  \label{e:spik}
\end{equation}

\noindent \textbf{Step 6: the bulk directions give the deterministic shift.}

For $q>K$, $\lambda_{p,q}=\sigma^2$.  Put
\[
  L_q=\sum_{k=1}^n z_{qk}\overline{z_{\ell k}}.
\]
Then
\begin{equation}
  (E_p)_{q\ell}
  =\frac{\sqrt{\sigma^2\lambda_{p,\ell}}}{n}L_q.
  \label{e:bulk}
\end{equation}
Let $m_p=p-K$.  Conditionally on $\mathcal F_p$ and on the row $z_\ell$, the variables $L_q$, $q>K$, are independent centered Gaussian variables.  In the real case their conditional variance is $\|z_\ell\|_2^2$; in the complex case the conditional second moment is also $\|z_\ell\|_2^2$.  In both cases
\[
  \mathbb E\left[ |L_q|^2\mid \mathcal F_p,z_\ell\right]
  =\|z_\ell\|_2^2,
\]
and
\[
  \operatorname{Var}\left( |L_q|^2\mid \mathcal F_p,z_\ell\right)
  \le C\|z_\ell\|_2^4.
\]
Therefore the conditional Chebyshev inequality gives
\begin{equation}
  V_p:=\frac{1}{m_p\|z_\ell\|_2^2}
  \sum_{q>K}|L_q|^2
  \xrightarrow{\mathbb P}1.
  \label{e:vbulk}
\end{equation}
Also, by the law of large numbers,
\begin{equation}
  U_p:=\frac{\|z_\ell\|_2^2}{n}
  \xrightarrow{\mathbb P}1.
  \label{e:ubulk}
\end{equation}
On $\Omega_p$, the bulk denominator is $\lambda_{p,\ell}-\sigma^2>0$.  Hence the bulk part satisfies
\begin{align}
  \Pi_{p,\ell}^{\text{bulk}}
  &:=p\sum_{q>K}
  \frac{|(E_p)_{q\ell}|^2}{\lambda_{p,\ell}-\sigma^2}\notag\\
  &=\frac{p m_p}{n}
  \frac{\sigma^2\lambda_{p,\ell}}{\lambda_{p,\ell}-\sigma^2}
  U_pV_p.
  \label{e:bfac}
\end{align}
Since $m_p/p\to1$, $p^2/n\to\tau$, and
\[
  \frac{\sigma^2\lambda_{p,\ell}}{\lambda_{p,\ell}-\sigma^2}
  \xrightarrow{\mathbb P}
  \frac{\sigma^2(\sigma^2+\theta_\ell)}{\theta_\ell},
\]
Slutsky's theorem gives
\begin{equation}
  \Pi_{p,\ell}^{\text{bulk}}
  \xrightarrow{\mathbb P}
  \tau\frac{\sigma^2(\sigma^2+\theta_\ell)}{\theta_\ell}.
  \label{e:bshf}
\end{equation}
This is the deterministic second-order bulk shift.

\noindent \textbf{Step 7: the third-order remainder vanishes.}

From \eqref{e:rbd} and \eqref{e:enrm}, on $\Omega_p$,
\begin{align*}
  p|R_{p,\ell}|
  &\le C d_\ell^{-2}p\|E_p\|_{\mathrm{op}}^3\\
  &=\mathrm O_{\mathbb P}\left(
      p\left(\frac pn\right)^{3/2}
    \right)\\
  &=\mathrm O_{\mathbb P}\left(
      \frac{(p^2/n)^{3/2}}{\sqrt p}
    \right)
  =o_{\mathbb P}(1),
\end{align*}
so
\begin{equation}
  pR_{p,\ell}=o_{\mathbb P}(1).
  \label{e:rvsh}
\end{equation}

\noindent \textbf{Step 8: the sample contribution.}

Multiplying \eqref{e:expn} by $p$ and using Steps~4--7 gives
\begin{equation}
  p\bigl(\widehat\lambda_{p,\ell}-\lambda_{p,\ell}\bigr)
  =
  \frac p{\sqrt n}\lambda_{p,\ell}S_{p,\ell}
  +\tau\frac{\sigma^2(\sigma^2+\theta_\ell)}{\theta_\ell}
  +o_{\mathbb P}(1).
  \label{e:wish}
\end{equation}
More precisely, before replacing the bulk part by its limit, the second term is $\Pi_{p,\ell}^{\text{bulk}}$ from \eqref{e:bfac}; equation \eqref{e:bshf} justifies the replacement.

\noindent \textbf{Step 9: adding the Haar fluctuation.}

Write
\begin{equation}
\begin{split}
  p\bigl(\widehat\lambda_{p,\ell}-\sigma^2-\theta_\ell\bigr)
  ={}&p\bigl(\lambda_{p,\ell}-\sigma^2-\theta_\ell\bigr)\\
  &+p\bigl(\widehat\lambda_{p,\ell}-\lambda_{p,\ell}\bigr).
\end{split}
  \label{e:last}
\end{equation}
Set
\[
  H_{p,\ell}=p\bigl(\lambda_{p,\ell}-\sigma^2-\theta_\ell\bigr),
  \qquad
  a_{p,\ell}=\frac p{\sqrt n}\lambda_{p,\ell}.
\]
The pair $(H_{p,\ell},a_{p,\ell})$ is $\mathcal F_p$-measurable.  Proposition~\ref{p:pop} and \eqref{e:lincoef} imply
\[
  (H_{p,\ell},a_{p,\ell})
  \xrightarrow{\ \mathrm{law}\ }
  \left(W_\ell,\sqrt\tau(\sigma^2+\theta_\ell)\right).
\]
The variable $S_{p,\ell}$ is a function of $Z_p$ only and is independent of $\mathcal F_p$.  Lemma~\ref{l:join} therefore gives the joint convergence
\[
  (H_{p,\ell},a_{p,\ell},S_{p,\ell})
  \xrightarrow{\ \mathrm{law}\ }
  \left(W_\ell,\sqrt\tau(\sigma^2+\theta_\ell),\sqrt\kappa G_\ell\right),
\]
where $G_\ell$ is independent of $W_\ell$.  Combining this joint convergence with \eqref{e:wish} and using Slutsky's theorem proves \eqref{e:main}.  Equation \eqref{e:mainc} follows by subtracting the deterministic center $\phi_{p/n}(\theta_\ell)$.
\qed

\subsection{Proof of Theorem \ref{t:infy}}
\label{proof of theorem infinity}

\subsubsection{A specific CLT when $p/n \to 0$}

We now state the deterministic spiked covariance input needed for the sample-dominated regime. This point deserves a separate proof. The standard Bai--Yao theorem is usually stated for a fixed limiting aspect ratio $p/n\to c\in(0,1)$. In the present paper the aspect ratio itself tends to zero, and the required statement is a triangular-array statement with $c_p=p/n\to0$. Quoting the fixed-$c$ theorem and then sending $c$ to zero would leave a uniformity gap.

\begin{prop}[Low-aspect-ratio spike CLT for a deterministic covariance]\label{p:dclt}
Let \(\sigma^2>0\), let \(K\) be fixed, and let
\[
  \Lambda_p
  =\operatorname{diag}(\lambda_{p,1},\ldots,\lambda_{p,K},
  \sigma^2,\ldots,\sigma^2),
  \qquad \lambda_{p,j}\ge0.
\]
Fix \(\ell\le K\). Assume
\[
  \lambda_{p,\ell}\to \lambda_\ell^*>\sigma^2.
\]
Assume also that there are constants \(\delta_0>0\) and \(p_0\) such that, for every \(p\ge p_0\),
\begin{equation}
  \lambda_{p,\ell}-\sigma^2\ge\delta_0,
  \qquad
  |\lambda_{p,\ell}-\lambda_{p,j}|\ge\delta_0
  \quad(j\ne\ell),
  \label{e:p64gap}
\end{equation}
and
\begin{equation}
  \sup_p\max_{1\le j\le K}\lambda_{p,j}<\infty.
  \label{e:A1}
\end{equation}
Let \(y_{p,1},\ldots,y_{p,n}\) be independent centered real or complex Gaussian vectors with covariance \(\Lambda_p\), and set
\[
  S_p=\frac1n\sum_{k=1}^n y_{p,k}y_{p,k}^*.
\]
For a Hermitian matrix \(M\), write its eigenvalues as
\[
  \mu_1(M)\ge\cdots\ge\mu_p(M).
\]
For \(p\ge p_0\), define
\[
  k_{p,\ell}
  =1+\#\{j\le K:\lambda_{p,j}>\lambda_{p,\ell}\},
\]
and define \(k_{p,\ell}\) in any way for the finitely many values \(p<p_0\). Put
\begin{equation}
  \widehat\lambda_{p,\ell}^{\mathrm{det}}
  =\mu_{k_{p,\ell}}(S_p).
  \label{e:p64label}
\end{equation}
If \(p,n\to\infty\) and \(p/n\to0\), then
\begin{equation}
  \sqrt n\left(
    \widehat\lambda_{p,\ell}^{\mathrm{det}}
    -\lambda_{p,\ell}
    -\frac pn\frac{\sigma^2\lambda_{p,\ell}}
    {\lambda_{p,\ell}-\sigma^2}
  \right)
  \xrightarrow{\ \mathrm{law}\ }
  N\!\left(0,\kappa(\lambda_\ell^*)^2\right).
  \label{e:dclt}
\end{equation}
\end{prop}

The proof of Proposition \ref{p:dclt} is provided in Appendix \ref{appendix: proof of p:dclt}

\begin{rem}[Why the $p^2/n \to \tau \in [0,\infty[$ regime proof is not sufficient here]\label{r:pert}
In the $p^2/n \to \tau \in [0,\infty[$ regime, expanding only up to second order is enough, because the third-order remainder is negligible after multiplying by $p$. In the regime $p^2/n\to\infty$, one works at the $\sqrt n$ scale, and a crude operator-norm bound for a fixed-order remainder is not uniform enough. Proposition~\ref{p:dclt} avoids this difficulty by using the Schur complement equation and the deterministic outlier location, which already resums the repeated interaction with the bulk.
\end{rem}

\subsubsection{Proof of Theorem \ref{t:infy}}

\begin{proof}
We condition on the Haar matrix, apply Proposition~\ref{p:dclt}, and then show that the Haar fluctuation disappears at the $\sqrt n$ scale.

\noindent \textbf{Step 1: Conditioning on the Haar matrix.}
Condition on $\mathcal F_p=\sigma(O_p)$. By Lemma~\ref{l:cond}, the conditional model is a deterministic spiked covariance model with population eigenvalues
\[
  \lambda_{p,1},\ldots,\lambda_{p,r+s},\sigma^2,\ldots,\sigma^2.
\]
By Proposition~\ref{p:pop},
\[
  \lambda_{p,q}\xrightarrow{\mathbb P}\sigma^2+\theta_q,
  \qquad q=1,\ldots,r+s.
\]
Since the limiting spikes are distinct, the separation assumptions in Proposition~\ref{p:dclt} hold with probability tending to one. The boundedness assumption \eqref{e:A1} holds for every \(p\), because
\[
  \|A_p+O_pB_pO_p^*\|_{\mathrm{op}}
  \le \|A_p\|_{\mathrm{op}}+\|B_p\|_{\mathrm{op}}.
\]

Let
\[
  X_{p,\ell}=
  \sqrt n\left(
    \widehat\lambda_{p,\ell}
    -\lambda_{p,\ell}
    -\frac pn\frac{\sigma^2\lambda_{p,\ell}}{\lambda_{p,\ell}-\sigma^2}
  \right),
\]
and let
\[
  \phi_p(t)
  =\mathbb E\left[e^{\mathrm{i}tX_{p,\ell}}\mid\mathcal F_p\right].
\]
We now justify the conditional limit without assuming a uniform version of Proposition~\ref{p:dclt}. Take any subsequence. Since the finite vector
\((\lambda_{p,1},\ldots,\lambda_{p,r+s})\) converges in probability, this subsequence has a further subsequence on which the vector converges almost surely. Fix a sample point in this probability-one event. Along the further subsequence, the population eigenvalues form a deterministic triangular array that satisfies Proposition~\ref{p:dclt}. Therefore, for every fixed \(t\),
\[
  \phi_p(t)
  \to 
  \exp\left(-\frac12\kappa(\sigma^2+\theta_\ell)^2t^2\right)
\]
almost surely along that further subsequence. Since this can be done for every subsequence, we conclude that
\[
  \phi_p(t)
  \xrightarrow{\mathbb P}
  \exp\left(-\frac12\kappa(\sigma^2+\theta_\ell)^2t^2\right).
\]
Lemma~\ref{l:deco} therefore gives the unconditional convergence
\begin{equation}
  X_{p,\ell}
  \xrightarrow{\ \mathrm{law}\ }
  \mathcal N (0,\kappa(\sigma^2+\theta_\ell)^2).
  \label{e:condi}
\end{equation}

\noindent \textbf{Step 2: Replacing the random center by a deterministic center.}
Let
\[
  f(x)=\frac{\sigma^2x}{x-
  \sigma^2}.
\]
Set $I_\ell:=[\sigma^2+\theta_\ell/2,\;\sigma^2+2\theta_\ell]$. This interval does not
contain $\sigma^2$, so $f$ is $C^1$ on $I_\ell$ and there is a constant $C_\ell$, depending
only on $\sigma^2$ and $\theta_\ell$, with $|f'|\le C_\ell$ on $I_\ell$. Proposition~\ref{p:pop} implies
\begin{equation}
  \lambda_{p,\ell}-\sigma^2-\theta_\ell=
  \mathrm{O}_{\mathbb P}(1/p).
  \label{e:popo}
\end{equation}
In particular $\lambda_{p,\ell}\xrightarrow{\mathbb P}\sigma^2+\theta_\ell$, so
$\mathbb P(\lambda_{p,\ell}\in I_\ell)\to1$. On this event, both $\lambda_{p,\ell}$ and
$\sigma^2+\theta_\ell$ lie in $I_\ell$, so the mean value theorem applies and gives
\begin{align}
 &\sqrt n\frac pn
 \left|
   \frac{\sigma^2\lambda_{p,\ell}}{\lambda_{p,\ell}-\sigma^2}
   -
   \frac{\sigma^2(\sigma^2+\theta_\ell)}{\theta_\ell}
 \right|\notag\\
 &\quad
 \leq
 C\sqrt n\frac pn
 |\lambda_{p,\ell}-\sigma^2-\theta_\ell|
 =
 C\frac1{\sqrt n}\,
 p|\lambda_{p,\ell}-\sigma^2-\theta_\ell|
 =o_{\mathbb P}(1).
  \label{e:centc}
\end{align}
We used $p|\lambda_{p,\ell}-\sigma^2-\theta_\ell|=\mathrm{O}_{\mathbb P}(1)$ and $n\to\infty$. This holds on an event of probability tending to one, hence
\eqref{e:centc} holds as stated.

\noindent \textbf{Step 3: The Haar population fluctuation vanishes at the $\sqrt n$ scale.}
Again by \eqref{e:popo},
\begin{equation}
  \sqrt n(\lambda_{p,\ell}-\sigma^2-\theta_\ell)
  =\frac{\sqrt n}{p}
  \cdot
  p(\lambda_{p,\ell}-\sigma^2-\theta_\ell)
  =o_{\mathbb P}(1),
  \label{e:hvan}
\end{equation}
because $p^2/n\to\infty$ is equivalent to $\sqrt n/p\to0$.

Finally, decompose
\begin{align*}
 &\widehat\lambda_{p,\ell}-\sigma^2-\theta_\ell
 -\frac pn\frac{\sigma^2(\sigma^2+\theta_\ell)}{\theta_\ell}\\
 &\quad=
 \left(
    \widehat\lambda_{p,\ell}
    -\lambda_{p,\ell}
    -\frac pn\frac{\sigma^2\lambda_{p,\ell}}{\lambda_{p,\ell}-\sigma^2}
 \right)\\
 &\qquad+
 (\lambda_{p,\ell}-\sigma^2-\theta_\ell)\\
 &\qquad+
 \frac pn
 \left(
   \frac{\sigma^2\lambda_{p,\ell}}{\lambda_{p,\ell}-\sigma^2}
   -
   \frac{\sigma^2(\sigma^2+\theta_\ell)}{\theta_\ell}
 \right).
\end{align*}
Multiplying by $\sqrt n$, the first term converges to the Gaussian law in \eqref{e:condi}, while the second and third terms vanish by \eqref{e:hvan} and \eqref{e:centc}. Slutsky's theorem proves \eqref{e:infy2}.
\end{proof}

\subsection{Proof of Theorem~\ref{t:prop}}\label{proof of theorem prop}

\subsubsection{A lemma to deal with the fact that the $\lambda_{p,\ell}$ are not fixed}

\label{s:A lemma}

\begin{lem}[The CLT holds under a small perturbation of the spikes]\label{l:spikecont}
Let $\vartheta_{p,1},\ldots,\vartheta_{p,K}$ be possibly random positive numbers, and
let $Z_p$ be as in Corollary~\ref{c:repr}. Put
\[
  \Lambda_p=\diag(\sigma^2+\vartheta_{p,1},\ldots,\sigma^2+\vartheta_{p,K},
  \sigma^2,\ldots,\sigma^2),
\]
and let $\Lambda_p^{\mathrm{id}}$ be as in Section~\ref{s:quoted}. Consider the two
sample covariance matrices built from the same $Z_p$,
\[
  \widehat\Sigma_p=\Lambda_p^{1/2}\Bigl(\frac1nZ_pZ_p^*\Bigr)\Lambda_p^{1/2},
  \qquad
  \widehat\Sigma_p^{\mathrm{id}}
  =(\Lambda_p^{\mathrm{id}})^{1/2}\Bigl(\frac1nZ_pZ_p^*\Bigr)
   (\Lambda_p^{\mathrm{id}})^{1/2},
\]
with ordered eigenvalues $\widehat\lambda_{p,\ell}$ and
$\widehat\lambda_{p,\ell}^{\mathrm{id}}$. Assume $p/n\to c\in(0,\infty)$ and
$\vartheta_{p,\ell}-\theta_\ell=\mathrm O_{\mathbb P}(1/p)$ for every $\ell\le K$.
Then
\[
  \max_{\ell\le p}
  \sqrt n\,\bigl|\widehat\lambda_{p,\ell}
  -\widehat\lambda_{p,\ell}^{\mathrm{id}}\bigr|
  =o_{\mathbb P}(1).
\]
\end{lem}

\begin{proof}
\noindent \textbf{Step 1: the two population matrices are close.} By construction $\Lambda_p$
and $\Lambda_p^{\mathrm{id}}$ differ only in their first $K$ diagonal entries, and
there the difference is $\vartheta_{p,\ell}-\theta_\ell$. Hence
\[
  \|\Lambda_p-\Lambda_p^{\mathrm{id}}\|_{\mathrm{op}}
  =\max_{\ell\le K}|\vartheta_{p,\ell}-\theta_\ell|
  =\mathrm O_{\mathbb P}(1/p).
\]
In particular, on an event $\Omega_p$ with
$\mathbb P(\Omega_p)\to1$, all the diagonal entries of $\Lambda_p$ and of
$\Lambda_p^{\mathrm{id}}$ lie in the fixed compact set
$[\sigma^2/2,\sigma^2+\theta_1+1]\subset(0,\infty)$. The same bound then holds for
the square roots: since $t\mapsto\sqrt t$ is Lipschitz on any interval bounded away
from $0$, there is a constant $L$ with
\[
  \|\Lambda_p^{1/2}-(\Lambda_p^{\mathrm{id}})^{1/2}\|_{\mathrm{op}}
  \le L\,\|\Lambda_p-\Lambda_p^{\mathrm{id}}\|_{\mathrm{op}}
  =\mathrm O_{\mathbb P}(1/p).
\]

\noindent \textbf{Step 2: the two sample matrices are close.} Put
$\Upsilon_p:=\Lambda_p^{1/2}-(\Lambda_p^{\mathrm{id}})^{1/2}$. Adding and subtracting,
\[
  \widehat\Sigma_p-\widehat\Sigma_p^{\mathrm{id}}
  =\Upsilon_p\Bigl(\frac1nZ_pZ_p^*\Bigr)\Lambda_p^{1/2}
  +(\Lambda_p^{\mathrm{id}})^{1/2}\Bigl(\frac1nZ_pZ_p^*\Bigr)\Upsilon_p ,
\]
so, using $\|MN\|_{\mathrm{op}}\le\|M\|_{\mathrm{op}}\|N\|_{\mathrm{op}}$,
\[
  \|\widehat\Sigma_p-\widehat\Sigma_p^{\mathrm{id}}\|_{\mathrm{op}}
  \le\|\Upsilon_p\|_{\mathrm{op}}
  \Bigl\|\frac1nZ_pZ_p^*\Bigr\|_{\mathrm{op}}
  \bigl(\|\Lambda_p^{1/2}\|_{\mathrm{op}}
  +\|(\Lambda_p^{\mathrm{id}})^{1/2}\|_{\mathrm{op}}\bigr).
\]
In the proportional regime $\|n^{-1}Z_pZ_p^*\|_{\mathrm{op}}
=\mathrm O_{\mathbb P}(1)$, since its largest eigenvalue
converges almost surely to $(1+\sqrt c)^2$ (Bai--Silverstein \cite{bai10}),
and the two square-root norms are $\mathrm O_{\mathbb P}(1)$ by Step~1. With
$\|\Upsilon_p\|_{\mathrm{op}}=\mathrm O_{\mathbb P}(1/p)$,
\[
  \|\widehat\Sigma_p-\widehat\Sigma_p^{\mathrm{id}}\|_{\mathrm{op}}
  =\mathrm O_{\mathbb P}(1/p).
\]

\noindent \textbf{Step 3: from matrices to eigenvalues.} By Weyl's inequality
(Lemma~\ref{a:weyl}), for every $\ell$,
\[
  |\widehat\lambda_{p,\ell}-\widehat\lambda_{p,\ell}^{\mathrm{id}}|
  \le\|\widehat\Sigma_p-\widehat\Sigma_p^{\mathrm{id}}\|_{\mathrm{op}}
  =\mathrm O_{\mathbb P}(1/p).
\]
Multiplying by $\sqrt n$ and using $\sqrt n/p=(1/\sqrt n)\cdot(n/p)\to0$ in the
proportional regime, since $n/p\to1/c$ and $1/\sqrt n\to0$, gives the claim.
\end{proof}

\subsubsection{Proof of Theorem \ref{t:prop}}

Put $K=r+s$. Part (i) follows from a finite-rank comparison with
a white Wishart matrix on the original probability space.
For parts (ii)--(iv), we use the following coupling with a model
whose spikes are fixed.

\noindent\textbf{The ideal model.}
For parts (ii)--(iv), we work with the representation of
Corollary~\ref{c:repr}: the eigenvalues of $\widehat\Sigma_p$
have the same law as those of
$\Lambda_p^{1/2}(n^{-1}Z_pZ_p^*)\Lambda_p^{1/2}$,
where $Z_p$ is independent of $\mathcal F_p$.
Using the same $Z_p$, we build the ideal model
$\widehat\Sigma_p^{\mathrm{id}}$ of Section~\ref{s:A lemma}.
Its population matrix $\Lambda_p^{\mathrm{id}}$ is deterministic
and has the fixed spikes $\theta_1,\ldots,\theta_K$.
Thus Theorems~\ref{t:bbp}, \ref{t:baiyao}
and~\ref{t:twedge} apply directly to this ideal model under
their respective hypotheses.
All comparisons below are made in this coupling.

\noindent\textbf{The comparison bound.}
By Proposition~\ref{p:pop},
$\vartheta_{p,\ell}:=\lambda_{p,\ell}-\sigma^2$ satisfies
$\vartheta_{p,\ell}-\theta_\ell=\mathrm O_{\mathbb P}(1/p)$
for every $\ell\le K$.
The proof of Lemma~\ref{l:spikecont} therefore gives
\begin{equation}
  \max_{\ell\le p}
  \bigl|\widehat\lambda_{p,\ell}
  -\widehat\lambda_{p,\ell}^{\mathrm{id}}\bigr|
  =\mathrm O_{\mathbb P}(1/p)
  =o_{\mathbb P}(n^{-1/2}).
  \label{e:propcomp}
\end{equation}
The last equality follows from
$\sqrt n/p=(n/p)/\sqrt n\to0$.

\begin{proof}[Proof of (i)]
Work on the original probability space and set
\[
  X_p=(x_{p,1},\ldots,x_{p,n}),
  \qquad Z_p^{(0)}=\Sigma_p^{-1/2}X_p.
\]
Conditionally on $O_p$, the entries of $Z_p^{(0)}$ are independent
standard real or circular complex Gaussian variables.
This conditional law does not depend on $O_p$, so $Z_p^{(0)}$
is also a standard Gaussian matrix unconditionally.
Put
\[
  M_p=\frac1n Z_p^{(0)}(Z_p^{(0)})^*,
  \qquad S_p^{(0)}=\sigma^2M_p,
  \qquad D_p=\Sigma_p^{1/2}-\sigma I_p,
\]
where $\sigma$ is the positive square root of $\sigma^2$.
Since $\Sigma_p-\sigma^2I_p$ has rank at most $K$,
the same is true of $D_p$.
Moreover,
\[
  \widehat\Sigma_p-S_p^{(0)}
  =D_pM_p\Sigma_p^{1/2}+\sigma M_pD_p,
\]
and hence
$\operatorname{rank}(\widehat\Sigma_p-S_p^{(0)})\le2K$.
For a Hermitian matrix $M$, let $F^M$ denote its empirical
spectral distribution function.
The rank inequality \cite{bai10} gives
\[
  \sup_{t\in\mathbb R}
  |F^{\widehat\Sigma_p}(t)-F^{S_p^{(0)}}(t)|
  \le\frac{\operatorname{rank}(\widehat\Sigma_p-S_p^{(0)})}{p}
  \le\frac{2K}{p}\longrightarrow0.
\]
By Theorem~\ref{a:mpesd}, the empirical spectral distribution
of $S_p^{(0)}$ converges weakly, almost surely, to $\mu_{c,\sigma^2}$.
The same conclusion therefore holds for $\widehat\Sigma_p$.
\end{proof}

\begin{proof}[Proof of (ii)]
Fix $\ell$ with $\theta_\ell>\sigma^2\sqrt c$; by assumption the spikes are distinct, so
$\theta_\ell$ is simple. We prove the central limit theorem; the convergence in
probability then follows, since the central limit theorem gives
$\widehat\lambda_{p,\ell}-\phi_{p/n}(\theta_\ell)=\mathrm O_{\mathbb P}(n^{-1/2})
=o_{\mathbb P}(1)$, and $\phi_{p/n}(\theta_\ell)\to\phi_c(\theta_\ell)$
because $p/n\to c$; adding the two gives
$\widehat\lambda_{p,\ell}\xrightarrow{\mathbb P}\phi_c(\theta_\ell)$.

\noindent \textbf{Step 1: the ideal model.} The ideal model
$\widehat\Sigma_p^{\mathrm{id}}$ of the preamble has fixed spikes
$\theta_1,\ldots,\theta_K$ and is built from the same Gaussian matrix $Z_p$. The spike
$\theta_\ell$ is supercritical and simple, so Theorem~\ref{t:baiyao} applies to the
ideal model and gives
\begin{equation}
  \sqrt n\,\bigl(\widehat\lambda_{p,\ell}^{\mathrm{id}}
  -\phi_{p/n}(\theta_\ell)\bigr)
  \xrightarrow{\ \mathrm{law}\ }\mathcal N(0,\kappa V_c(\theta_\ell)).
  \label{e:propideal}
\end{equation}

\noindent \textbf{Step 2: comparing the ideal model and our model.} By
Proposition~\ref{p:pop}, the true spikes satisfy
$(\lambda_{p,\ell}-\sigma^2)-\theta_\ell=\mathrm O_{\mathbb P}(1/p)$, and both models are
built from the same $Z_p$, so Lemma~\ref{l:spikecont} applies with $\vartheta_{p,\ell}=\lambda_{p,\ell}-\sigma^2$ and gives, through \eqref{e:propcomp},
\[
  \sqrt n\,\bigl(\widehat\lambda_{p,\ell}
  -\widehat\lambda_{p,\ell}^{\mathrm{id}}\bigr)
  =o_{\mathbb P}(1).
\]
Writing
\[
  \sqrt n\,\bigl(\widehat\lambda_{p,\ell}-\phi_{p/n}(\theta_\ell)\bigr)
  =\sqrt n\,\bigl(\widehat\lambda_{p,\ell}^{\mathrm{id}}
   -\phi_{p/n}(\theta_\ell)\bigr)
  +\sqrt n\,\bigl(\widehat\lambda_{p,\ell}
   -\widehat\lambda_{p,\ell}^{\mathrm{id}}\bigr),
\]
and combining \eqref{e:propideal} with the $o_{\mathbb P}(1)$ term, Slutsky's theorem
gives
\[
  \sqrt n\,\bigl(\widehat\lambda_{p,\ell}-\phi_{p/n}(\theta_\ell)\bigr)
  \xrightarrow{\ \mathrm{law}\ }\mathcal N(0,\kappa V_c(\theta_\ell)).
\]
Both inputs used here are unconditional. The ideal model $\widehat\Sigma_p^{\mathrm{id}}$
has deterministic spikes and is a function of $Z_p$ alone, so Theorem~\ref{t:baiyao}
applies to it directly; and Lemma~\ref{l:spikecont} is stated for possibly random
spikes, its hypothesis
$\vartheta_{p,\ell}-\theta_\ell=\mathrm O_{\mathbb P}(1/p)$ being supplied by
Proposition~\ref{p:pop}. So the randomness of $O_p$ is already absorbed by
Lemma~\ref{l:spikecont}, and no conditioning on $\mathcal F_p$ is needed.

\end{proof}

\begin{proof}[Proof of (iii)]
Fix $\ell$ with $0<\theta_\ell\le\sigma^2\sqrt c$.
Theorem~\ref{t:bbp}, applied to the ideal model with fixed spikes,
gives
\[
  \widehat\lambda_{p,\ell}^{\mathrm{id}}
  \xrightarrow{\ \mathrm{a.s.}\ }
  \sigma^2(1+\sqrt c)^2.
\]
By \eqref{e:propcomp},
\[
  \widehat\lambda_{p,\ell}
  -\widehat\lambda_{p,\ell}^{\mathrm{id}}
  =\mathrm O_{\mathbb P}(1/p)=o_{\mathbb P}(1).
\]
Combining these statements yields
$\widehat\lambda_{p,\ell}
\xrightarrow{\mathbb P}\sigma^2(1+\sqrt c)^2$.
\end{proof}

\begin{proof}[Proof of (iv)]
Assume $\theta_K>\sigma^2\sqrt c$.
From \eqref{e:john} and $p/n\to c\in(0,\infty)$, we obtain
\[
  s_{p,n}^+
  =\sigma^2(1+\sqrt c)^{4/3}c^{-1/6}
    n^{-2/3}(1+o(1)).
\]
Also, the centering satisfies the exact identity
\[
  m_{p,n}^+-\sigma^2\bigl(1+\sqrt{p/n}\bigr)^2
  =-\frac{\sigma^2}{n}
    \left(1+\frac{2\sqrt p}{\sqrt n+\sqrt{n-1}}\right)
  =\mathrm O(n^{-1}).
\]
All spikes of the ideal model are supercritical, so
Theorem~\ref{t:twedge} gives
\[
  \frac{\widehat\lambda_{p,K+1}^{\mathrm{id}}-m_{p,n}^+}{s_{p,n}^+}
  \xrightarrow{\ \mathrm{law}\ }\mathrm{TW}.
\]
On the other hand, \eqref{e:propcomp} implies
\[
  \frac{\widehat\lambda_{p,K+1}
       -\widehat\lambda_{p,K+1}^{\mathrm{id}}}{s_{p,n}^+}
  =\mathrm O_{\mathbb P}\left(\frac{n^{2/3}}p\right)
  =o_{\mathbb P}(1).
\]
Slutsky's theorem therefore yields
\[
  \frac{\widehat\lambda_{p,K+1}-m_{p,n}^+}{s_{p,n}^+}
  \xrightarrow{\ \mathrm{law}\ }\mathrm{TW}.
\]
Convergence in distribution implies tightness, so
$\widehat\lambda_{p,K+1}-m_{p,n}^+
=\mathrm O_{\mathbb P}(n^{-2/3})$.
Combining this with the centering identity and
$n^{-1}=o(n^{-2/3})$ proves
\[
  \widehat\lambda_{p,K+1}
  =\sigma^2\bigl(1+\sqrt{p/n}\bigr)^2
   +\mathrm O_{\mathbb P}(n^{-2/3}).
\]
\end{proof}

\subsection{Proofs of the divergence of the fluctuations with a wrong scaling}
\label{s:scalprf}

\begin{proof}[Proof of Proposition~\ref{p:diverge}]
We apply Lemma~\ref{a:scale}(ii) three times, and we check its two hypotheses in each
case: convergence in distribution under one scaling, with a limit having no atom at
zero, and divergence of the ratio of the two scalings.

(i) The assumption $p^2/n\to0$ is the case $\tau=0$ of Theorem~\ref{t:main}, so
\[
  p\,\bigl(\widehat\lambda_{p,\ell}-\phi_{p/n}(\theta_\ell)\bigr)
  \xrightarrow{\ \mathrm{law}\ }W_\ell ,
\]
and Corollary~\ref{a:atom} with $\tau=0$ gives $\mathbb P(W_\ell=0)=0$. Moreover
$(p/\sqrt n)^2=p^2/n\to0$ and $p/\sqrt n>0$, so $\sqrt n/p\to\infty$.
Lemma~\ref{a:scale}(ii) with $u_p=p$ and $v_p=\sqrt n$ gives the claim.

(ii) By Theorem~\ref{t:infy},
\[
  \sqrt n\,\bigl(\widehat\lambda_{p,\ell}-\phi_{p/n}(\theta_\ell)\bigr)
  \xrightarrow{\ \mathrm{law}\ }\mathcal N\bigl(0,\kappa V_0(\theta_\ell)\bigr) .
\]
The variance $\kappa V_0(\theta_\ell)=\kappa(\sigma^2+\theta_\ell)^2$ is positive, so
this law has no atom at zero by Corollary~\ref{a:atom}. Moreover
$p/\sqrt n=\sqrt{p^2/n}\to\infty$. Lemma~\ref{a:scale}(ii) with $u_p=\sqrt n$ and
$v_p=p$ gives the claim.

(iii) By Theorem~\ref{t:prop}(ii),
\[
  \sqrt n\,\bigl(\widehat\lambda_{p,\ell}-\phi_{p/n}(\theta_\ell)\bigr)
  \xrightarrow{\ \mathrm{law}\ }\mathcal N\bigl(0,\kappa V_c(\theta_\ell)\bigr) .
\]
The variance is positive: the assumption $\theta_\ell>\sigma^2\sqrt c$ gives
$\theta_\ell^2>c\sigma^4$, hence $1-c\sigma^4/\theta_\ell^2>0$ and
$V_c(\theta_\ell)>0$. So this law has no atom at zero by Corollary~\ref{a:atom}.
Moreover $p/\sqrt n=(p/n)\sqrt n\to\infty$, because $p/n\to c>0$ and $n\to\infty$.
Lemma~\ref{a:scale}(ii) with $u_p=\sqrt n$ and $v_p=p$ gives the claim.

The last assertion follows from Lemma~\ref{a:nolimit} in the three cases.
\end{proof}

\begin{proof}[Proof of Proposition~\ref{p:mixed}]
Theorem~\ref{t:main} gives
\[
  p\,\bigl(\widehat\lambda_{p,\ell}-\phi_{p/n}(\theta_\ell)\bigr)
  \xrightarrow{\ \mathrm{law}\ }
  W_\ell+\sqrt{\kappa\tau V_0(\theta_\ell)}\,G_\ell .
\]
The assumption $p^2/n\to\tau$ means $(p/\sqrt n)^2\to\tau$; since $p/\sqrt n>0$ and the
square root is continuous on $[0,\infty)$, this gives $p/\sqrt n\to\sqrt\tau$, and
$\tau>0$ gives $\sqrt n/p\to\tau^{-1/2}\in(0,\infty)$. Lemma~\ref{a:scale}(i) with
$u_p=p$, $v_p=\sqrt n$ and $\gamma=\tau^{-1/2}$ therefore gives
\[
  \sqrt n\,\bigl(\widehat\lambda_{p,\ell}-\phi_{p/n}(\theta_\ell)\bigr)
  \xrightarrow{\ \mathrm{law}\ }
  \tau^{-1/2}\Bigl(W_\ell+\sqrt{\kappa\tau V_0(\theta_\ell)}\,G_\ell\Bigr)
  =\tau^{-1/2}W_\ell+\sqrt{\kappa V_0(\theta_\ell)}\,G_\ell ,
\]
because $\tau^{-1/2}\sqrt{\kappa\tau V_0(\theta_\ell)}=\sqrt{\kappa V_0(\theta_\ell)}$.
\end{proof}

\subsection{Proofs of the properties of the spectral test}
\label{proofs of the properties of the test}

In section \ref{subsubsection: Proof of consistency and level of the test}, we prove the consistency and level of our test under assumptions (A1)-(A7). In section \ref{subsubsection:Proofs of the properties of the test when relaxing assumption (A6}, we prove the properties of the test when relaxing assumption (A6).

\subsubsection{Proof of consistency and level of the test}
\label{subsubsection: Proof of consistency and level of the test}

We first isolate the continuity property needed to replace the unknown reference-law parameters by consistent estimates. The statement also covers a vanishing Gaussian component and weights of either sign.
\begin{lem}[Continuity of the reference law and its quantile]
\label{lem:quantile-continuity}
Fix $\gamma>0$, $a\in(0,1)$, and an integer $m\ge1$. Let $Z_1,\ldots,Z_m$ be independent $\chi^2_1$ variables in the real case or independent $\mathrm{Exp}(1)$ variables in the complex case, and let $G\sim\mathcal N(0,1)$ be independent of them. For $b_j\notin\{0,\gamma\}$ and $v\ge0$, define
\[
  L(b,v)=\sum_{j=1}^m\frac{\gamma b_j}{\gamma-b_j}Z_j+vG,
  \qquad F_{b,v}(x)=\mathbb P(L(b,v)\le x),
\]
and $q(b,v)=\inf\{x:F_{b,v}(x)\ge1-a\}$.
The law of $L(b,v)$ is absolutely continuous, and its cdf is strictly increasing on the interior of its support. Its quantile lies in that interior and satisfies $F_{b,v}(q(b,v))=1-a$. If $(b,v)\to(\beta,v_0)$, where
$\beta_j\notin\{0,\gamma\}$ and $v_0\ge0$, then
\[
  F_{b,v}(x)\to F_{\beta,v_0}(x)\quad\text{for every }x\in\mathbb R,
  \qquad q(b,v)\to q(\beta,v_0).
\]
\end{lem}
\begin{proof}
\emph{(A) Regularity.}
All coefficients are nonzero, so Lemma~\ref{a:conv} gives absolute continuity. The support is $\mathbb R$ if $v>0$ or if the coefficients have both signs; otherwise it is one of the two closed half-lines with endpoint zero. Every nonempty open interval in the interior of this support has positive probability, since the $Z_j$ have positive densities on $(0,\infty)$ and $G$ has a positive density on $\mathbb R$. Thus the cdf is continuous and strictly increasing on that interior. Since $1-a\in(0,1)$, its quantile is an interior point and its cdf value there is $1-a$.

\emph{(B) Continuity.}
Using the same auxiliary variables for all parameters gives
$L(b,v)\to L(\beta,v_0)$ almost surely. Absolute continuity of the limiting law implies convergence of the cdfs at every fixed $x$. Write $q=q(\beta,v_0)$. For every sufficiently small $\varepsilon>0$,
\[
  F_{\beta,v_0}(q-\varepsilon)<1-a
  <F_{\beta,v_0}(q+\varepsilon).
\]
The same inequalities hold with $F_{b,v}$ for all parameters sufficiently close to $(\beta,v_0)$, and imply $q-\varepsilon<q(b,v)\le q+\varepsilon$. This proves quantile continuity.
\end{proof}
We now prove Theorem~\ref{thm:test}. Steps~1--3 establish consistency of the nuisance-parameter estimates. Steps~4--5 prove convergence of the plug-in and Monte-Carlo critical values, and Steps~6--7 establish level and consistency.

\paragraph{Notation and tools.}
We write $\lambda_{p,k}:=\lambda_k(\Sigma_p)$ for the $k$-th largest eigenvalue of the
population covariance. Since $A_p+O_pB_pO_p^*$ has rank at most $r+s$, we have
$\lambda_{p,k}=\sigma^2$ for $k>r+s$, and for $k\le r+s$ the number
$\lambda_{p,k}-\sigma^2$ is the $k$-th largest eigenvalue of $A_p+O_pB_pO_p^*$, which
converges to $\theta_k$. We use three facts repeatedly.

(T1) Weyl's inequality (Lemma~\ref{a:weyl}):
$\max_k|\widehat\lambda_{p,k}-\lambda_{p,k}|\le\|\widehat\Sigma_p-\Sigma_p\|_{\mathrm{op}}$.

(T2) Operator-norm bound: by Corollary~\ref{a:wishartop} and Remark~\ref{a:wishartrand}, since
$\|\Sigma_p\|_{\mathrm{op}}$ is bounded by a deterministic constant, and since $p/n\to0$ in this regime,
$\|\widehat\Sigma_p-\Sigma_p\|_{\mathrm{op}}=\mathrm O_{\mathbb P}(\sqrt{p/n})$.

(T3) Eigenvalue limit: by Theorem~\ref{t:main}, $\widehat\lambda_{p,\ell}\xrightarrow{\mathbb P}\sigma^2+\theta_\ell$.

\paragraph{Step 1: $\tilde\sigma^2\xrightarrow{\mathbb P}\sigma^2$ and $\hat\sigma^2-\sigma^2=\mathrm O_{\mathbb P}(1/\sqrt{pn})$.}
Split $\tilde\sigma^2-\sigma^2=\frac1p(\operatorname{tr}\widehat\Sigma_p-\operatorname{tr}\Sigma_p)+\frac1p(\operatorname{tr}\Sigma_p-p\sigma^2)$.
The deterministic part equals
$\frac1p\sum_{\ell\le r+s}(\lambda_{p,\ell}-\sigma^2)=\mathrm O(1/p)$, a sum of $r+s$
bounded terms. For the random part, using
$\operatorname{tr}\widehat\Sigma_p=\frac1n\sum_{k=1}^n\|x_{p,k}\|_2^2$, this is,
conditionally on $\mathcal F_p$, an average of $n$ i.i.d.\ terms of common mean
$\mathbb E[\|x_{p,1}\|_2^2\mid\mathcal F_p]=\operatorname{tr}\Sigma_p$. Writing
$x_{p,1}=\Sigma_p^{1/2}z$ with $z\sim\mathcal N_{\mathbb F}(0,I_p)$ and diagonalising,
$\|x_{p,1}\|_2^2=\sum_k\lambda_{p,k}|\tilde z_k|^2$ with $|\tilde z_k|^2$ i.i.d.\ of
variance $\kappa$, so
$\operatorname{Var}(\|x_{p,1}\|_2^2\mid\mathcal F_p)=\kappa\operatorname{tr}(\Sigma_p^2)$
and $\operatorname{Var}(\operatorname{tr}\widehat\Sigma_p\mid\mathcal F_p)=\frac\kappa n\operatorname{tr}(\Sigma_p^2)$.
Since $\operatorname{tr}(\Sigma_p^2)=(p-r-s)\sigma^4+\sum_{k\le r+s}\lambda_{p,k}^2=\mathrm O(p)$
with a constant depending only on $\sigma^2$ and the $\theta_k$, the conditional variance is bounded uniformly in
$\mathcal F_p$:
\[
  \operatorname{Var}\bigl(\operatorname{tr}\widehat\Sigma_p\mid\mathcal F_p\bigr)
  =\frac\kappa n\operatorname{tr}(\Sigma_p^2)\le C\,\frac pn,
\]
for a constant $C$ not depending on the realisation of $O_p$. The conditional Chebyshev inequality gives, for every
$t>0$,
\[
  \mathbb P\bigl(|\operatorname{tr}\widehat\Sigma_p-\operatorname{tr}\Sigma_p|\ge t\mid\mathcal F_p\bigr)
  \le\frac{\operatorname{Var}(\operatorname{tr}\widehat\Sigma_p-\operatorname{tr}\Sigma_p\mid\mathcal F_p)}{t^2} = \frac{\operatorname{Var}(\operatorname{tr}\widehat\Sigma_p\mid\mathcal F_p)}{t^2}\le\frac{C\,p/n}{t^2},
\]
and the right-hand side is deterministic; taking expectations over $\mathcal F_p$ (the tower
property) yields $\mathbb P(|\operatorname{tr}\widehat\Sigma_p-\operatorname{tr}\Sigma_p|\ge t)\le C\,(p/n)/t^2$. With $t=K\sqrt{p/n}$ this is
$\le C/K^2$, uniformly in $p$, so $\operatorname{tr}\widehat\Sigma_p-\operatorname{tr}\Sigma_p=\mathrm O_{\mathbb P}(\sqrt{p/n})$. Dividing by $p$,
\[
  \tfrac1p\bigl(\operatorname{tr}\widehat\Sigma_p-\operatorname{tr}\Sigma_p\bigr)
  =\mathrm O_{\mathbb P}\bigl(\tfrac1p\sqrt{p/n}\bigr)=\mathrm O_{\mathbb P}(1/\sqrt{pn}).
\]
Together with the deterministic part $\frac1p(\operatorname{tr}\Sigma_p-p\sigma^2)=\mathrm O(1/p)$,
this gives $\tilde\sigma^2-\sigma^2=\mathrm O_{\mathbb P}(1/\sqrt{pn})+\mathrm O(1/p)\to0$.

For $\hat\sigma^2$, we work on $\{\widehat{r+s}=r+s\}$, whose probability tends to $1$ by
Step~2. On this event, $\hat\sigma^2=\frac1{p-r-s}\sum_{k>r+s}\widehat\lambda_{p,k}$, and since
$\lambda_{p,k}=\sigma^2$ for $k>r+s$, subtracting gives
$(p-r-s)(\hat\sigma^2-\sigma^2)=\sum_{k>r+s}(\widehat\lambda_{p,k}-\lambda_{p,k})$.
Writing $\sum_{k>r+s}=\sum_{k=1}^p-\sum_{k\le r+s}$,
\[
  (p-r-s)(\hat\sigma^2-\sigma^2)
  =\underbrace{(\operatorname{tr}\widehat\Sigma_p-\operatorname{tr}\Sigma_p)}_{=:S_1}
  -\underbrace{\sum_{k\le r+s}(\widehat\lambda_{p,k}-\lambda_{p,k})}_{=:S_2}.
\]
For $S_1$, the bound above gives $S_1=\mathrm O_{\mathbb P}(\sqrt{p/n})$. For $S_2$, each of
the $r+s$ terms is $\le\|\widehat\Sigma_p-\Sigma_p\|_{\mathrm{op}}$ by (T1), so
$|S_2|\le(r+s)\|\widehat\Sigma_p-\Sigma_p\|_{\mathrm{op}}=\mathrm O_{\mathbb P}(\sqrt{p/n})$
by (T2). Dividing by $p-r-s$, which is of order $p$,
\[
  |\hat\sigma^2-\sigma^2|\le\frac{|S_1|+|S_2|}{p-r-s}=\mathrm O_{\mathbb P}(\sqrt{p/n}/p)=\mathrm O_{\mathbb P}(1/\sqrt{pn}).
\]
This holds on $\{\widehat{r+s}=r+s\}$, hence unconditionally. We record the consequence used
below, at scale $p$: $p(\hat\sigma^2-\sigma^2)=\mathrm O_{\mathbb P}(\sqrt{p/n})\to0$.

\paragraph{Step 2: $\mathbb P(\widehat{r+s}=r+s)\to1$.}
Write $\tilde\theta_p:=\tilde\sigma^2(1+\sqrt{p/n})^2+\epsilon_p$. The count equals $r+s$ when
the $r+s$ spike eigenvalues exceed $\tilde\theta_p$ and the $p-r-s$ bulk eigenvalues fall
below it. From Step~1, we can write
\[\tilde\sigma^2-\sigma^2=R_p+D_p,\qquad R_p:=\tfrac1p\bigl(\operatorname{tr}\widehat\Sigma_p-\operatorname{tr}\Sigma_p\bigr),\quad D_p:=\tfrac1p\operatorname{tr}\Sigma_p-\sigma^2=\tfrac1p\textstyle\sum_{k\le r+s}(\lambda_{p,k}-\sigma^2).\]

with $D_p=\mathrm O(1/p)$, $R_p = \mathrm O_\mathbb P(1/\sqrt{p n})$, so
\begin{equation}
  \tilde\theta_p-\sigma^2=(D_p+R_p)(1+\sqrt{p/n})^2+\sigma^2(2\sqrt{p/n}+p/n)+\epsilon_p.
  \label{eq:threshold-decomp}
\end{equation}

For every bulk index $k>r+s$, $\lambda_{p,k}=\sigma^2$, so by (T1)--(T2)
$\widehat\lambda_{p,k}-\sigma^2\le\|\widehat\Sigma_p-\Sigma_p\|_{\mathrm{op}}=\mathrm O_{\mathbb P}(\sqrt{p/n})$
uniformly in $k$. Dropping the nonnegative terms in \eqref{eq:threshold-decomp},
$\tilde\theta_p-\sigma^2\ge R_p(1+\sqrt{p/n})^2+\epsilon_p$, hence
\[
  \widehat\lambda_{p,k}-\tilde\theta_p=\widehat\lambda_{p,k}-\sigma^2+\sigma^2 -\tilde\theta_p\le-\epsilon_p+\mathrm O_{\mathbb P}(\sqrt{p/n})+\mathrm O_{\mathbb P}(1/\sqrt{pn}).
\]
Both error terms are $\mathrm o_{\mathbb P}(\epsilon_p)$, since $\sqrt{p/n}=\mathrm o(\epsilon_p)$ by
(A3) and $1/\sqrt{pn}\le\sqrt{p/n}$. So with probability tending to one every bulk
eigenvalue is below $\tilde\theta_p$.

For $k\le r+s$, $\widehat\lambda_{p,k}\xrightarrow{\mathbb P}\sigma^2+\theta_k$ by
Theorem~\ref{t:main}.
On the other hand \eqref{eq:threshold-decomp} gives
$\tilde\theta_p\xrightarrow{\mathbb P}\sigma^2$, since $D_p\to0$, $R_p\xrightarrow{\mathbb P}0$,
$p/n\to0$ and $\epsilon_p\to0$. Hence
$\widehat\lambda_{p,k}-\tilde\theta_p\xrightarrow{\mathbb P}\theta_k>0$, so every spike
exceeds the threshold with probability tending to one.

So $\widehat{r+s}=r+s$ with probability tending to one, and, $r$ being known,
$\hat s=\widehat{r+s}-r=s$ on the same event.

\paragraph{Step 3: the centered eigenvalues converge.}
On $\{\widehat{r+s}=r+s\}$, for each $\ell\in[r+s]$,
\[
  \hat\mu_\ell-\theta_\ell
  =\underbrace{(\widehat\lambda_{p,\ell}-\lambda_{p,\ell})}_{(\mathrm i)}
  +\underbrace{\bigl((\lambda_{p,\ell}-\sigma^2)-\theta_\ell\bigr)}_{(\mathrm{ii})}
  +\underbrace{(\sigma^2-\hat\sigma^2)}_{(\mathrm{iii})},
\]
where $\hat\mu_\ell=\widehat\lambda_{p,\ell}-\hat\sigma^2$. Term (i) is
$\mathrm O_{\mathbb P}(\sqrt{p/n})$ by (T1)--(T2); term (ii) is the gap between the $\ell$-th
eigenvalue of $A_p+O_pB_pO_p^*$ and its limit $\theta_\ell$, which is
$\mathrm O_{\mathbb P}(1/p)$ by Proposition~\ref{p:pop}; term (iii) is
$\mathrm O_{\mathbb P}(1/\sqrt{pn})$ by Step~1. All three tend to $0$, so
$\hat\mu_\ell\xrightarrow{\mathbb P}\theta_\ell$. Applying this to the ranks labelled as
batch effects gives $\hat\beta_j\xrightarrow{\mathbb P}\beta_j$ for every $j\in[s]$.

\paragraph{Step 4: convergence of the plug-in reference law.}
Relabel the batch strengths in decreasing order, and order their
estimates in the same way. This does not change the reference law.
By Steps~1--3 and (A6), with probability tending to one, $\hat s=s$
and the selected batch estimates satisfy
\[
  \max_{1\le j\le s}|\hat\beta_j-\beta_j|
  \xrightarrow{\mathbb P}0,
  \qquad
  \hat v:=\sqrt{\kappa p^2/n}\,(\hat\sigma^2+\gamma)
  \xrightarrow{\mathbb P}
  v_0:=\sqrt{\kappa\tau}\,(\sigma^2+\gamma).
\]
Assumption (A2) gives
$\min_j\min\{|\beta_j|,|\beta_j-\gamma|\}>0$.
Hence there is an event $E_p$, with $\mathbb P(E_p)\to1$, on which
$\hat s=s$, all required ranks are available, and
$\hat\beta_j\notin\{0,\gamma\}$ for every $j$.
For the convergence statements below, extend the fitted parameters
by $(\beta,v_0)$ on $E_p^c$; this is only a proof convention.

Put $F_0=F_{\beta,v_0}$, $q=q(\beta,v_0)$, and
$\widehat F_p=F_{\hat\beta,\hat v}$.
Lemma~\ref{lem:quantile-continuity} and the continuous mapping theorem give
\[
  q(\hat\beta,\hat v)\xrightarrow{\mathbb P}q,
  \qquad
  \widehat F_p(x)\xrightarrow{\mathbb P}F_0(x)
  \quad\text{for every fixed }x\in\mathbb R.
\]
These conclusions hold under both hypotheses. Under $H_0$, $F_0$ is
the cdf of the limiting test statistic before the negligible
noise-estimation correction.

\paragraph{Step 5: convergence of the Monte-Carlo critical value.}
Write $u=1-a$, and let $\mathcal D_p$ be the sigma-field generated by
the data. Conditionally on $\mathcal D_p$, generate the $M=M_p$ draws
independently from $\widehat L$ using fresh auxiliary randomness.
Their empirical cdf
\[
  F_{p,M}(x)=\frac1M\sum_{m=1}^M\mathbf1_{\{c_m\le x\}}
\]
satisfies, on $E_p$,
\[
  \mathbb E[F_{p,M}(x)\mid\mathcal D_p]=\widehat F_p(x),
  \qquad
  \operatorname{Var}(F_{p,M}(x)\mid\mathcal D_p)\le\frac1{4M}.
\]
Fix a sufficiently small $\varepsilon>0$. By
Lemma~\ref{lem:quantile-continuity},
$F_0(q-\varepsilon)<u<F_0(q+\varepsilon)$.
Choose $\delta>0$ smaller than both of these cdf gaps.
Step~4 implies that
\[
  H_p=E_p\cap
  \{\widehat F_p(q-\varepsilon)\le u-\delta,\ 
    \widehat F_p(q+\varepsilon)\ge u+\delta\}
\]
has probability tending to one. Conditional Chebyshev inequalities
at $q-\varepsilon$ and $q+\varepsilon$, followed by integration over
the data, give
\[
  \mathbb P(|\hat q_{1-a}-q|>\varepsilon)
  \le\mathbb P(H_p^c)+\frac1{2M\delta^2}.
\]
Indeed, if the empirical cdf is below $u$ at $q-\varepsilon$ and
at least $u$ at $q+\varepsilon$, then
$c_{(\lceil Mu\rceil)}\in(q-\varepsilon,q+\varepsilon]$.
Since $M_p\to\infty$, this proves
$\hat q_{1-a}\xrightarrow{\mathbb P}q$.
We henceforth write $q_{1-a}=q$. No rate condition relating $M_p$ to $p$ or $n$ is required.

\paragraph{Step 6: asymptotic level.}
Under $H_0$ ($\theta_\ell=\gamma$),
\[
  T_\ell=\underbrace{p(\widehat\lambda_{p,\ell}-\phi_{p/n}(\gamma))}_{(\mathrm A)}
  +\underbrace{p(\phi_{p/n}(\gamma)-\hat\phi(\gamma))}_{(\mathrm B)}.
\]
Since $\theta_\ell=\gamma$, Theorem~\ref{t:main} applies at rank $\ell$.
By (A5), $\theta_\ell=\alpha_i$ for some signal index $i$, and
\[
  W_\ell(\gamma)=
  \sum_{j=1}^s\frac{\gamma\beta_j}{\gamma-\beta_j}|g_{ij}|^2.
\]
Thus
\[
  (\mathrm A)\xrightarrow{\ \mathrm{law}\ }L(\gamma)=W_\ell(\gamma)+v_0G.
\]
For (B), with $\phi_c(y)=\sigma^2+y+c\,\sigma^2(\sigma^2+y)/y$, $c=p/n$, and
$\hat\phi(y)=\hat\sigma^2+y+\frac{p\,\hat\sigma^2(\hat\sigma^2+y)}{n y}$,
\[
  (\mathrm B)=p(\sigma^2-\hat\sigma^2)+\frac{p^2}{n}\cdot\frac1\gamma\bigl(\sigma^2(\sigma^2+\gamma)-\hat\sigma^2(\hat\sigma^2+\gamma)\bigr).
\]

By Step~1, $\hat\sigma^2-\sigma^2=\mathrm O_{\mathbb P}(1/\sqrt{pn})$, so
$p(\sigma^2-\hat\sigma^2)=\mathrm O_{\mathbb P}(\sqrt{p/n})\xrightarrow{\mathbb P}0$ since
$p/n\to0$. For the second term, the identity $a(a+\gamma)-b(b+\gamma)=(a-b)(a+b+\gamma)$
gives
$\hat\sigma^2(\hat\sigma^2+\gamma)-\sigma^2(\sigma^2+\gamma)=(\hat\sigma^2-\sigma^2)(\hat\sigma^2+\sigma^2+\gamma)=\mathrm O_{\mathbb P}(1/\sqrt{pn})=\mathrm o_{\mathbb P}(1)$,
and $p^2/n\to\tau<\infty$ is $\mathrm O(1)$, so the second term is
$\mathrm O(1)\cdot\mathrm o_{\mathbb P}(1)\xrightarrow{\mathbb P}0$. Thus
$(\mathrm B)\xrightarrow{\mathbb P}0$, and by Slutsky
$T_\ell\xrightarrow{\ \mathrm{law}\ }L(\gamma)$.

By Step~5, $\hat q_{1-a}\xrightarrow{\mathbb P}q_{1-a}$, a constant, so Slutsky gives
$T_\ell-\hat q_{1-a}\xrightarrow{\ \mathrm{law}\ }L(\gamma)-q_{1-a}$. Now
$\mathbb P(T_\ell>\hat q_{1-a})=\mathbb E[h(T_\ell-\hat q_{1-a})]$ with
$h(x)=\mathbf 1_{\{x>0\}}$, which is bounded and continuous except at $0$; and the limit
$L(\gamma)-q_{1-a}$ has no atom at $0$ because $L(\gamma)$ has a density (part (A) of the
lemma). The Portmanteau theorem then gives
\[
  \mathbb P(T_\ell>\hat q_{1-a})\to\mathbb P(L(\gamma)>q_{1-a})=a,
\]
by the definition of $q_{1-a}$ and the continuity of $F_{L(\gamma)}$.

\paragraph{Step 7: consistency.}
Under $H_1$, $\theta_\ell>\gamma$ is a fixed positive gap and $p/n\to0$. By
Theorem~\ref{t:main}, $\widehat\lambda_{p,\ell}\xrightarrow{\mathbb P}\sigma^2+\theta_\ell$,
and $\hat\phi(\gamma)\xrightarrow{\mathbb P}\sigma^2+\gamma$, so
$\widehat\lambda_{p,\ell}-\hat\phi(\gamma)\xrightarrow{\mathbb P}\theta_\ell-\gamma>0$, and
multiplying by $p\to\infty$ gives $T_\ell\xrightarrow{\mathbb P}+\infty$. Moreover
$\hat q_{1-a}=\mathrm O_{\mathbb P}(1)$, since Step~5 gives
$\hat q_{1-a}\xrightarrow{\mathbb P}q_{1-a}$ using only $\hat\beta_j\to\beta_j$,
$\widehat{r+s}=r+s$ and $\hat\sigma^2\to\sigma^2$, which hold under $H_1$ as well. For any
fixed $\xi>0$,
$\mathbb P(T_\ell>\hat q_{1-a})\ge\mathbb P(T_\ell>\xi)-\mathbb P(\hat q_{1-a}>\xi)$; the
first term tends to $1$ for every $\xi$, and the second can be made arbitrarily small by taking $\xi$ large; hence $\mathbb P(T_\ell>\hat q_{1-a})\to1$.

\qed

\subsubsection{Proofs of the properties of the test when relaxing assumption (A6)}

\label{subsubsection:Proofs of the properties of the test when relaxing assumption (A6}

All constructions rely on a single fact: the centered eigenvalues are consistent
estimators of the spikes. We isolate it as a lemma.

\begin{lem}[Consistency of the centered eigenvalues]
\label{lem:centered-consistency}
Under assumptions (A1)--(A4), the event $\{\widehat{r+s}=r+s\}$ has probability tending to
$1$, and on this event
\[
  \hat\mu_\ell=\widehat\lambda_{p,\ell}-\hat\sigma^2\ \xrightarrow{\mathbb P}\ \theta_\ell
  \qquad\text{for every }\ell\in\{1,\dots,r+s\}.
\]
\end{lem}

\begin{proof}
That $\mathbb P(\widehat{r+s}=r+s)\to1$ is Step~2 of the proof of Theorem~\ref{thm:test}. We
work on this event, so that $\hat\sigma^2$ is the average of the $p-(r+s)$ smallest
eigenvalues.

Fix $\ell\in\{1,\dots,r+s\}$ and let $\lambda_{p,\ell}$ be the $\ell$-th largest eigenvalue
of the population covariance $\Sigma_p$. We decompose the error into three parts by adding
and subtracting $\lambda_{p,\ell}$ and $\sigma^2$:
\[
  \hat\mu_\ell-\theta_\ell
  =\underbrace{\bigl(\widehat\lambda_{p,\ell}-\lambda_{p,\ell}\bigr)}_{(\mathrm i)}
  +\underbrace{\bigl((\lambda_{p,\ell}-\sigma^2)-\theta_\ell\bigr)}_{(\mathrm{ii})}
  +\underbrace{\bigl(\sigma^2-\hat\sigma^2\bigr)}_{(\mathrm{iii})}.
\]
We bound each term separately.

\emph{Term (i):} By Weyl's inequality (Lemma~\ref{a:weyl}),
$|\widehat\lambda_{p,\ell}-\lambda_{p,\ell}|\le\|\widehat\Sigma_p-\Sigma_p\|_{\mathrm{op}}$.
By the operator-norm bound (Corollary~\ref{a:wishartop}) and since $p/n\to0$ in this regime,
$\|\widehat\Sigma_p-\Sigma_p\|_{\mathrm{op}}=\mathrm O_{\mathbb P}(\sqrt{p/n})$. Hence term
(i) is $\mathrm O_{\mathbb P}(\sqrt{p/n})\to0$.

\emph{Term (ii):} The number $\lambda_{p,\ell}-\sigma^2$ is the
$\ell$-th largest eigenvalue of $A_p+O_pB_pO_p^*$, and by Proposition~\ref{p:pop},
$p\bigl((\lambda_{p,\ell}-\sigma^2)-\theta_\ell\bigr)$ converges in law to $W_\ell$.
A sequence that converges in law after multiplication by $p$ is $\mathrm O_{\mathbb P}(1/p)$,
so term (ii) is $\mathrm O_{\mathbb P}(1/p)\to0$.

\emph{Term (iii):} By Step~1 of the proof of
Theorem~\ref{thm:test}, $\hat\sigma^2-\sigma^2=\mathrm O_{\mathbb P}(1/\sqrt{pn})$, so term
(iii) is $\mathrm O_{\mathbb P}(1/\sqrt{pn})\to0$.

The three terms tend to $0$ in probability, hence so does their sum, and
$\hat\mu_\ell\xrightarrow{\mathbb P}\theta_\ell$.
\end{proof}

We now use this lemma to prove the properties of our test with prior knowledge of the scale of the signals and under known scale separation.

\begin{proof}[Proof of Proposition \ref{prop:calibrated}]
The proof has two parts: first we show that the estimated ranks are correct with probability
tending to $1$; then we transfer the guarantees of Theorem~\ref{thm:test} by a coupling with
an oracle that knows the true ranks.

\medskip
\emph{Part 1: the estimated ranks are correct.} Fix a signal index $i\in\{1,\dots,r\}$.
Assumption (C1) says that $\alpha_i$ is strictly closer to $\alpha_i^\circ$ than any other
spike, so the separation margin
\[
  \rho_i:=\min_{\theta\in\mathcal S,\ \theta\ne\alpha_i}|\alpha_i^\circ-\theta|
  -|\alpha_i^\circ-\alpha_i|
\]
is a minimum over finitely many terms, each strictly larger than $|\alpha_i^\circ-\alpha_i|$,
hence $\rho_i>0$. By definition of $\rho_i$, every spike $\theta_\ell\ne\alpha_i$ (that is,
every $\ell\ne k_i$) satisfies
\begin{equation}
  |\alpha_i^\circ-\theta_\ell|\ \ge\ |\alpha_i^\circ-\alpha_i|+\rho_i.
  \label{eq:margin}
\end{equation}

Work on $\{\widehat{r+s}=r+s\}$, where there are exactly $r+s$ centered eigenvalues. For
each $\ell$ set
\[
  \varepsilon_\ell:=|\hat\mu_\ell-\alpha_i^\circ|-|\theta_\ell-\alpha_i^\circ|.
\]
By the reverse triangle inequality $\bigl||u|-|v|\bigr|\le|u-v|$, applied with
$u=\hat\mu_\ell-\alpha_i^\circ$ and $v=\theta_\ell-\alpha_i^\circ$, we get
$|\varepsilon_\ell|\le|\hat\mu_\ell-\theta_\ell|$, which tends to $0$ in probability by
Lemma~\ref{lem:centered-consistency}. As there are $r+s$ indices (a fixed finite number),
$\max_\ell|\varepsilon_\ell|\xrightarrow{\mathbb P}0$.

Now compare, for any $\ell\ne k_i$, the distance of $\hat\mu_\ell$ to $\alpha_i^\circ$ with
that of $\hat\mu_{k_i}$. Writing each distance as its population value plus $\varepsilon$,
and using $\theta_{k_i}=\alpha_i$,
\[
  |\hat\mu_\ell-\alpha_i^\circ|-|\hat\mu_{k_i}-\alpha_i^\circ|
  =\bigl(|\theta_\ell-\alpha_i^\circ|-|\alpha_i-\alpha_i^\circ|\bigr)
  +(\varepsilon_\ell-\varepsilon_{k_i})
  \ \overset{\eqref{eq:margin}}{\ge}\ \rho_i-|\varepsilon_\ell|-|\varepsilon_{k_i}|.
\]
Consider the event $E_p^{(i)}:=\bigl\{\max_\ell|\varepsilon_\ell|<\rho_i/2\bigr\}$. On
$E_p^{(i)}$ the right-hand side is $>\rho_i-\rho_i/2-\rho_i/2=0$, so $\hat\mu_{k_i}$ is the
strict unique minimizer of $\ell\mapsto|\hat\mu_\ell-\alpha_i^\circ|$, that is $\hat k_i=k_i$.
Since $\max_\ell|\varepsilon_\ell|\xrightarrow{\mathbb P}0$, we have
$\mathbb P(E_p^{(i)})\to1$.

Intersecting over the $r$ signals, the event
$\{\widehat{r+s}=r+s\}\cap\bigcap_{i=1}^r E_p^{(i)}$ is a finite intersection of events each
of probability tending to $1$, hence has probability tending to $1$. On it, $\hat k_i=k_i$
for every $i$. The true ranks $k_1,\dots,k_r$ are distinct, because the spikes are distinct
by (A2); therefore $\{\hat k_1,\dots,\hat k_r\}=\{k_1,\dots,k_r\}=\mathcal K_{\mathrm{sig}}$,
which proves Part~1.

\medskip
\emph{Part 2: coupling with the oracle.} Consider two decision rules run on the same data
and the same Monte-Carlo draws:
\begin{itemize}
\item $\psi_p$, our test, which uses the estimated signal ranks $\{\hat k_1,\dots,\hat k_r\}$
to decide which centered eigenvalues become the batch estimates $\hat\beta_1,\dots,\hat\beta_{\hat s}$;
\item $\psi_p^\star$, an oracle test, identical in every respect but using the true signal
ranks $\{k_1,\dots,k_r\}$.
\end{itemize}
The oracle $\psi_p^\star$ is exactly the test of Theorem~\ref{thm:test}, which assumes the
ranks known; so under $H_0$, $\mathbb P(\psi_p^\star=1)\to a$, and under $H_1$,
$\mathbb P(\psi_p^\star=1)\to1$.

Let $R_p:=\{\hat k_i=k_i\ \forall i\}\cap\{\widehat{r+s}=r+s\}$. On $R_p$, both rules keep the
same centered eigenvalues as batch estimates, so they build the same $\hat\beta_j$, the same
estimated law $\widehat L$, draw from it the same $M$ copies (the draws are shared), obtain
the same quantile $\hat q_{1-a}$, and compare it to the same statistic $T_\ell$; hence
$\psi_p=\psi_p^\star$ on $R_p$. Therefore $\{\psi_p\ne\psi_p^\star\}\subseteq R_p^c$, and
\[
  \bigl|\mathbb P(\psi_p=1)-\mathbb P(\psi_p^\star=1)\bigr|
  \le\mathbb P(\psi_p\ne\psi_p^\star)
  \le\mathbb P(R_p^c)\ \longrightarrow\ 0
\]
by Part~1. Combining with the limits of $\mathbb P(\psi_p^\star=1)$ gives
$\mathbb P(\psi_p=1)\to a$ under $H_0$ and $\mathbb P(\psi_p=1)\to1$ under $H_1$.
\end{proof}

The proof technique is similar when we assume a known scale separation between batch spikes and signal spikes.

\begin{proof}[Proof of Proposition \ref{prop:scale-separation}]
Under (S1) every signal exceeds $\eta$ and every batch effect is below $\eta$, so the $r$
signals are the $r$ largest spikes. Since the spikes are sorted decreasingly, their ranks
are $1,\dots,r$, that is $\mathcal K_{\mathrm{sig}}=\{1,\dots,r\}$.

Let
\[
  \delta:=\min\Bigl(\min_{i\in\{1,\dots,r\}}(\alpha_i-\eta),\
  \min_{j\in\{1,\dots,s\}}(\eta-\beta_j)\Bigr),
\]
a minimum of finitely many strictly positive numbers by (S1), so $\delta>0$. Then
$\theta_\ell-\eta\ge\delta$ for every signal rank $\ell\le r$, and $\eta-\theta_\ell\ge\delta$
for every batch rank $\ell>r$.

By Lemma~\ref{lem:centered-consistency}, $\hat\mu_\ell-\theta_\ell\xrightarrow{\mathbb P}0$
for each $\ell$. Hence, for a signal rank $\ell\le r$,
\[
  \mathbb P(\hat\mu_\ell\le\eta)
  =\mathbb P\bigl(\hat\mu_\ell-\theta_\ell\le\eta-\theta_\ell\bigr)
  \le\mathbb P\bigl(\hat\mu_\ell-\theta_\ell\le-\delta\bigr)
  \le\mathbb P\bigl(|\hat\mu_\ell-\theta_\ell|\ge\delta\bigr)\ \longrightarrow\ 0,
\]
and symmetrically, for a batch rank $\ell>r$,
$\mathbb P(\hat\mu_\ell>\eta)\le\mathbb P(|\hat\mu_\ell-\theta_\ell|\ge\delta)\to0$.

Consider the event
\[
  R_p:=\{\widehat{r+s}=r+s\}\cap\bigcap_{\ell\le r}\{\hat\mu_\ell>\eta\}
  \cap\bigcap_{\ell>r}\{\hat\mu_\ell\le\eta\}.
\]
By a union bound over the $r+s$ complementary events, each of probability tending to $0$,
and since $\mathbb P(\widehat{r+s}=r+s)\to1$, we get $\mathbb P(R_p)\to1$. On $R_p$, the
centered eigenvalues above $\eta$ are exactly those of the signals, so
$\widehat{\mathcal K}_{\mathrm{sig}}=\{1,\dots,r\}=\mathcal K_{\mathrm{sig}}$.

For level and consistency we argue by coupling with an oracle, as in the proof of
Proposition~\ref{prop:calibrated}. Let $\psi_p$ be the present test, which uses the
estimated signal ranks $\widehat{\mathcal K}_{\mathrm{sig}}$, and let $\psi_p^\star$ be the
oracle test of Theorem~\ref{thm:test}, identical but using the true ranks
$\mathcal K_{\mathrm{sig}}$, run on the same data and Monte-Carlo draws. On the event
$R_p:=\{\widehat{\mathcal K}_{\mathrm{sig}}=\mathcal K_{\mathrm{sig}}\}\cap\{\widehat{r+s}=r+s\}$
the two rules select the same batch estimates, hence build the same $\widehat L$, the same
quantile and the same statistic, so $\psi_p=\psi_p^\star$. Therefore
$\{\psi_p\ne\psi_p^\star\}\subseteq R_p^c$, and
\[
  \bigl|\mathbb P(\psi_p=1)-\mathbb P(\psi_p^\star=1)\bigr|
  \le\mathbb P(\psi_p\ne\psi_p^\star)
  \le\mathbb P(R_p^c)\ \longrightarrow\ 0,
\]
since $\mathbb P(R_p)\to1$ by the first part of the proof. As
$\mathbb P(\psi_p^\star=1)\to a$ under $H_0$ and $\to1$ under $H_1$ by
Theorem~\ref{thm:test}, the same limits hold for $\psi_p$.
\end{proof}

\section{Discussion}
\label{s:disc}

\subsection{What the transition means for estimation}
\label{ss:disc-estim}

Table~\ref{t:four} can be viewed as a statement on the estimation of a spike. The sample
eigenvalue does not estimate $\sigma^2+\theta_\ell$ directly: as in the classical spiked
model, it is centered at $\phi_{p/n}(\theta_\ell)$, which differs from $\sigma^2+\theta_\ell$
by a deterministic bias of order $p/n$. Once this bias is substracted, two errors remain: the sampling error, of order $n^{-1/2}$, and the Haar
fluctuation, of order $1/p$. The size of the total error
$\widehat\lambda_{p,\ell}-\phi_{p/n}(\theta_\ell)$ is that of the larger of the two, and
comparing $n^{-1/2}$ with $1/p$ amounts to comparing $n$ with $p^2$. This is why the
correct normalization is $\sqrt n$ when $p^2/n\to\tau\in(0,\infty]$, as in the classical
spiked model, and $p$ when $p^2/n\to0$. The last case has no counterpart in the classical
spiked model, where the population matrix is fixed and where the sampling error is the only
one.

The shape of the limiting law is also important. In the two regimes where the Haar
fluctuation is not negligible compared to the sampling fluctuation, that is $p^2/n\to\tau\in[0,\infty)$, the limit contains the
variable $W_\ell$, which is neither Gaussian nor centered. Hence, a confidence interval for the value of a spike $\theta_\ell$ built from
a Gaussian approximation would then be misspecified. Our results give the
law to be used instead, and Algorithm~\ref{alg:spectral_test} shows how to simulate it.

\subsection{Limits of the model and open questions}
\label{ss:disc-open}

We list the assumptions that we would like to remove.

\emph{Non-Gaussian samples.}
The results of Bai and Yao \cite{bai08,bai12} cover non-Gaussian spiked models under suitable moment and independence assumptions.
Applying such results to the present random population model, however, requires additional work.
Lemma~\ref{l:cond} uses Gaussian rotational invariance, which generally fails when the independent Gaussian entries are replaced by non-Gaussian entries in a fixed coordinate system.
Moreover, the sampling fluctuation variance may depend on fourth cumulants and on the population eigenvectors; in general it is not obtained by simply replacing $\kappa$ in the Gaussian formula.
The population fluctuation $W_\ell$ is unaffected as long as the population model is unchanged.
We leave these extensions to future work.

\emph{Fixed ranks.} The ranks $r$ and $s$ are fixed, and our results do not hold anymore when letting them grow with $p$. Proposition~\ref{p:pop} treats the $r\times s$ block of overlaps
$a_i^*O_pb_j$ as a small perturbation, and this block is small only because it has a fixed
number of entries, each of size $p^{-1/2}$: if $r$ and $s$ grow too fast, it is no longer
small and Proposition~\ref{p:pop} does not hold anymore. Another
difficulty is that the strengths $\alpha_1,\ldots,\alpha_r$ and $\beta_1,\ldots,\beta_s$
are assumed to be distinct, while the weights $\alpha_i\beta_j/(\alpha_i-\beta_j)$
appearing in $W_\ell$ grow to infinity when two of them become close: a growing number of strengths in
a bounded range would have to be accompanied by an assumption on how they separate. We do
not know which limiting law to expect in that setting.

\emph{Beyond Haar rotations.}
Exact Haar invariance can be weakened by imposing conditions directly on the overlaps. Let $Q_p$ be an orthogonal or unitary random matrix, and let $a_1,\ldots,a_r$ and $b_1,\ldots,b_s$ be orthonormal eigenvectors associated with the nonzero eigenvalues of $A_p$ and $B_p$.
Set
\[
  H_p=(a_i^*Q_pb_j)_{1\le i\le r,\,1\le j\le s}.
\]
Suppose that $H_p\to0$ almost surely and that $\sqrt p\,H_p$ converges jointly in distribution to a matrix of independent standard real or circular complex Gaussian entries.
Then the same finite-rank perturbation calculation gives the population limits in Proposition~\ref{p:pop}, with the same random variables $W_\ell$.
If the conditional Gaussian sampling model is retained, with the standard Gaussian noise independent of $Q_p$, Lemma~\ref{l:cond} and the sampling arguments remain valid.
Thus the asymptotic conclusions of Theorems~\ref{t:main}, \ref{t:infy} and~\ref{t:prop} are unchanged.
For the population fluctuation limit alone, the joint weak convergence of $\sqrt p\,H_p$ is sufficient; the almost-sure condition is used to retain the almost-sure population convergence.
Different overlap scales or dependence in their limiting law may change the population fluctuation and its competition with sampling.

\emph{Eigenvectors.} We only study eigenvalues. The sample eigenvector associated with
$\widehat\lambda_{p,\ell}$ carries the same two sources of the randomness, coming from the sampling and the Haar matrix, so one may ask whether its behaviour is also governed by $p^2/n$, as for sample eigenvalues. We did not treat
this question and we keep it for further work.

\emph{The test.} Algorithm~\ref{alg:spectral_test} has two limitations. The first one is
assumption (A6): the reference law $W_\ell$ is a sum over the batch strengths only, so one
has to know which of the estimated spikes are signals and which are batch effects.
Section~\ref{sss:relax assumption A6} weakens this assumption but does not remove it, and we
do not know how to decide, from the sample eigenvalues alone, whether a given spike comes
from $A_p$ or from $O_pB_pO_p^*$. The second one is that the test assumes we are in the regime where $p^2/n \to \tau < \infty$. In reality, a practitioner often only has access to $n,p$ and not always to the rate at which $p^2$ grows compared to $n$. A procedure which does not require this
knowledge, and which would adapt to the regime (whether $p^2/n \to \infty$ or $p^2/n \to \tau <\infty$), would be more useful in practice.

\newpage
\appendix
\section{Appendix: frequently used facts}\label{s:appfacts}

This appendix collects standard facts that are used several times in the paper.  They are stated in the form needed in the main text.
\subsection{A fixed interval for a separated spike}\label{s:appsep}

\begin{lem}[A fixed interval around a separated spike]\label{a:sepint}
Let \(K\) be fixed. Suppose that
\[
  \lambda_p\to\lambda_*>\sigma^2,
\]
and that, for some \(\delta_0>0\) and all large \(p\),
\[
  \lambda_p-\sigma^2\ge\delta_0,
  \qquad
  |\lambda_p-\lambda_{p,j}|\ge\delta_0
  \quad(2\le j\le K).
\]
Set
\[
  \nu=\frac1{16}\min\{\delta_0,\lambda_* -\sigma^2\},
\]
\[
  J=[\lambda_*-\nu,\lambda_*+\nu],
  \qquad
  I=[\lambda_*-2\nu,\lambda_*+2\nu].
\]
Then, for all large \(p\),
\[
  \lambda_p\in J\subset\operatorname{int}I,
  \qquad
  \operatorname{dist}(\lambda_p,\partial I)\ge\nu,
\]
\[
  \operatorname{dist}\!\left(
    I,\{\sigma^2\}\cup\{\lambda_{p,j}:2\le j\le K\}
  \right)\ge4\nu.
\]
Moreover, there is a fixed \(\varepsilon>0\) such that
\[
  I\subset(\sigma^2+\varepsilon,\infty).
\]
\end{lem}

\begin{proof}
Since \(\lambda_p\to\lambda_*\), for all large \(p\),
\(|\lambda_p-\lambda_*|\le\nu\). This gives the first line.
For \(2\le j\le K\),
\[
  |\lambda_* -\lambda_{p,j}|
  \ge |\lambda_p-\lambda_{p,j}|-|\lambda_p-\lambda_*|
  \ge\delta_0-\nu.
\]
Every point of \(I\) is at distance at most \(2\nu\) from \(\lambda_*\). Hence
\[
  \operatorname{dist}(I,\lambda_{p,j})
  \ge\delta_0-3\nu\ge4\nu.
\]
The same argument with \(\sigma^2\) gives
\(\operatorname{dist}(I,\sigma^2)\ge4\nu\).
Finally, if \(a=\lambda_* -\sigma^2\), then
\(2\nu\le a/8\). Thus the left end of \(I\) is at least
\(\sigma^2+7a/8\), so one may take, for example,
\(\varepsilon=a/2\).
\end{proof}

\subsection{Gaussian sample covariance bounds}\label{s:appgauss}

\begin{prop}[Koltchinskii--Lounici bound in the Gaussian case]\label{a:klcov}
Let \(X_1,\ldots,X_n\) be independent centered Gaussian vectors in \(\mathbb R^m\), with covariance matrix \(\Gamma\).  Put
\[
  \widehat\Gamma=\frac1n\sum_{k=1}^n X_kX_k^*.
\]
There is a universal constant \(C>0\) such that, for every \(u\ge1\), with probability at least \(1-e^{-u}\),
\[
  \|\widehat\Gamma-\Gamma\|_{\mathrm{op}}
  \le C\|\Gamma\|_{\mathrm{op}}
  \left(
    \sqrt{\frac{r(\Gamma)}{n}}
    +\frac{r(\Gamma)}{n}
    +\sqrt{\frac{u}{n}}
    +\frac{u}{n}
  \right),
\]
where the effective rank is
\[
  r(\Gamma)=\frac{\operatorname{Tr}\Gamma}{\|\Gamma\|_{\mathrm{op}}}
\]
when \(\Gamma\ne0\), and \(r(0)=0\).
The same conclusion holds in the complex circular Gaussian case, up to a change of the universal constant.
\end{prop}

\begin{proof}
This is the Gaussian finite-dimensional form of Koltchinskii--Lounici \cite[Theorem~4 and Corollary~2]{kol17}.  In the complex case, one applies the real result to the realification of the complex vector; this only changes absolute constants.
\end{proof}

\begin{cor}[Operator-norm estimates used in the paper]\label{a:wishartop}
Let \(Z_p\) be a \(p\times n\) standard real or complex Gaussian matrix.  If \(p/n\to0\), then
\[
  \left\|\frac1nZ_pZ_p^*-I_p\right\|_{\mathrm{op}}
  =\mathrm O_{\mathbb P}\left(\sqrt{\frac pn}+\frac pn\right)
  =\mathrm O_{\mathbb P}\left(\sqrt{\frac pn}\right).
\]
More generally, if \(D_p\) is deterministic positive semidefinite, \(\sup_p\|D_p\|_{\mathrm{op}}<\infty\), and \(\operatorname{rank}D_p\le Cp\), then
\[
  \left\|D_p^{1/2}\left(\frac1nZ_pZ_p^*\right)D_p^{1/2}-D_p\right\|_{\mathrm{op}}
  =\mathrm O_{\mathbb P}\left(\sqrt{\frac pn}+\frac pn\right).
\]
In particular, for \(E_p=\Lambda_p^{1/2}(n^{-1}Z_pZ_p^*)\Lambda_p^{1/2}-\Lambda_p\) with \(\sup_p\|\Lambda_p\|_{\mathrm{op}}<\infty\),
\[
  \|E_p\|_{\mathrm{op}}
  =\mathrm O_{\mathbb P}\left(\sqrt{\frac pn}\right).
\]
\end{cor}

\begin{proof}
Apply Proposition~\ref{a:klcov}.  For \(D_p\), the covariance is \(D_p\), and its effective rank is bounded by a constant times \(p\) under the stated assumptions.
\end{proof}

\begin{rem}[The case of a random $\Lambda_p$]\label{a:wishartrand}
In the paper, Corollary~\ref{a:wishartop} is used with $\Lambda_p$ random and
$\mathcal F_p$-measurable, while $Z_p$ is independent of $\mathcal F_p$. The conclusion
still holds, for the following reason. In our model
\[
  \|\Lambda_p\|_{\mathrm{op}}=\|\Sigma_p\|_{\mathrm{op}}
  \le M_0:=\sigma^2+\max_{i\le r}\alpha_i+\max_{j\le s}\beta_j ,
\]
for every $p$ and every realization of $O_p$, so the bound is deterministic.
Conditionally on $\mathcal F_p$, the matrix $\Lambda_p$ is deterministic and $Z_p$ is
still a standard Gaussian matrix, so Corollary~\ref{a:wishartop} applies conditionally.
The constant it produces depends only on $M_0$ and on the universal constant of
Proposition~\ref{a:klcov}, and not on the realization of $O_p$. Hence, for every
$\varepsilon>0$ there is a constant $C_\varepsilon$, not depending on $p$ nor on the
realization, such that for all large $p$
\[
  \mathbb P\Bigl(\|E_p\|_{\mathrm{op}}>C_\varepsilon\sqrt{p/n}
  \ \Big|\ \mathcal F_p\Bigr)\le\varepsilon .
\]
The right-hand side is a constant, so taking expectations gives the same bound
unconditionally, that is $\|E_p\|_{\mathrm{op}}=\mathrm O_{\mathbb P}(\sqrt{p/n})$.
\end{rem}

\begin{lem}[A chi-square tail consequence]\label{a:chisq}
Let \(x_1,\ldots,x_n\) be independent \(N_{\mathbb F}(0,\sigma^2I_m)\) vectors, real or complex.  If \(m/n\to0\), then
\[
  \max_{k\le n}\frac{\|x_k\|^2}{n}=o_{\mathbb P}(1).
\]
\end{lem}

\begin{proof}
In the real case, \(\|x_k\|^2/\sigma^2\sim\chi_m^2\).  Laurent--Massart's bound \cite[Lemma~1]{laurent00} gives
\[
  \mathbb P(\chi_m^2\ge m+2\sqrt{mu}+2u)\le e^{-u}.
\]
Taking \(u=2\log n\) and using a union bound over \(k\le n\) proves the result, because \((m+\log n)/n\to0\).  The complex case is the same after writing a complex standard Gaussian vector as two real Gaussian vectors with variance \(1/2\).
\end{proof}

\subsection{Marchenko--Pastur facts}\label{s:appmp}

For \(c>0\), let
\[
  a_c=\sigma^2(1-\sqrt c)^2,
  \qquad
  b_c=\sigma^2(1+\sqrt c)^2.
\]
The Marchenko--Pastur distribution with aspect ratio \(c\) and variance level \(\sigma^2\) is
\[
\begin{split}
  \mu_{c,\sigma^2}
  ={}&\max\left\{1-\frac1c,0\right\}\delta_0 \\
  &+\frac{\sqrt{(b_c-x)(x-a_c)}}{2\pi c\sigma^2x}
  \mathbf1_{[a_c,b_c]}(x)\,dx.
\end{split}
\]
When \(c=0\), we use the convention \(\mu_{0,\sigma^2}=\delta_{\sigma^2}\).  In the main proof we only use \(c\to0\), so eventually \(0<c<1\), and there is no atom at zero.

\begin{thm}[Marchenko--Pastur empirical law]\label{a:mpesd}
Let \(Z\) be an \(m\times n\) standard real or complex Gaussian matrix and let
\[
  C_{m,n}=\frac{\sigma^2}{n}ZZ^*.
\]
If \(m/n\to c\in[0,\infty)\), then the empirical eigenvalue distribution of \(C_{m,n}\) converges weakly, in probability and almost surely, to \(\mu_{c,\sigma^2}\).
\end{thm}

\begin{proof}
This is the Marchenko--Pastur theorem \cite{mp67}; see also Bai--Silverstein \cite[Chapter~3]{bai10} or Yao--Zheng--Bai \cite[Chapter~2]{yao15}.
\end{proof}

For \(0<c<1\) and \(t>b_c\), define the positive-resolvent transforms
\begin{equation}
  s_c(t)=c\int\frac1{t-x}\,\mu_{c,\sigma^2}(dx),
  \qquad
  h_c(t)=c\int\frac{x}{t-x}\,\mu_{c,\sigma^2}(dx).
  \label{e:appsh}
\end{equation}
For \(c=0\), set \(s_0(t)=h_0(t)=0\).  Since \(x/(t-x)=t/(t-x)-1\),
\begin{equation}
  h_c(t)=t s_c(t)-c.
  \label{e:apphts}
\end{equation}

\begin{prop}[Marchenko--Pastur self-consistent equation]\label{a:mpself}
Let \(0<c<1\) and \(t>b_c\).  Then \(s=s_c(t)\) satisfies
\begin{equation}
  G_t(s):=\sigma^2t s^2-\{t-\sigma^2(1-c)\}s+c=0.
  \label{e:appG}
\end{equation}
Equivalently, for \(s\neq1/\sigma^2\),
\begin{equation}
  F_t(s):=ts-c-\frac{\sigma^2s}{1-\sigma^2s}=0
  \label{e:appF}
\end{equation}
if and only if \(G_t(s)=0\).  More precisely,
\[
  (1-\sigma^2s)F_t(s)=-G_t(s).
\]
\end{prop}

\begin{proof}
Let \(m_{c,\sigma^2}(z)=\int(x-z)^{-1}\mu_{c,\sigma^2}(dx)\) be the usual Stieltjes transform.  It satisfies
\[
  c\sigma^2 z m_{c,\sigma^2}(z)^2
  +\{z-\sigma^2(1-c)\}m_{c,\sigma^2}(z)+1=0
\]
for \(z\in\mathbb C\setminus[a_c,b_c]\); see Bai--Silverstein \cite[Chapter~3]{bai10}.  For real \(t>b_c\), \(s_c(t)=-c m_{c,\sigma^2}(t)\).  Substitution gives \eqref{e:appG}.  The identity relating \(F_t\) and \(G_t\) follows by direct algebra.
\end{proof}

\begin{prop}[Root identification]\label{a:mproot}
Let \(0<c<1\) and \(t>b_c\).  The equation \(G_t(s)=0\) has two distinct real roots
\[
  0<s_-(t)<s_+(t)<\frac1{\sigma^2}.
\]
Moreover,
\[
  s_c(t)=s_-(t).
\]
Equivalently, among the roots of \(F_t(s)=0\) in \((0,1/\sigma^2)\), the Marchenko--Pastur value \(s_c(t)\) is the one for which \(F_t'(s)>0\).
\end{prop}

\begin{proof}
The discriminant of \(G_t\) is
\[
  \{t-\sigma^2(1-c)\}^2-4\sigma^2ct
  =(t-a_c)(t-b_c)>0.
\]
Thus the two roots are real and distinct.  Their sum and product are positive, so both roots are positive.  Also
\[
  G_t(1/\sigma^2)=1>0,
\]
and the midpoint of the two roots is smaller than \(1/(2\sigma^2)\).  Hence both roots lie in \((0,1/\sigma^2)\).

On \((-\infty,1/\sigma^2)\),
\[
  F_t''(s)=-\frac{2\sigma^4}{(1-\sigma^2s)^3}<0.
\]
So \(F_t\) is strictly concave.  Since \(F_t(0)=-c<0\) and the two roots are \(s_-(t)\) and \(s_+(t)\), strict concavity implies
\[
  F_t'(s_-(t))>0,
  \qquad
  F_t'(s_+(t))<0.
\]
Finally, \(t\mapsto s_c(t)\) is continuous on \((b_c,\infty)\), and it must equal one of these two continuous roots.  For large \(t\),
\[
  s_c(t)\le\frac{c}{t-b_c}\to 0,
\]
whereas \(s_+(t)\) stays bounded away from zero.  Therefore \(s_c(t)=s_-(t)\) for large \(t\), and by continuity it holds for every \(t>b_c\).
\end{proof}

\subsection{Gaussian quadratic forms}\label{s:appquad}

\begin{lem}[Gaussian quadratic form identities]\label{a:qform}
Let \(z=(z_1,\ldots,z_n)\) be a standard Gaussian row vector, real or complex circular.  Let \(A=A^*\in M_n(\mathbb F)\) be deterministic, or independent of \(z\).  Then, conditionally on \(A\),
\[
  \mathbb E_z[zAz^*]=\operatorname{Tr}A,
\]
and
\[
  \operatorname{Var}_z(zAz^*)=\kappa\operatorname{Tr}(A^2),
  \qquad
  \kappa=
  \begin{cases}
  2,&\mathbb F=\mathbb R,\\
  1,&\mathbb F=\mathbb C.
  \end{cases}
\]
Consequently,
\[
  \operatorname{Var}_z(zAz^*)
  \le \kappa\operatorname{rank}(A)\|A\|_{\mathrm{op}}^2.
\]
If \(x\sim N_{\mathbb F}(0,\sigma^2I_n)\), the same formulas hold with the right-hand sides multiplied by \(\sigma^2\) for the mean and \(\sigma^4\) for the variance.
\end{lem}

\begin{proof}
Diagonalize \(A=U\operatorname{diag}(a_1,\ldots,a_n)U^*\).  By orthogonal or unitary invariance of the standard Gaussian vector, \(zU\) has the same distribution as \(z\).  Thus
\[
  zAz^*\overset{d}{=}\sum_{j=1}^n a_j|z_j|^2.
\]
The mean is \(\sum_j a_j=\operatorname{Tr}A\).  In the real case, \(|z_j|^2\sim\chi_1^2\), whose variance is \(2\).  In the complex circular case, \(|z_j|^2\sim\operatorname{Exp}(1)\), whose variance is \(1\).  This proves the variance formula.  The final bound follows from \(\operatorname{Tr}(A^2)=\sum_j a_j^2\le\operatorname{rank}(A)\|A\|_{\mathrm{op}}^2\).
\end{proof}

\subsection{Linear algebra facts}\label{s:applinear}

We use the following elementary facts without further comment.
\begin{lem}[Weyl's eigenvalue inequality]\label{a:weyl}
Let \(A\) and \(B\) be Hermitian \(d\times d\) matrices, and write their eigenvalues in decreasing order. Then
\[
  \max_{1\le j\le d}|\mu_j(A)-\mu_j(B)|
  \le\|A-B\|_{\mathrm{op}}.
\]
\end{lem}

\begin{lem}[Schur complement eigenvalue equation]\label{a:schur}
Let
\[
  S=\begin{pmatrix}a&b^*\\ b&C\end{pmatrix},
\]
where \(a\) is a scalar and \(C\) is Hermitian. If \(tI-C\) is invertible, then
\[
  \det(tI-S)
  =\det(tI-C)\{t-a-b^*(tI-C)^{-1}b\}.
\]
In particular, if \(t\notin\operatorname{Spec}(C)\), then \(t\) is an eigenvalue of \(S\) if and only if
\[
  t=a+b^*(tI-C)^{-1}b.
\]
\end{lem}

\begin{proof}
Multiply \(tI-S\) on the left by
\[
  \begin{pmatrix}
    1&b^*(tI-C)^{-1}\\
    0&I
  \end{pmatrix}.
\]
Its determinant is one, and the product is block lower triangular with diagonal blocks
\(t-a-b^*(tI-C)^{-1}b\) and \(tI-C\). Taking determinants gives the formula.
\end{proof}

\begin{lem}[Sherman--Morrison formula]\label{a:sm}
Let \(M\) be invertible and let \(w\) be a vector.  If \(1-w^*M^{-1}w\neq0\), then
\[
  (M-ww^*)^{-1}
  =M^{-1}+\frac{M^{-1}ww^*M^{-1}}{1-w^*M^{-1}w}.
\]
\end{lem}

\begin{lem}[Resolvent and finite-rank trace bounds]\label{a:resrank}
If \(A\) and \(B\) are Hermitian and \(tI-A\), \(tI-B\) are invertible, then
\[
  (tI-A)^{-1}-(tI-B)^{-1}
  =(tI-A)^{-1}(A-B)(tI-B)^{-1}.
\]
Also, if \(M\) has rank at most \(r\), then
\[
  |\operatorname{Tr}M|\le r\|M\|_{\mathrm{op}}.
\]
\end{lem}
\begin{lem}[Finite-rank resolvent trace comparison]\label{a:frtrace}
Let \(A\) and \(B\) be Hermitian matrices of the same size. Let \(I\subset\mathbb R\) be compact, and assume that, for some \(d>0\),
\[
  \operatorname{dist}(I,\operatorname{Spec}(A))\ge d,
  \qquad
  \operatorname{dist}(I,\operatorname{Spec}(B))\ge d.
\]
If
\[
  \operatorname{rank}(A-B)\le r,
  \qquad
  \|A-B\|_{\mathrm{op}}\le M,
\]
then
\[
  \sup_{t\in I}
  \left|
    \operatorname{Tr}\{(tI-A)^{-1}-(tI-B)^{-1}\}
  \right|
  \le \frac{rM}{d^2}.
\]
\end{lem}

\begin{proof}
The resolvent identity gives
\[
  (tI-A)^{-1}-(tI-B)^{-1}
  =(tI-A)^{-1}(A-B)(tI-B)^{-1}.
\]
The matrix on the right has rank at most \(r\) and operator norm at most \(M/d^2\). The result follows from
\(|\operatorname{Tr}H|\le\operatorname{rank}(H)\|H\|_{\mathrm{op}}\).
\end{proof}

\subsection{Convergence lemmas}\label{s:appconv}

\begin{lem}[Localization]\label{l:localization}
Let \(\mathbb P(\Omega_p)\to1\) and let \(u_p>0\).  If
\[
  \mathbb E\left[\mathbf 1_{\Omega_p}|X_p|\right]=o(u_p),
\]
then \(X_p=o_{\mathbb P}(u_p)\).
\end{lem}

\begin{proof}
For every \(\delta>0\), Markov's inequality gives
\[
  \mathbb P(|X_p|>\delta u_p)
  \leq \mathbb P(\Omega_p^c)
  +\frac{\mathbb E[\mathbf 1_{\Omega_p}|X_p|]}{\delta u_p}.
\]
The right-hand side tends to zero.
\end{proof}

\begin{lem}[Joint convergence for independent arrays]\label{l:join}
Let $X_p$ be $\mathcal F_p$-measurable and let $Y_p$ be independent of
$\mathcal F_p$. If $X_p\xrightarrow{\ \mathrm{law}\ } X$ and
$Y_p\xrightarrow{\ \mathrm{law}\ } Y$, then
\[
  (X_p,Y_p)\xrightarrow{\ \mathrm{law}\ }(X,Y^{\circ}),
\]
where $Y^{\circ}$ has the same distribution as $Y$ and is independent of $X$.
\end{lem}

\begin{proof}
For real $t,s$, independence gives
\[
  \mathbb E\exp\bigl(\mathrm{i}tX_p+\mathrm{i}sY_p\bigr)
  =\mathbb E\exp(\mathrm{i}tX_p)\,\mathbb E\exp(\mathrm{i}sY_p).
\]
The right-hand side converges to $\mathbb E e^{\mathrm{i}tX}\,\mathbb E e^{\mathrm{i}sY}$,
which is the characteristic function of an independent pair. The result follows from
L\'{e}vy's continuity theorem.
\end{proof}

\begin{lem}[Deconditioning by characteristic functions]\label{l:deco}
Let $X_p$ be real-valued random variables and let $\mathcal F_p$ be sigma-fields.
Suppose that, for every $t\in\mathbb R$,
\[
  \mathbb E\left[\exp(\mathrm{i}tX_p)\mid\mathcal F_p\right]
  \xrightarrow{\mathbb P} \phi(t),
\]
where $\phi$ is the characteristic function of a probability distribution. Then $X_p$
converges in distribution to that probability distribution.
\end{lem}

\begin{proof}
The conditional characteristic functions have modulus at most one. Convergence in
probability together with boundedness implies convergence of expectations. Hence
\[
  \mathbb E\exp(\mathrm{i}tX_p)
  =\mathbb E\left[\mathbb E\left[\exp(\mathrm{i}tX_p)\mid\mathcal F_p\right]\right]
  \longrightarrow \phi(t)
\]
for every $t$. L\'{e}vy's continuity theorem proves the assertion.
\end{proof}

\begin{lem}[Change of scaling]\label{a:scale}
Let $X_p$ be real random variables, and let $u_p>0$ and $v_p>0$ be deterministic. Assume
that $u_pX_p\xrightarrow{\ \mathrm{law}\ }L$ for some real random variable $L$.
\begin{enumerate}
  \item[\textnormal{(i)}] If $v_p/u_p\to\gamma\in[0,\infty)$, then
  $v_pX_p\xrightarrow{\ \mathrm{law}\ }\gamma L$.
  \item[\textnormal{(ii)}] If $v_p/u_p\to\infty$ and $\mathbb P(L=0)=0$, then
  $|v_pX_p|\xrightarrow{\mathbb P}\infty$.
\end{enumerate}
\end{lem}

\begin{proof}
(i) Write $v_pX_p=(v_p/u_p)(u_pX_p)$. The first factor is a deterministic sequence
converging to $\gamma$, and the second one converges in distribution to $L$. Slutsky's
theorem gives the claim.

(ii) Fix $M>0$ and $\varepsilon>0$. Since $v_p/u_p\to\infty$, we have $u_p/v_p\to0$, so
there is $p_0$ with $Mu_p/v_p\le\varepsilon$ for every $p\ge p_0$. For such $p$, since
$u_p,v_p>0$,
\[
  \bigl\{|v_pX_p|\le M\bigr\}
  =\Bigl\{|u_pX_p|\le M\,\frac{u_p}{v_p}\Bigr\}
  \subseteq\bigl\{u_pX_p\in[-\varepsilon,\varepsilon]\bigr\} .
\]
The set $[-\varepsilon,\varepsilon]$ is closed, so the portmanteau theorem gives
\[
  \limsup_{p\to\infty}\mathbb P\bigl(|v_pX_p|\le M\bigr)
  \le\limsup_{p\to\infty}
  \mathbb P\bigl(u_pX_p\in[-\varepsilon,\varepsilon]\bigr)
  \le\mathbb P\bigl(|L|\le\varepsilon\bigr).
\]
This holds for every $\varepsilon>0$. The events $\{|L|\le1/k\}$ decrease to $\{L=0\}$
as $k\to\infty$, so, by continuity from above,
$\mathbb P(|L|\le1/k)\to\mathbb P(L=0)=0$. Taking $\varepsilon=1/k$ and letting
$k\to\infty$ gives $\limsup_p\mathbb P(|v_pX_p|\le M)=0$.
\end{proof}

\begin{lem}[Divergence in probability excludes weak limits]\label{a:nolimit}
Let $Y_p$ be real random variables with $|Y_p|\xrightarrow{\mathbb P}\infty$. Then no
subsequence of $(Y_p)$ converges in distribution to a real random variable. In
particular, $(Y_p)$ is not tight.
\end{lem}

\begin{proof}
Suppose that $Y_{p_k}\xrightarrow{\ \mathrm{law}\ }Y$ along a subsequence $(p_k)$, for
some real random variable $Y$. The events $\{|Y|<M\}$, $M\in\mathbb N$, increase to the
whole space, so there is $M>0$ with $\mathbb P(|Y|<M)\ge1/2$. The set $(-M,M)$ is open,
so the portmanteau theorem gives
\[
  \liminf_{k\to\infty}\mathbb P\bigl(Y_{p_k}\in(-M,M)\bigr)
  \ge\mathbb P\bigl(Y\in(-M,M)\bigr)\ge\frac12 .
\]
On the other hand $\mathbb P(Y_{p_k}\in(-M,M))\le\mathbb P(|Y_{p_k}|\le M)\to0$ by
assumption, a contradiction. Since every tight sequence has a subsequence converging in
distribution, by Prokhorov's theorem, $(Y_p)$ is not tight.
\end{proof}

\begin{lem}[Convolution with an absolutely continuous law]\label{a:conv}
Let $\xi$ and $\eta$ be independent real random variables, and assume that the law of
$\xi$ is absolutely continuous with respect to the Lebesgue measure, with density $f$.
Then the law of $\xi+\eta$ is absolutely continuous, with density
\[
  h(x)=\int_{\mathbb R}f(x-y)\,\mathbb P_\eta(\mathrm dy),
\]
where $\mathbb P_\eta$ is the law of $\eta$. In particular
$\mathbb P(\xi+\eta=x)=0$ for every $x\in\mathbb R$.
\end{lem}

\begin{proof}
Let $A\subseteq\mathbb R$ be a Borel set. By independence and by Tonelli's theorem, all
the integrands being non-negative and measurable,
\[
  \mathbb P(\xi+\eta\in A)
  =\int_{\mathbb R}\mathbb P(\xi\in A-y)\,\mathbb P_\eta(\mathrm dy)
  =\int_{\mathbb R}\Bigl(\int_Af(x-y)\,\mathrm dx\Bigr)\mathbb P_\eta(\mathrm dy)
  =\int_Ah(x)\,\mathrm dx .
\]
So $h$ is a density of $\xi+\eta$. A single point has Lebesgue measure zero, which gives
the last assertion.
\end{proof}

\begin{cor}[The limit laws have no atom at zero]\label{a:atom}
Fix $\ell\in\{1,\ldots,r+s\}$, let $W_\ell$ be given by \eqref{e:well}, and let
$G_\ell\sim\mathcal N(0,1)$ be independent of $W_\ell$. For every $\tau\in[0,\infty)$,
the law of
\[
  W_\ell+\sqrt{\kappa\tau V_0(\theta_\ell)}\,G_\ell
\]
is absolutely continuous with respect to the Lebesgue measure; in particular it puts no
mass at zero. The same holds for $\mathcal N(0,v)$ for every $v>0$.
\end{cor}

\begin{proof}
We first treat $W_\ell$, and we assume $\theta_\ell=\alpha_i$ for some $i$; the case
$\theta_\ell=\beta_j$ is identical, with the roles of $r$ and $s$ exchanged. Write
\[
  W_\ell=\sum_{j=1}^s c_j|g_{ij}|^2,
  \qquad
  c_j=\frac{\alpha_i\beta_j}{\alpha_i-\beta_j} .
\]
The numbers $\alpha_1,\ldots,\alpha_r,\beta_1,\ldots,\beta_s$ are positive and all
distinct, so every $c_j$ is well defined and nonzero. The variable $|g_{i1}|^2$ is
absolutely continuous: it follows the chi-square law with one degree of freedom in the
real case, and the exponential law of parameter one in the complex case. Multiplying by
the nonzero constant $c_1$ preserves absolute continuity. Put
\[
  \xi=c_1|g_{i1}|^2,
  \qquad
  \eta=\sum_{j=2}^sc_j|g_{ij}|^2 ,
\]
with $\eta=0$ if $s=1$, which is allowed in Lemma~\ref{a:conv}. The variables $\xi$ and
$\eta$ are independent because the $g_{ij}$ are, so Lemma~\ref{a:conv} shows that
$W_\ell$ is absolutely continuous.

If $\tau=0$ the claim is proved. If $\tau>0$, then
$\sqrt{\kappa\tau V_0(\theta_\ell)}G_\ell$ is a Gaussian variable with positive
variance, since $\kappa\ge1$, $\tau>0$ and $V_0(\theta_\ell)=(\sigma^2+\theta_\ell)^2>0$.
It is therefore absolutely continuous, and independent of $W_\ell$, so
Lemma~\ref{a:conv} applies again. The last assertion is the absolute continuity of a
Gaussian law with positive variance.
\end{proof}

\subsection{Second order eigenvalue perturbation lemma}
\label{appendix: second order perturbation lemma}

We recall Lemma \ref{l:pert} before proving it.

\begin{lem}[Second-order eigenvalue perturbation]
Let $M$ be Hermitian, let $\lambda$ be a simple eigenvalue of $M$, and let $v$ be a
corresponding unit eigenvector. Let $\{v_q\}$ be an orthonormal eigenbasis of $M$,
with $Mv_q=\lambda_qv_q$, and let $\ell$ be the index such that $v_\ell=v$. Set
\[
  \delta=\min_{\mu\in\operatorname{Spec}(M)\setminus\{\lambda\}}|\lambda-\mu|>0.
\]
If $E$ is Hermitian and $\|E\|_{\mathrm{op}}\leq\delta/4$, then the eigenvalue
$\widehat\lambda$ of $M+E$ associated with $\lambda$ satisfies
\begin{equation}
  \widehat\lambda-\lambda
  =v^*Ev
  +\sum_{q\neq\ell}
  \frac{|v_q^*Ev|^2}{\lambda-\lambda_q}
  +R,
  \label{e:pert-app}
\end{equation}
and there exists a universal constant $C>0$ such that
\begin{equation}
  |R|\leq C\delta^{-2}\|E\|_{\mathrm{op}}^3.
  \label{e:remd-app}
\end{equation}
\end{lem}

\begin{proof}
Let $P=vv^*$ and $Q=I-P$. With respect to the decomposition
$\operatorname{Ran}(P)\oplus\operatorname{Ran}(Q)$,
\[
  M+E=
  \begin{pmatrix}
    \lambda+v^*Ev & v^*EQ\\
    QEv & Q(M+E)Q
  \end{pmatrix}.
\]
Weyl's inequality gives $|\widehat\lambda-\lambda|\leq\|E\|_{\mathrm{op}}$. Since the
spectrum of $QMQ$ is at distance at least $\delta$ from $\lambda$ and
$\|E\|_{\mathrm{op}}\leq\delta/4$,
\[
  \operatorname{dist}\bigl(\widehat\lambda,\operatorname{Spec}(Q(M+E)Q)\bigr)
  \geq\delta-2\|E\|_{\mathrm{op}}
  \geq\frac\delta2.
\]
Thus
\[
  \left\|\bigl(\widehat\lambda I-Q(M+E)Q\bigr)^{-1}\right\|_{\mathrm{op}}
  \leq\frac2\delta.
\]
The Schur complement equation gives
\begin{equation}
  \widehat\lambda-\lambda
  =v^*Ev
  +v^*EQ\bigl(\widehat\lambda I-Q(M+E)Q\bigr)^{-1}QEv.
  \label{e:schu}
\end{equation}
Let
\[
  R_0=(\lambda I-QMQ)^{-1}.
\]
Then $\|R_0\|_{\mathrm{op}}\leq1/\delta$. The resolvent identity gives
\begin{align*}
 &\bigl(\widehat\lambda I-Q(M+E)Q\bigr)^{-1}-R_0\\
 &\quad=
 \bigl(\widehat\lambda I-Q(M+E)Q\bigr)^{-1}
 \bigl((\lambda-\widehat\lambda)I+QEQ\bigr)R_0.
\end{align*}
Therefore
\[
  \left\|\bigl(\widehat\lambda I-Q(M+E)Q\bigr)^{-1}-R_0\right\|_{\mathrm{op}}
  \leq\frac{4}{\delta^2}\|E\|_{\mathrm{op}}.
\]
Replacing the resolvent in \eqref{e:schu} by $R_0$ produces an error bounded by
$4\delta^{-2}\|E\|_{\mathrm{op}}^3$. Finally,
\[
  v^*EQR_0QEv
  =\sum_{q\neq\ell}\frac{|v_q^*Ev|^2}{\lambda-\lambda_q}.
\]
This proves the claim.
\end{proof}

\newpage
\section{Proof of Proposition \ref{p:pop}}
\label{appendix: proof of p:pop}

\subsection{Notations}

Throughout this section, $O_p$ is Haar distributed on $U(p)$ or $O(p)$, and
$\hat O^{(p)}=(o_{ij}^{(p)})_{i\in[r],\,j\in[s]}$ is its rectangular truncation. Recall
from Section~\ref{ss:setup} the diagonal matrices $A_p,B_p$ with nonzero eigenvalues
$\alpha_1,\dots,\alpha_r$ and $\beta_1,\dots,\beta_s$, and set
$\hat A=\diag(\alpha_1,\dots,\alpha_r)$, $\hat B=\diag(\beta_1,\dots,\beta_s)$. Let
$P(x,y)$ be a noncommutative polynomial in two variables with $P(0,0)=0$. Define
\[
  P_1(x):=P(x,0),\qquad Q_1(y):=P(0,y),
\]
and decompose
\[
  P(x,y)=P_1(x)+Q_1(y)+P_2(x,y)+Q_2(x,y)+P_3(x,y)+Q_3(x,y)+R(x,y),
\]
where
\[
  P_2(x,y)=\sum_{k,\ell\ge1}a_{k,\ell}\,x^ky^\ell,
  \qquad
  Q_2(x,y)=\sum_{k,\ell\ge1}b_{k,\ell}\,y^\ell x^k,
\]
\[
  P_3(x,y)=\sum_{k,\ell,m\ge1}c_{k,\ell,m}\,x^ky^\ell x^m,
  \qquad
  Q_3(x,y)=\sum_{k,\ell,m\ge1}d_{k,\ell,m}\,y^\ell x^k y^m,
\]
and $R$ is a linear combination of monomials of length greater than three (the length
counts each $x^k$ and each $y^\ell$ as one). We use the conventions
$a_{k,0}=$ coefficient of $x^k$ in $P_1$, $b_{0,\ell}=$ coefficient of $y^\ell$ in $Q_1$,
$c_{k,\ell,0}:=a_{k,\ell}$ and $d_{k,\ell,0}:=b_{k,\ell}$ for $k,\ell\ge1$. For a matrix
$X$ we write $\|X\|:=\sqrt{\operatorname{Tr}(X^*X)}$ for the Hilbert--Schmidt norm.

\subsection{The proposition in the general case}

\begin{prop}[Fluctuations of the eigenvalues of a polynomial in $A_p$ and
$O_pB_pO_p^*$; extension of Theorem~1.3 of \cite{collins24}]\label{a:polyfluc}
With the notation above, the matrix $P(A_p,O_pB_pO_p^*)$ has $p-r-s$ trivial zero
eigenvalues. The remaining $r+s$ eigenvalues, called the non-trivial eigenvalues,
converge almost surely to $\{P_1(\alpha_i)\}_{i\in[r]}\cup\{Q_1(\beta_j)\}_{j\in[s]}$ as
$p\to\infty$. If these $r+s$ limiting values are pairwise distinct, then the non-trivial
eigenvalues are of the form
\begin{align}
  &P_1(\alpha_i)
   +\frac1p\sum_{j\in[s]}p_{i,j}\bigl|\sqrt p\,o_{ij}^{(p)}\bigr|^2
   +\frac{\xi_i^{(p)}}{p^{3/2}},
   \qquad i\in[r],\label{eq:mu-alpha}\\
  &Q_1(\beta_j)
   +\frac1p\sum_{i\in[r]}q_{i,j}\bigl|\sqrt p\,o_{ij}^{(p)}\bigr|^2
   +\frac{\zeta_j^{(p)}}{p^{3/2}},
   \qquad j\in[s],\label{eq:mu-beta}
\end{align}
where
\begin{align}
  p_{i,j}&:=P_2(\alpha_i,\beta_j)+P_3(\alpha_i,\beta_j)
   +\frac{\bigl[P_1(\alpha_i)+P_2(\alpha_i,\beta_j)\bigr]
          \bigl[Q_1(\beta_j)+Q_2(\alpha_i,\beta_j)\bigr]}
         {P_1(\alpha_i)-Q_1(\beta_j)},\label{eq:pij}\\
  q_{i,j}&:=Q_2(\alpha_i,\beta_j)+Q_3(\alpha_i,\beta_j)
   -\frac{\bigl[P_1(\alpha_i)+P_2(\alpha_i,\beta_j)\bigr]
          \bigl[Q_1(\beta_j)+Q_2(\alpha_i,\beta_j)\bigr]}
         {P_1(\alpha_i)-Q_1(\beta_j)}.\label{eq:qij}
\end{align}
The constants $p_{i,j},q_{i,j}$ are the same for the two groups, and
$\xi_i^{(p)},\zeta_j^{(p)}$ converge in law to $\mathbb C$-valued random variables. The
family $\bigl\{|\sqrt p\,o_{ij}^{(p)}|^2:i\in[r],\,j\in[s]\bigr\}$ converges in law to a
family of i.i.d.\ random variables, whose common law is $\chi_1^2$ in the orthogonal case
and $\mathrm{Exp}(1)=\tfrac12\chi_2^2$ in the unitary case.
\end{prop}

Proposition~\ref{p:pop} is the special case $P(x,y)=x+y$, that is,
$\Sigma_p-\sigma^2I_p=A_p+O_pB_pO_p^*$. Then $P_1(x)=x$, $Q_1(y)=y$, and
$P_2,P_3,Q_2,Q_3$ all vanish, so $P_1(\alpha_i)=\alpha_i$, $Q_1(\beta_j)=\beta_j$, and
\eqref{eq:pij}--\eqref{eq:qij} reduce to
\[
  p_{i,j}=\frac{\alpha_i\beta_j}{\alpha_i-\beta_j},
  \qquad
  q_{i,j}=\frac{\alpha_i\beta_j}{\beta_j-\alpha_i}.
\]
By Proposition~\ref{p:haar} and the continuous mapping theorem,
\[
  \bigl(|\sqrt p\,o_{ij}^{(p)}|^2\bigr)_{i\in[r],\,j\in[s]}
  \xrightarrow{\ \mathrm{law}\ }
  \bigl(|g_{ij}|^2\bigr)_{i\in[r],\,j\in[s]}.
\]
Since $r,s$ are fixed and $\xi_i^{(p)},\zeta_j^{(p)}$ are tight, all the remainders $\xi_i^{(p)}/\sqrt p$ and $\zeta_j^{(p)}/\sqrt p$ converge jointly to zero in probability. After centering and multiplying \eqref{eq:mu-alpha}--\eqref{eq:mu-beta} by $p$, the continuous mapping theorem and Slutsky's theorem therefore give
\[
  \bigl(p(\lambda_{p,\ell}-\sigma^2-\theta_\ell)\bigr)_{\ell=1}^{r+s}
  \xrightarrow{\ \mathrm{law}\ }
  (W_\ell)_{\ell=1}^{r+s}.
\]
The limiting spike values are distinct, so the eigenvalue branches agree with the decreasing ordering for all sufficiently large $p$, almost surely. All the variables $W_\ell$ use the same family $(g_{ij})$ in \eqref{e:well}. Together with the almost sure convergence in Proposition~\ref{a:polyfluc}, this proves Proposition~\ref{p:pop}.

\subsection{Auxiliary lemmas for the population fluctuations}\label{s:apppop}

\begin{lem}[Density of the rectangular truncation]\label{a:density}
Let $r,s\ge1$ with $r+s\le p$, and let $O_p$ be Haar distributed on $U(p)$ or $O(p)$. Let
$\hat O^{(p)}=(o_{ij}^{(p)})_{i\in[r],\,j\in[s]}$ be the rectangular truncation of $O_p$.
\begin{enumerate}
  \item \textbf{Orthogonal case ($O(p)$).} Assume without loss of generality that
        $r\le s$ (otherwise swap the roles using
        $\det(I_r-\hat O\hat O^\top)=\det(I_s-\hat O^\top\hat O)$). Then $\hat O^{(p)}$ has
        probability density
        \[
          c_{p,r,s}\,\det\bigl(I_r-\hat O\hat O^\top\bigr)^{(p-r-s-1)/2}\,
          \mathbf 1_{\{\|\hat O\|_{\mathrm{op}}\le1\}}
        \]
        on $M_{r\times s}(\mathbb R)$, with respect to the Lebesgue measure.
  \item \textbf{Unitary case ($U(p)$).} The truncation $\hat O^{(p)}$ has probability
        density
        \[
          c'_{p,r,s}\,\det\bigl(I_r-\hat O\hat O^*\bigr)^{p-r-s}\,
          \mathbf 1_{\{\|\hat O\|_{\mathrm{op}}\le1\}}
        \]
        on $M_{r\times s}(\mathbb C)$, with respect to the Lebesgue measure
        $\prod_{i=1}^r\prod_{j=1}^s dx_{ij}\,dy_{ij}$, where
        $x_{ij}=\Re(o_{ij})$ and $y_{ij}=\Im(o_{ij})$.
\end{enumerate}
In particular, in both cases $\hat O^{(p)}$ admits a density with respect to the Lebesgue
measure. The asymptotic behavior of $\sqrt p\,\hat O^{(p)}$ is given by
Proposition~\ref{p:haar}.
\end{lem}

\begin{proof}
    The results come from \cite[Lemma 2.5]{jiang06} for the orthogonal case and \cite[Lemma~2.1 ]{collins24} for the unitary case.
\end{proof}

\begin{lem}[Haar projection]\label{a:haarproj}
Let $O_p\in U(p)$ be Haar distributed, and let $o_j$ be its $j$-th column. Let
$v\in\mathbb F^p$ be a fixed deterministic unit vector ($\|v\|_2=1$). Then:
\begin{enumerate}
  \item The vector $o_j$ is uniformly distributed on the unit sphere
        $\mathbb S(\mathbb F^p):=\{x\in\mathbb F^p:\|x\|_2=1\}$.
  \item The inner product $\langle v,o_j\rangle$ satisfies
        $\langle v,o_j\rangle=\mathrm O_{\mathbb P}(p^{-1/2})$.
  \item Let $(u_1,\dots,u_k)$ be a fixed orthonormal family in $\mathbb F^p$ and set
        $v_i:=O_pu_i$ for $i\in\{1,\dots,k\}$. Then for all fixed deterministic unit
        vectors $w_1,\dots,w_k\in\mathbb F^p$ (in particular for $w_i=u_i$), and all
        $i,j$, $\langle w_i,v_j\rangle=\mathrm O_{\mathbb P}(p^{-1/2})$.
\end{enumerate}
\end{lem}

\begin{proof}
\textbf{Part 1: Uniform distribution on the sphere.}
By the definition of the Haar measure on $U(p)$, the distribution of $O_p$ is
invariant under left-multiplication by any fixed deterministic $R\in U(p)$: that
is, $RO_p\stackrel{d}{=}O_p$.

The $j$-th column of $O_p$ is $o_j=O_pe_j$, where $e_j$ is the $j$-th vector of the
canonical basis. Applying the deterministic rotation $R$ to this column yields
\[
  Ro_j=R(O_pe_j)=(RO_p)e_j\stackrel{d}{=}O_pe_j=o_j.
\]
So the distribution of $o_j$ is invariant under any $R\in U(p)$. The group
$U(p)$ acts transitively on the unit sphere $\mathbb S(\mathbb F^p)$ (for
$\mathbb F=\mathbb R$ this is the action of $O(p)$ on $S^{p-1}$; for $\mathbb F=\mathbb C$
it is the action of $U(p)$ on the complex unit sphere), and $o_j$ has unit norm because
$O_p\in U(p)$ has orthonormal columns. A probability measure on the sphere that is
invariant under a transitive group of isometries is the uniform measure. Hence $o_j$ is
uniform on $\mathbb S(\mathbb F^p)$.

\textbf{Part 2: Asymptotic order of the projection.}
Since $v$ is a deterministic unit vector, there exists a fixed $R_v\in U(p)$ such
that $R_vv=e_1$, where $e_1=(1,0,\dots,0)^\top$ (such an $R_v$ exists by transitivity
from Part 1). By the invariance of the inner product under $R_v$ (which holds because
$R_v^*R_v=I$),
\[
  \langle v,o_j\rangle=\langle R_vv,R_vo_j\rangle=\langle e_1,R_vo_j\rangle.
\]
Since $o_j$ is uniform on $\mathbb S(\mathbb F^p)$, the rotated vector $u\equiv R_vo_j$ is
also uniform on $\mathbb S(\mathbb F^p)$. The inner product is then the first coordinate
of this uniform random vector:
\[
  \langle v,o_j\rangle\stackrel{d}{=}\langle e_1,u\rangle=\bar u_1.
\]
We compute the second moment of $u_1$. Since $u$ lies on the unit sphere,
$\sum_{k=1}^p|u_k|^2=1$ almost surely, so $\sum_{k=1}^p\mathbb E[|u_k|^2]=1$. By the
permutation symmetry of the uniform distribution on the sphere, all coordinates have the
same second moment, so $p\,\mathbb E[|u_1|^2]=1$, which gives
\[
  \mathbb E\bigl[|\langle v,o_j\rangle|^2\bigr]=\mathbb E[|u_1|^2]=\frac1p.
\]
By Markov's inequality applied to $|\langle v,o_j\rangle|^2$: for any $M>0$,
\[
  \mathbb P\bigl(\sqrt p\,|\langle v,o_j\rangle|\ge M\bigr)
  =\mathbb P\Bigl(|\langle v,o_j\rangle|^2\ge\tfrac{M^2}p\Bigr)
  \le\frac{p\,\mathbb E[|\langle v,o_j\rangle|^2]}{M^2}=\frac1{M^2}.
\]
This bound is uniform in $p$ and tends to $0$ as $M\to\infty$, so
$\sqrt p\,|\langle v,o_j\rangle|$ is bounded in probability, i.e.\
$\langle v,o_j\rangle=\mathrm O_{\mathbb P}(p^{-1/2})$.

\textbf{Part 3.}
Choose a fixed $R_i\in U(p)$ with $R_iu_j=e_1$; then
\[
  \langle w_i,v_j\rangle=\langle w_i,O_pu_j\rangle=\langle w_i,O_pR_i^*e_1\rangle
  =\langle(O_pR_i^*)^*w_i,e_1\rangle=\overline{\langle e_1,(O_pR_i^*)^*w_i\rangle}.
\]
Now \(O_pR_i^*\stackrel{d}{=}O_p\) by right-invariance of the Haar measure.
Hence
\[
  (O_pR_j^*)^*w_i
  \stackrel{d}{=}
  O_p^*w_i.
\]
Since $O_p^*$ is also Haar distributed, $O_p^*w_i$ is uniformly
distributed on $\mathbb S(\mathbb F^p)$ by Part 1. Therefore, writing
$z:=(O_pR_j^*)^*w_i$, we get
\[
  \langle w_i,v_j\rangle\stackrel{d}{=}\overline{\langle e_1,z\rangle}=\overline{z_1}.
\]
By the computation of Part 2, $\mathbb E[|z_1|^2]=1/p$, so by Markov's inequality
$\langle w_i,v_j\rangle=\mathrm O_{\mathbb P}(p^{-1/2})$. Taking $w_i=u_i$ gives
$\langle u_i,v_j\rangle=\mathrm O_{\mathbb P}(p^{-1/2})$. Notice that this means that the
fixed family $(u_i)$ and its random rotation $(v_j)=(O_pu_j)$ are asymptotically
orthogonal.
\end{proof}

\begin{lem}[Almost sure decay of the truncated blocks]\label{a:asdecay}
For $O_p$ Haar distributed on $U(p)$, we have $o_{ij}^{(p)}\to0$ almost surely as
$p\to\infty$, for every fixed pair $(i,j)$.
\end{lem}

\begin{proof}
The first column of $O_p$ is uniformly distributed on the unit sphere of $\mathbb F^p$
(this is Part 1 of Lemma~\ref{a:haarproj}, valid for the two groups). We compute the
fourth moment of one entry.

\textbf{Orthogonal case ($O(p)$).} The first column is uniform on the real sphere
$S^{p-1}\subset\mathbb R^p$, so $(o_{11}^{(p)})^2$ is $\mathrm{Beta}(1/2,(p-1)/2)$-distributed and
\[
  \mathbb E[(o_{ij}^{(p)})^4]=\frac3{p(p+2)}=\mathrm O(p^{-2}).
\]

\textbf{Unitary case ($U(p)$).} The first column is uniform on the complex sphere of
$\mathbb C^p$. The fourth moment is
\[
  \mathbb E[|o_{ij}^{(p)}|^4]=\frac2{p(p+1)}=\mathrm O(p^{-2}),
\]
which is the formula used in \cite{collins24} (proof of Lemma 2.3).

In both cases $\mathbb E[|o_{ij}^{(p)}|^4]=\mathrm O(p^{-2})$, so
\[
  \mathbb E\Bigl[\sum_{p=1}^\infty|o_{ij}^{(p)}|^4\Bigr]
  =\sum_{p=1}^\infty\mathbb E\bigl[|o_{ij}^{(p)}|^4\bigr]
  \le C\sum_{p=1}^\infty\frac1{p^2}<\infty.
\]
So the random series $\sum_p|o_{ij}^{(p)}|^4$ is finite almost surely, which forces
$|o_{ij}^{(p)}|\to0$ almost surely.
\end{proof}

\begin{lem}[Perturbation series for simple eigenvalues]\label{a:pertseries}
Let $\Lambda=\diag(\lambda_1,\dots,\lambda_d)\in M_d(\mathbb C)$ with
$\lambda_1,\dots,\lambda_d$ pairwise distinct. There exist homogeneous polynomials
$\Pi_{q,k}(X)$, $q\in[d]$, $k\ge1$, of degree $k$ in the $d^2$ entries of
$X=\{x_{ij}\}_{i,j\in[d]}$, and a constant $C>0$, such that for every $X\in M_d(\mathbb C)$
with $\|X\|<C$ the eigenvalues of $\Lambda+X$ are
\[
  \mu_q(X)=\lambda_q+\sum_{k=1}^\infty\Pi_{q,k}(X),\qquad q\in[d],
\]
with absolutely convergent series. The first two terms are
\[
  \Pi_{q,1}(X)=x_{qq},
  \qquad
  \Pi_{q,2}(X)=\sum_{i\ne q}\frac{x_{iq}\,x_{qi}}{\lambda_q-\lambda_i}.
\]
\end{lem}

\begin{proof}
See Collins et al. \cite{collins24}, where this is Lemma~2.5.
\end{proof}

\subsection{Proof of Proposition~\ref{p:pop}}\label{s:apppopprf}

\noindent \textbf{Step 1: matrix representation in an adapted basis.}
Let $o_j$ denote the $j$-th column of $O_p$. The family
\[
  \mathcal B=(e_1,\dots,e_r,\;o_1,\dots,o_s,\;e_{r+1},\dots,e_{p-s})
\]
is a basis of $\mathbb F^p$ with probability one. Indeed, the change-of-basis matrix has
determinant equal, up to sign, to the determinant of the $s\times s$ block of $O_p$ made
of rows $\{p-s+1,\dots,p\}$ and columns $\{1,\dots,s\}$; by
Lemma \ref{a:density} this block admits a density with respect to
Lebesgue measure (on $M_s(\mathbb R)$ in the orthogonal case, on $M_s(\mathbb C)$ in the
unitary case), so its determinant is non-zero a.s. Throughout the proof we work on this
event.

A direct calculation, using $Ae_i=\alpha_ie_i$ for $i\in[r]$, $Ae_i=0$ for $i>r$, the
relation $O_p^*o_j=e_j$ (which holds because $O_p^*O_p=I$, valid for both groups),
$Be_j=\beta_je_j$ for $j\in[s]$ and $Be_j=0$ for $j>s$, gives
\[
  \tilde A:=[A]_{\mathcal B}
   =\begin{pmatrix}\hat A&\hat A\hat O&0\\0&0&0\\0&0&0\end{pmatrix},
  \qquad
  \tilde B:=[O_pBO_p^*]_{\mathcal B}
   =\begin{pmatrix}0&0&0\\\hat B\hat O^*&\hat B&\ast\\0&0&0\end{pmatrix},
\]
where $\hat O^*=(\overline{o_{ij}})_{j\in[s],i\in[r]}$ is the $s\times r$ adjoint of
$\hat O$, and the $\ast$ in $\tilde B$ denotes a block depending on entries of $O_p$
outside the truncation $\hat O$ (it will play no role).

The only place where the field matters here is the adjoint: in the orthogonal case
$O_p^*=O_p^\top$ so $\hat O^*=\hat O^\top$ has the same real entries as $\hat O$ with
indices swapped; in the unitary case $\hat O^*$ has the conjugated entries
$\overline{o_{ij}}$. The block shapes are identical.

Consequently $P(A,O_pBO_p^*)$ has the same eigenvalues as $P(\tilde A,\tilde B)$, and in
particular $p-r-s$ trivial zero eigenvalues coming from the last block of rows/columns
(which is zero in both $\tilde A$ and $\tilde B$, so it annihilates every monomial of
degree $\ge1$, and $P(0,0)=0$ removes the constant term). The remaining $r+s$
non-trivial eigenvalues coincide with those of the upper-left $(r+s)\times(r+s)$ corner
of $P(\tilde A,\tilde B)$.

\noindent \textbf{Step 2: products in the algebra generated by $\tilde A,\tilde B$.}
By induction on $k,\ell\ge1$,
\[
  \tilde A^k=\begin{pmatrix}\hat A^k&\hat A^k\hat O&0\\0&0&0\\0&0&0\end{pmatrix},
  \qquad
  \tilde B^\ell=\begin{pmatrix}0&0&0\\\hat B^\ell\hat O^*&\hat B^\ell&\ast\\0&0&0\end{pmatrix}.
\]
For $k,\ell,m\ge1$, a direct block multiplication yields
\begin{align*}
  \tilde A^k\tilde B^\ell&=\begin{pmatrix}
   \hat A^k\hat O\hat B^\ell\hat O^*&\hat A^k\hat O\hat B^\ell&\ast\\
   0&0&0\\0&0&0\end{pmatrix},\\
  \tilde B^\ell\tilde A^k&=\begin{pmatrix}
   0&0&0\\
   \hat B^\ell\hat O^*\hat A^k&\hat B^\ell\hat O^*\hat A^k\hat O&0\\
   0&0&0\end{pmatrix},\\
  \tilde A^k\tilde B^\ell\tilde A^m&=\begin{pmatrix}
   \hat A^k\hat O\hat B^\ell\hat O^*\hat A^m&\hat A^k\hat O\hat B^\ell\hat O^*\hat A^m\hat O&0\\
   0&0&0\\0&0&0\end{pmatrix},\\
  \tilde B^\ell\tilde A^k\tilde B^m&=\begin{pmatrix}
   0&0&0\\
   \hat B^\ell\hat O^*\hat A^k\hat O\hat B^m\hat O^*&\hat B^\ell\hat O^*\hat A^k\hat O\hat B^m&\ast\\
   0&0&0\end{pmatrix}.
\end{align*}
These computations are the same in the two groups: they use only the block shapes from
Step 1 and the associativity of matrix multiplication. The only difference is again that
$\hat O^*$ stands for $\hat O^\top$ (orthogonal) or the conjugate transpose (unitary).

Any monomial of length $>3$ in $\tilde A,\tilde B$ has all of its blocks of order
$\mathrm O(\|\hat O\|^3)$ or zero. Hence, using the conventions for
$c_{k,\ell,0},d_{k,\ell,0}$,
\begin{multline*}
  P(\tilde A,\tilde B)=
  \begin{pmatrix}P_1(\hat A)&0&0\\0&Q_1(\hat B)&\ast\\0&0&0\end{pmatrix}
  +\!\!\sum_{k\ge1,\,\ell\ge0}\!\!
  \begin{pmatrix}0&a_{k,\ell}\hat A^k\hat O\hat B^\ell&\ast\\0&0&\ast\\0&0&0\end{pmatrix}\\
  +\!\!\sum_{k\ge0,\,\ell\ge1}\!\!
  \begin{pmatrix}0&0&\ast\\b_{k,\ell}\hat B^\ell\hat O^*\hat A^k&0&\ast\\0&0&0\end{pmatrix}
  +\!\!\sum_{\substack{k,\ell\ge1\\m\ge0}}\!\!
  \begin{pmatrix}c_{k,\ell,m}\hat A^k\hat O\hat B^\ell\hat O^*\hat A^m&0&0\\
   0&d_{k,\ell,m}\hat B^\ell\hat O^*\hat A^k\hat O\hat B^m&\ast\\0&0&0\end{pmatrix}\\
  +\begin{pmatrix}\mathrm O(\|\hat O\|^3)&\mathrm O(\|\hat O\|^2)&\ast\\
   \mathrm O(\|\hat O\|^2)&\mathrm O(\|\hat O\|^3)&\ast\\0&0&0\end{pmatrix}.
\end{multline*}
In $M:=\begin{pmatrix}P_1(\hat A)&0\\0&Q_1(\hat B)\end{pmatrix}$, the polynomials
$P_1,Q_1$ are applied entry-wise to $\hat A$ and $\hat B$.

\noindent \textbf{Step 3: a.s.\ convergence of non-trivial eigenvalues.}
The upper-left $(r+s)\times(r+s)$ corner of $P(\tilde A,\tilde B)$ equals $M+V_p$ where
\[
  V_p\to0\text{ a.s.\ (by Lemma~\ref{a:asdecay})}.
\]
This holds for the two groups, since Lemma~\ref{a:asdecay} gives $\hat O^{(p)}\to0$ a.s.\
in both of them, and every entry of $V_p$ is a polynomial in the entries of
$\hat O,\hat O^*$ with no constant term. By the argument principle applied to the
characteristic polynomials (see Lemma 4.1 from \cite{collins24}), the eigenvalues of
$M+V_p$ converge a.s.\ to those of $M$, namely
$\{P_1(\alpha_i)\}_{i\in[r]}\cup\{Q_1(\beta_j)\}_{j\in[s]}$. This proves the a.s.\
convergence stated in Proposition~\ref{a:polyfluc}.

\noindent \textbf{Step 4: perturbation expansion under the simplicity assumption.}
Assume from now on that the $r+s$ limits $\{P_1(\alpha_i)\}\cup\{Q_1(\beta_j)\}$ are
pairwise distinct. Then Lemma~\ref{a:pertseries} applies to $M+V_p$ (for $p$ large enough
so that $\|V_p\|<C$, which happens a.s.). Label the non-trivial eigenvalues so that
$\mu_i^{(p)}\to P_1(\alpha_i)$ ($i\in[r]$) and $\mu_{r+j}^{(p)}\to Q_1(\beta_j)$
($j\in[s]$). Lemma~\ref{a:pertseries} gives, up to order two,
\begin{align}
  \mu_i^{(p)}
  &=P_1(\alpha_i)+(V_p)_{i,i}
   +\sum_{\substack{q\in[r+s]\\q\ne i}}
    \frac{(V_p)_{q,i}(V_p)_{i,q}}{P_1(\alpha_i)-\lambda_q}
   +\mathrm O(\|V_p\|^3),\label{eq:mu-i}
\end{align}
where $\lambda_q\in\{P_1(\alpha_{i'})\}\cup\{Q_1(\beta_j)\}$ is the eigenvalue
corresponding to index $q$. This expansion is the same in the two groups
(Lemma~\ref{a:pertseries} is deterministic and the entries of $V_p$ are complex numbers
in all cases).

From the block decomposition of $P(\tilde A,\tilde B)$ in Step~2:
\begin{itemize}
\item The diagonal entry $(V_p)_{i,i}$ for $i\in[r]$ is
      $\sum_{k,\ell\ge1,\,m\ge0}c_{k,\ell,m}\bigl(\hat A^k\hat O\hat B^\ell\hat O^*\hat A^m\bigr)_{i,i}
      +\mathrm O(\|\hat O\|^3)$.
\item For $q=i'\in[r]\setminus\{i\}$, $(V_p)_{i',i}$ and $(V_p)_{i,i'}$ are
      $\mathrm O(\|\hat O\|^2)$, so their product is $\mathrm O(\|\hat O\|^4)$ and
      contributes to the $\mathrm O(\|\hat O\|^3)$ remainder.
\item For $q=r+j$ ($j\in[s]$), the entries $(V_p)_{r+j,i}$ and $(V_p)_{i,r+j}$ are, up to
      $\mathrm O(\|\hat O\|^3)$,
      \[
        \sum_{k\ge0,\,\ell\ge1}b_{k,\ell}\bigl(\hat B^\ell\hat O^*\hat A^k\bigr)_{j,i},
        \qquad
        \sum_{k\ge1,\,\ell\ge0}a_{k,\ell}\bigl(\hat A^k\hat O\hat B^\ell\bigr)_{i,j}.
      \]
\end{itemize}

\noindent \textbf{Step 5: entry-wise evaluation.}
Since $\hat A,\hat B$ are diagonal, a direct calculation gives, for any
$k,\ell,m,k',\ell'$ (using $(\hat O^*)_{j,i}=\overline{o_{ij}}$):
\begin{align}
  \bigl(\hat A^k\hat O\hat B^\ell\hat O^*\hat A^m\bigr)_{i,i}
   &=\alpha_i^{k+m}\sum_{j\in[s]}\beta_j^\ell\,o_{ij}\,\overline{o_{ij}}
    =\alpha_i^{k+m}\sum_{j\in[s]}\beta_j^\ell\,|o_{ij}|^2,\label{eq:diag}\\
  \bigl(\hat A^k\hat O\hat B^\ell\bigr)_{i,j}
   &=\alpha_i^k\,\beta_j^\ell\,o_{ij},\label{eq:Lij}\\
  \bigl(\hat B^{\ell'}\hat O^*\hat A^{k'}\bigr)_{j,i}
   &=\alpha_i^{k'}\,\beta_j^{\ell'}\,\overline{o_{ij}}.\label{eq:Rji}
\end{align}
The product of \eqref{eq:Lij} and \eqref{eq:Rji} equals
\[
  \alpha_i^{k+k'}\beta_j^{\ell+\ell'}\,o_{ij}\,\overline{o_{ij}}
   =\alpha_i^{k+k'}\beta_j^{\ell+\ell'}\,|o_{ij}|^2.
\]
The orthogonal and unitary cases differ, and it is the reason why the building block is
$o_{ij}^2$ in one case and $|o_{ij}|^2$ in the other. In the orthogonal case
($\mathbb F=\mathbb R$) the entries are real, $\overline{o_{ij}}=o_{ij}$, so
$o_{ij}\,\overline{o_{ij}}=o_{ij}^2=|o_{ij}|^2$. In the unitary case
($\mathbb F=\mathbb C$) the adjoint $\hat O^*$ carries the complex conjugate, so the
product is $o_{ij}\,\overline{o_{ij}}=|o_{ij}|^2$. In both cases the result is the same
non-negative real quantity $|o_{ij}|^2$.

Substituting \eqref{eq:diag}--\eqref{eq:Rji} into \eqref{eq:mu-i} and factoring out
$|o_{ij}|^2$:
\begin{align*}
  \mu_i^{(p)}
  &=P_1(\alpha_i)
   +\sum_{j\in[s]}\Bigl[\!\!\sum_{\substack{k,\ell\ge1\\m\ge0}}\!\!
     c_{k,\ell,m}\,\alpha_i^{k+m}\beta_j^\ell\Bigr]|o_{ij}|^2\\
  &\quad+\sum_{j\in[s]}\frac1{P_1(\alpha_i)-Q_1(\beta_j)}
     \Bigl[\!\!\sum_{\substack{k\ge1\\\ell\ge0}}\!\!a_{k,\ell}\,\alpha_i^k\beta_j^\ell\Bigr]
     \Bigl[\!\!\sum_{\substack{k'\ge0\\\ell'\ge1}}\!\!b_{k',\ell'}\,\alpha_i^{k'}\beta_j^{\ell'}\Bigr]
     |o_{ij}|^2\\
  &\quad+\mathrm O(\|\hat O\|^3).
\end{align*}
Using the conventions $c_{k,\ell,0}=a_{k,\ell}$ and $a_{k,0}=$ (coeff.\ of $x^k$ in
$P_1$), we recognize:
\begin{align*}
  \sum_{k,\ell\ge1,\,m\ge0}c_{k,\ell,m}\,\alpha_i^{k+m}\beta_j^\ell
   &=P_2(\alpha_i,\beta_j)+P_3(\alpha_i,\beta_j),\\
  \sum_{k\ge1,\,\ell\ge0}a_{k,\ell}\,\alpha_i^k\beta_j^\ell
   &=P_1(\alpha_i)+P_2(\alpha_i,\beta_j),\\
  \sum_{k\ge0,\,\ell\ge1}b_{k,\ell}\,\alpha_i^k\beta_j^\ell
   &=Q_1(\beta_j)+Q_2(\alpha_i,\beta_j).
\end{align*}
Hence
\begin{equation}\label{eq:mu-final}
  \mu_i^{(p)}=P_1(\alpha_i)+\sum_{j\in[s]}p_{i,j}\,|o_{ij}|^2+\mathrm O(\|\hat O\|^3),
\end{equation}
The analogous calculation for $j\in[s]$ yields
\[
  \mu_{r+j}^{(p)}=Q_1(\beta_j)+\sum_{i\in[r]}q_{i,j}\,|o_{ij}|^2+\mathrm O(\|\hat O\|^3),
\]
where $p_{i,j}$ and $q_{i,j}$ are given by:
\begin{align*}
  p_{i,j}&:=P_2(\alpha_i,\beta_j)+P_3(\alpha_i,\beta_j)
   +\frac{\bigl[P_1(\alpha_i)+P_2(\alpha_i,\beta_j)\bigr]
          \bigl[Q_1(\beta_j)+Q_2(\alpha_i,\beta_j)\bigr]}
         {P_1(\alpha_i)-Q_1(\beta_j)},\\
  q_{i,j}&:=Q_2(\alpha_i,\beta_j)+Q_3(\alpha_i,\beta_j)
   -\frac{\bigl[P_1(\alpha_i)+P_2(\alpha_i,\beta_j)\bigr]
          \bigl[Q_1(\beta_j)+Q_2(\alpha_i,\beta_j)\bigr]}
         {P_1(\alpha_i)-Q_1(\beta_j)}.
\end{align*}
The sign change in the second term of $q_{i,j}$ comes from
$Q_1(\beta_j)-P_1(\alpha_i)=-(P_1(\alpha_i)-Q_1(\beta_j))$ in Lemma~\ref{a:pertseries}.
The constants $p_{i,j},q_{i,j}$ do not depend on the group: only the building block
$|o_{ij}|^2$ does (it is $o_{ij}^2$ for $O(p)$ and $|o_{ij}|^2$ for $U(p)$).

\noindent \textbf{Step 6: rescaling and limiting law.}

\emph{Power series structure.} By Lemma~\ref{a:pertseries}, the eigenvalue $\mu_i^{(p)}$
admits an absolutely convergent expansion
\[
  \mu_i^{(p)}-P_1(\alpha_i)=\sum_{k=1}^\infty\Pi_{i,k}(V_p),
\]
true a.s.\ for $p$ large enough so that $\|V_p\|<C$, which holds by
Lemma~\ref{a:asdecay}. Each $\Pi_{i,k}$ is a homogeneous polynomial of degree $k$ in the
entries of $V_p$. Since the matrix $V_p$ is a polynomial in the entries
$\{o_{ij}\}_{i\in[r],\,j\in[s]}$ of $\hat O$ and their conjugates
$\{\overline{o_{ij}}\}$, with no constant term (as shown in Step~2),
$\sum_{k=1}^\infty\Pi_{i,k}(V_p)$ is a polynomial in the entries of $\hat O$ and their
conjugates. Substituting these expressions and regrouping by total degree in the real
coordinates of the entries $\{o_{ij}\}$, we obtain
\begin{equation}\label{eq:homog-exp}
  \mu_i^{(p)}-P_1(\alpha_i)=\sum_{n=2}^\infty\Phi_{i,n}(\hat O),
\end{equation}
where each $\Phi_{i,n}$ is a homogeneous polynomial of degree $n$ in the real
coordinates of the entries of $\hat O$ ($rs$ real coordinates in the orthogonal case,
$2rs$ in the unitary case, since each complex entry has a real and an imaginary part),
and the series converges absolutely for $\|\hat O\|$ small enough. The sum starts at
$n=2$ because: the diagonal entries $(V_p)_{q,q}$ are $\mathrm O(\|\hat O\|^2)$ (Step~4),
and the off-diagonal entries $(V_p)_{i,r+j}$, $(V_p)_{r+j,i}$ are
$\mathrm O(\|\hat O\|)$ but appear in $\Pi_{i,2}$ only through products, which are
$\mathrm O(\|\hat O\|^2)$. Hence no degree-$1$ term remains: $\Phi_{i,1}\equiv0$.

The degree-$2$ term $\Phi_{i,2}$ is exactly the explicit $\mathrm O(\|\hat O\|^2)$
contribution computed in Step~5:
\[
  \Phi_{i,2}(\hat O)=\sum_{j\in[s]}p_{i,j}\,|o_{ij}|^2,
\]
with $p_{i,j}$ given by \eqref{eq:pij}. Note that $|o_{ij}|^2$ is indeed a homogeneous
polynomial of degree $2$ in the real coordinates of $o_{ij}$ in both cases: it is
$o_{ij}^2$ for $O(p)$, and $\Re(o_{ij})^2+\Im(o_{ij})^2$ for $U(p)$.

\emph{Rescaling.} Since $\Phi_{i,n}$ is homogeneous of degree $n$, we have
$\Phi_{i,n}(\hat O)=p^{-n/2}\,\Phi_{i,n}(\sqrt p\,\hat O)$. Define
\[
  \xi_i^{(p)}
  :=p^{3/2}\Bigl(\mu_i^{(p)}-P_1(\alpha_i)-\sum_{j\in[s]}p_{i,j}\,|o_{ij}|^2\Bigr)
   =p^{3/2}\sum_{n=3}^\infty\Phi_{i,n}(\hat O)
   =\sum_{n=3}^\infty p^{(3-n)/2}\,\Phi_{i,n}\bigl(\sqrt p\,\hat O\bigr).
\]
This allows us to rewrite \eqref{eq:mu-final} as
\[
  \mu_i^{(p)}=P_1(\alpha_i)
   +\frac1p\sum_{j\in[s]}p_{i,j}\bigl|\sqrt p\,o_{ij}^{(p)}\bigr|^2
   +\frac{\xi_i^{(p)}}{p^{3/2}},
\]
which is the claimed expansion for $\mu_i^{(p)}$. The proof of the corresponding expansion for $\mu_{r+j}^{(p)}$ is identical, with
$q_{i,j}$ instead of $p_{i,j}$.

\emph{Convergence in law of $\xi_i^{(p)}$.} We show that $\xi_i^{(p)}$ converges in law.
By Proposition~\ref{p:haar}, $\sqrt p\,\hat O^{(p)}\xrightarrow{\ \mathrm{law}\ }G$,
where $G=(g_{ij})$ has i.i.d.\ entries (real $\mathcal N(0,1)$ in the orthogonal case,
standard complex Gaussian in the unitary case). We show that
$\xi_i^{(p)}\xrightarrow{\ \mathrm{law}\ }\Phi_{i,3}(G)$.

Write $\xi_i^{(p)}=T_p+R_p$ where $T_p:=\Phi_{i,3}\bigl(\sqrt p\,\hat O\bigr)$ and
$R_p:=\sum_{n=4}^\infty p^{(3-n)/2}\,\Phi_{i,n}\bigl(\sqrt p\,\hat O\bigr)$.

By the continuous mapping theorem applied to $\sqrt p\,\hat O$ and the continuous
function $\Phi_{i,3}$, $T_p=\Phi_{i,3}\bigl(\sqrt p\,\hat O\bigr)\xrightarrow{\ \mathrm{law}\ }\Phi_{i,3}(G)$.

We show $R_p\xrightarrow{\mathbb P}0$. Since $\sum_n\Phi_{i,n}(\hat O)$ converges
absolutely for $\|\hat O\|<\delta$ for some $\delta>0$, there exist constants $C_n>0$
with $|\Phi_{i,n}(X)|\le C_n\,\|X\|^n$ for all $X$ and
$\sum_{n=2}^\infty C_n\,\delta^n<\infty$. Hence, on the event $\{\|\hat O\|<\delta\}$
(which has probability tending to $1$, since $\|\hat O\|\xrightarrow{\mathrm{a.s.}}0$ by
Lemma~\ref{a:asdecay}),
\[
  |R_p|\le\sum_{n=4}^\infty p^{(3-n)/2}\,C_n\,\bigl(\sqrt p\,\|\hat O\|\bigr)^n
   \le p^{-1/2}\,D\,\bigl(\sqrt p\,\|\hat O\|\bigr)^4,
\]
where $D:=\sum_{n\ge4}C_n\delta^{n-4}<\infty$ (we used $p^{(3-n)/2}\le p^{-1/2}$ for
$n\ge4$, and $\|\hat O\|<\delta$ to bound $\|\hat O\|^{n-4}\le\delta^{n-4}$). Given that
$\sqrt p\,\|\hat O\|\xrightarrow{\ \mathrm{law}\ }\|G\|$, the continuous mapping theorem
yields $(\sqrt p\,\|\hat O\|)^4=\mathrm O_{\mathbb P}(1)$, hence
$p^{-1/2}D(\sqrt p\,\|\hat O\|)^4\xrightarrow{\mathbb P}0$. This gives
$R_p\xrightarrow{\mathbb P}0$.

As a result, by the Slutsky lemma,
\[
  \xi_i^{(p)}=T_p+R_p\xrightarrow{\ \mathrm{law}\ }\Phi_{i,3}(G),
\]
which is a $\mathbb C$-valued random variable. The same argument applies to
$\zeta_j^{(p)}$. None of these steps use the field: they only use the homogeneity of
$\Phi_{i,n}$, the a.s.\ decay of $\|\hat O\|$ (Lemma~\ref{a:asdecay}) and the convergence
in law of $\sqrt p\,\hat O$ (Proposition~\ref{p:haar}), all of which
hold for the two groups.

\emph{Limiting law of the building blocks.} Finally, by the continuous mapping theorem
applied to $\sqrt p\,\hat O^{(p)}\xrightarrow{\ \mathrm{law}\ }G$ and the map
$X\mapsto(|X_{ij}|^2)_{i\in[r],\,j\in[s]}$, we get the joint convergence
\[
  \bigl\{|\sqrt p\,o_{ij}^{(p)}|^2:\,i\in[r],\,j\in[s]\bigr\}
   \xrightarrow{\ \mathrm{law}\ }
   \bigl\{|g_{ij}|^2:\,i\in[r],\,j\in[s]\bigr\}.
\]
The law of $|g_{ij}|^2$ depends on the group:
\begin{itemize}
\item Orthogonal case: the $g_{ij}$ are i.i.d.\ real $\mathcal N(0,1)$, so
      $|g_{ij}|^2=g_{ij}^2$ are i.i.d.\ $\chi^2_1$;
\item Unitary case: the $g_{ij}$ are i.i.d.\ standard complex Gaussian, so
      $|g_{ij}|^2=\Re(g_{ij})^2+\Im(g_{ij})^2$ is a sum of two independent
      $\mathcal N(0,1/2)^2$ variables, which follows the exponential law
      $\mathrm{Exp}(1)=\tfrac12\chi^2_2$.
\end{itemize}

This completes the proof. \qed

\newpage
\section{Proof of Proposition \ref{p:dclt}}

\label{appendix: proof of p:dclt}

We first explain why resolvents and the Marchenko--Pastur transform
enter the proof. Let $\lambda_p$ denote the target population
eigenvalue, including the noise level $\sigma^2$, and let
$\widehat\lambda_p$ be the corresponding sample eigenvalue.
Although $p/n\to0$ gives consistency, the deterministic sampling
shift is of order $p/n$, which need not be negligible at the
$n^{-1/2}$ fluctuation scale.

After placing the target coordinate first, the sample matrix has
the block form
\[
  S_p=\begin{pmatrix}\lambda_p q_p&b_p^*\\ b_p&C_p\end{pmatrix},
  \qquad q_p=\frac1n\sum_{k=1}^n|z_k|^2,
\]
where the $z_k$ are independent standard real or circular complex
Gaussian variables, $C_p$ is the remaining sample covariance block,
and $b_p$ is the off-diagonal block.
On an event of probability tending to one, a fixed interval around
the target spike is separated from the spectrum of $C_p$.
The Schur complement then gives the exact equation
\[
  \widehat\lambda_p=\lambda_p q_p+B_p(\widehat\lambda_p),
  \qquad B_p(t)=b_p^*(tI-C_p)^{-1}b_p.
\]
The first term has a one-dimensional central limit theorem.
The second term accounts for the interaction with the remaining
coordinates. For fixed $t$, its conditional mean, given the
remaining Gaussian rows, is
\[
  \frac{\lambda_p}{n}
  \operatorname{Tr}\{C_p(tI-C_p)^{-1}\}.
\]
The Gaussian quadratic-form estimate in Step~5 below shows that
replacing $B_p(t)$ by this mean has a uniform error
$\mathrm O_{\mathbb P}(\sqrt p/n)=o_{\mathbb P}(n^{-1/2})$
on the fixed interval. Uniformity is needed because the equation
is evaluated at the random point $t=\widehat\lambda_p$.

Let $K$ be the fixed dimension of the spike block in
Proposition~\ref{p:dclt}, and put $c_p=(p-K)/n$,
the aspect ratio of the pure-bulk block.
Step~6 removes the other $K-1$ spike coordinates with an
$\mathrm O(n^{-1})$ trace error on the same high-probability event.
Lemma~\ref{l:resb} then approximates the pure-bulk trace by
$h_{c_p}(t)$, defined in \eqref{e:appsh}, with uniform error
$o_{\mathbb P}(n^{-1/2})$.
Lemma~\ref{l:mpid} identifies the deterministic solution
\[
  \rho_p=\lambda_p\left(1+c_p\frac{\sigma^2}{\lambda_p-\sigma^2}\right)
  =\lambda_p\{1+h_{c_p}(\rho_p)\}.
\]
Finally, the bound $h_{c_p}'=\mathrm O(c_p)$ allows us to deduce
\[
  \sqrt n(\widehat\lambda_p-\rho_p)
  =\lambda_p\sqrt n(q_p-1)+o_{\mathbb P}(1).
\]
The desired CLT is therefore reduced to the usual CLT for $n^{-1}\sum_k|z_k|^2$.
Replacing $c_p=(p-K)/n$ by $p/n$ changes the centering by $\mathrm O(n^{-1})$, which is negligible at this scale.
We now state and prove the two lemmas, followed by the full proof.

\begin{lem}[Pure-bulk resolvent trace estimate]\label{l:resb}
Let \(m=m(p)\) and \(n=n(p)\) satisfy \(m/n\to0\).  Let \(Z\) be an \(m\times n\) real or complex standard Gaussian matrix and set
\[
  C=C_{m,n}=\frac{\sigma^2}{n}ZZ^*.
\]
Put \(c=c_{m,n}=m/n\).  Let \(I\) be a fixed compact interval contained in \((\sigma^2+\varepsilon,\infty)\) for some \(\varepsilon>0\).  Then
\begin{equation}
  \sup_{t\in I}
  \left|
    \frac1n\operatorname{Tr}\bigl(C(tI-C)^{-1}\bigr)
    -h_c(t)
  \right|
  =o_{\mathbb P}(n^{-1/2}),
  \label{e:resb}
\end{equation}
where \(h_c\) is defined in Appendix~\ref{s:appmp}.  Moreover,
\begin{equation}
  \sup_{t\in I}|h_c'(t)|=\mathrm O(c).
  \label{e:hder}
\end{equation}
\end{lem}

\begin{proof}
We give the proof for the real and complex cases at the same time.  Constants may change from line to line.  Write the columns of \(\sigma Z\) as \(x_1,\ldots,x_n\in\mathbb F^m\), where \(\mathbb F=\mathbb R\) or \(\mathbb C\).  Thus
\[
  C=\frac1n\sum_{k=1}^n x_kx_k^*,
  \qquad
  x_k\sim N_{\mathbb F}(0,\sigma^2I_m),
\]
and the vectors \(x_k\) are independent.

For \(t\notin\operatorname{Spec}(C)\), define
\[
  Q(t)=(tI-C)^{-1},
  \qquad
  s_n(t)=\frac1n\operatorname{Tr}Q(t),
  \qquad
  r_n(t)=\frac1n\operatorname{Tr}\bigl(CQ(t)\bigr).
\]
The identity \(CQ(t)=tQ(t)-I_m\) gives
\begin{equation}
  r_n(t)=t s_n(t)-c.
  \label{e:l62alg}
\end{equation}

Fix \(0<\eta<\varepsilon/8\).  Let
\[
  \Omega_n^{op}=\{\|C-\sigma^2I_m\|_{\mathrm{op}}\le\eta\},
  \qquad
  \Omega_n^{col}=\left\{\max_{k\le n}\frac{\|x_k\|^2}{n}\le\eta\right\},
\]
and put \(\Omega_n=\Omega_n^{op}\cap \Omega_n^{col}\).  By the Gaussian covariance bound and the chi-square tail estimate in Appendix~\ref{s:appgauss},
\begin{equation}
  \mathbb P(\Omega_n)\to1.
  \label{e:l62omega}
\end{equation}
For each \(k\), define the leave-one-out matrix
\[
  C^{(k)}=C-\frac1n x_kx_k^*,
  \qquad
  Q_k(t)=(tI-C^{(k)})^{-1},
\]
and the event
\[
  \Omega_n^{(k)}=\{\|C^{(k)}-\sigma^2I_m\|_{\mathrm{op}}\le2\eta\}.
\]
The event \(\Omega_n^{(k)}\) is independent of \(x_k\), and \(\Omega_n\subset\Omega_n^{(k)}\).  Indeed, on \(\Omega_n\),
\[
  \|C^{(k)}-\sigma^2I_m\|_{\mathrm{op}}
  \leq \|C-\sigma^2I_m\|_{\mathrm{op}}+\frac{\|x_k\|^2}{n}
  \leq2\eta.
\]

On \(\Omega_n^{op}\), every eigenvalue of \(C\) lies in \([\sigma^2-\eta,\sigma^2+\eta]\).  Since \(I\subset(\sigma^2+\varepsilon,\infty)\),
\begin{equation}
  \sup_{t\in I}\|Q(t)\|_{\mathrm{op}}\le C_0
  \quad\text{on }\Omega_n^{op}.
  \label{e:l62Q}
\end{equation}
Similarly,
\begin{equation}
  \sup_{t\in I}\|Q_k(t)\|_{\mathrm{op}}\le C_0
  \quad\text{on }\Omega_n^{(k)}.
  \label{e:l62Qk}
\end{equation}
Here \(C_0\) depends only on \(\varepsilon\).  Put
\[
  q_k(t)=\frac1n x_k^*Q_k(t)x_k,
  \qquad
  s_k(t)=\frac1n\operatorname{Tr}Q_k(t).
\]
On \(\Omega_n\), using \eqref{e:l62Qk},
\begin{equation}
  \sup_{t\in I}\max_{k\le n}|q_k(t)|
  \le C_0\max_{k\le n}\frac{\|x_k\|^2}{n}
  \le C_0\eta.
  \label{e:l62qsmall0}
\end{equation}
By taking \(\eta\) smaller if needed, we may and do assume that the right-hand side is at most \(1/6\).  Also,
\begin{equation}
  \sup_{t\in I}|s_n(t)|\le C_0c,
  \qquad
  \sup_{t\in I}\max_{k\le n}|s_k(t)|\le C_0c
  \label{e:l62ssmall}
\end{equation}
on the relevant events.  Since \(c\to0\), these quantities are small for all large \(n\).

Let \(\psi(u)=u/(1-u)\).  To avoid problems away from the good event, choose a smooth function \(\chi\) such that \(\chi(u)=1\) for \(|u|\le1/2\), \(\chi(u)=0\) for \(|u|\ge3/4\), and \(0\le\chi\le1\).  Set
\[
  \widetilde\psi(u)=\chi(u)\frac{u}{1-u}.
\]
Then \(\widetilde\psi\), \(\widetilde\psi'\), and \(\widetilde\psi''\) are bounded on \(\mathbb R\), and \(\widetilde\psi=\psi\) on \([-1/2,1/2]\).

By the Sherman--Morrison formula in Appendix~\ref{s:applinear}, on \(\Omega_n\),
\[
  \frac1n x_k^*Q(t)x_k
  =\frac{q_k(t)}{1-q_k(t)}
  =\psi(q_k(t)).
\]
Hence, still on \(\Omega_n\),
\begin{equation}
  r_n(t)
  =\frac1n\sum_{k=1}^n\psi(q_k(t))
  =\frac1n\sum_{k=1}^n\widetilde\psi(q_k(t)).
  \label{e:l62rep}
\end{equation}
Define
\[
  E_n(t)=\frac1n\sum_{k=1}^n\widetilde\psi(q_k(t))
  -\widetilde\psi(\sigma^2s_n(t)).
\]
By \eqref{e:l62rep}, \eqref{e:l62alg}, and the smallness of \(s_n(t)\), on \(\Omega_n\),
\begin{equation}
  E_n(t)
  =r_n(t)-\frac{\sigma^2s_n(t)}{1-\sigma^2s_n(t)}
  =F_t(s_n(t)),
  \label{e:l62EequalsF}
\end{equation}
where
\[
  F_t(s)=ts-c-\frac{\sigma^2s}{1-\sigma^2s}.
\]

We next estimate \(E_n(t)\).  Let \(\mathbb E_k\) denote expectation with respect to \(x_k\), with all other columns fixed.  Since \(Q_k(t)\) is independent of \(x_k\), the Gaussian quadratic-form formula in Appendix~\ref{s:appquad} gives
\begin{equation}
  \mathbb E_k q_k(t)=\sigma^2s_k(t),
  \qquad
  \mathbb E_k|q_k(t)-\sigma^2s_k(t)|^2
  \le C\frac{m}{n^2}
  \label{e:l62qvar}
\end{equation}
on \(\Omega_n^{(k)}\), uniformly in \(t\in I\).  Because \(\Omega_n\subset\Omega_n^{(k)}\) and \(\Omega_n^{(k)}\) is independent of \(x_k\), Cauchy's inequality gives
\begin{equation}
  \mathbb E\left[\mathbf1_{\Omega_n}|q_k(t)-\sigma^2s_k(t)|\right]
  \le C\frac{\sqrt m}{n}.
  \label{e:l62qL1}
\end{equation}
Also, by the resolvent identity and the rank-one trace bound in Appendix~\ref{s:applinear},
\begin{equation}
  \sup_{t\in I}\max_{k\le n}|s_k(t)-s_n(t)|\le \frac{C}{n}
  \quad\text{on }\Omega_n.
  \label{e:l62sk}
\end{equation}
Since \(\widetilde\psi\) is Lipschitz, \eqref{e:l62qL1} and \eqref{e:l62sk} imply, for each fixed \(t\in I\),
\begin{equation}
  \mathbb E\left[\mathbf1_{\Omega_n}|E_n(t)|\right]
  \le C\left(\frac{\sqrt m}{n}+\frac1n\right)
  =o(n^{-1/2}).
  \label{e:l62Epoint}
\end{equation}

We need the estimate uniformly in \(t\).  On \(\Omega_n\), the functions above are \(C^1\) in \(t\), and
\[
  Q'(t)=-Q(t)^2,
  \qquad
  Q_k'(t)=-Q_k(t)^2.
\]
Thus
\[
  q_k'(t)=-\frac1n x_k^*Q_k(t)^2x_k,
  \qquad
  s_k'(t)=-\frac1n\operatorname{Tr}Q_k(t)^2.
\]
The same quadratic-form estimate as above, with \(Q_k(t)^2\) in place of \(Q_k(t)\), gives
\begin{equation}
  \mathbb E\left[\mathbf1_{\Omega_n}|q_k'(t)-\sigma^2s_k'(t)|\right]
  \le C\frac{\sqrt m}{n}.
  \label{e:l62qprime}
\end{equation}
The same rank argument gives
\begin{equation}
  \sup_{t\in I}\max_{k\le n}|s_k'(t)-s_n'(t)|\le \frac{C}{n}
  \quad\text{on }\Omega_n,
  \label{e:l62skprime}
\end{equation}
and \(\sup_{t\in I}|s_n'(t)|\le C\) on \(\Omega_n\).  Since \(\widetilde\psi'\) is bounded and Lipschitz, differentiating \(E_n(t)\) and using \eqref{e:l62qL1}--\eqref{e:l62skprime} gives
\begin{equation}
  \sup_{t\in I}\mathbb E\left[\mathbf1_{\Omega_n}|E_n'(t)|\right]
  \le C\left(\frac{\sqrt m}{n}+\frac1n\right)
  =o(n^{-1/2}).
  \label{e:l62Eprime}
\end{equation}
Fix \(t_0\in I\).  By the fundamental theorem of calculus,
\[
  \sup_{t\in I}|E_n(t)|
  \le |E_n(t_0)|+\int_I |E_n'(u)|\,du.
\]
Combining \eqref{e:l62Epoint}, \eqref{e:l62Eprime}, Fubini's theorem, and Lemma~\ref{l:localization}, we obtain
\begin{equation}
  \sup_{t\in I}|E_n(t)|=o_{\mathbb P}(n^{-1/2}).
  \label{e:l62Eunif}
\end{equation}
Together with \eqref{e:l62EequalsF}, this gives
\begin{equation}
  \sup_{t\in I}|F_t(s_n(t))|=o_{\mathbb P}(n^{-1/2}).
  \label{e:l62Funif}
\end{equation}

Now compare \(s_n(t)\) with the Marchenko--Pastur quantity \(s_c(t)\) from Appendix~\ref{s:appmp}.  For all large \(n\), the right edge \(b_c=\sigma^2(1+\sqrt c)^2\) is at most \(\sigma^2+\varepsilon/2\), so \(I\subset(b_c,\infty)\).  By Proposition~\ref{a:mpself},
\[
  F_t(s_c(t))=0,
  \qquad t\in I.
\]
Moreover,
\[
  0\le s_c(t)
  =c\int\frac1{t-x}\,\mu_{c,\sigma^2}(dx)
  \le \frac{2c}{\varepsilon},
  \qquad t\in I.
\]
On \(\Omega_n\), \eqref{e:l62ssmall} gives \(0\le s_n(t)\le Cc\).  Hence both \(s_n(t)\) and \(s_c(t)\) lie in \([0,\delta_n]\), where \(\delta_n\to0\).  For \(s\in[0,\delta_n]\) and \(t\in I\),
\[
  F_t'(s)=t-\frac{\sigma^2}{(1-\sigma^2s)^2}
  \ge \frac\varepsilon2
\]
for all large \(n\).  The mean value theorem gives
\[
  \sup_{t\in I}|s_n(t)-s_c(t)|
  \le C\sup_{t\in I}|F_t(s_n(t))|
  =o_{\mathbb P}(n^{-1/2}).
\]
Finally, by \(r_n(t)=ts_n(t)-c\) and \(h_c(t)=ts_c(t)-c\),
\[
  \sup_{t\in I}|r_n(t)-h_c(t)|
  \le (\max I)\sup_{t\in I}|s_n(t)-s_c(t)|
  =o_{\mathbb P}(n^{-1/2}).
\]
This proves \eqref{e:resb}.

It remains only to prove \eqref{e:hder}.  Since the support of \(\mu_{c,\sigma^2}\) is contained in \([a_c,b_c]\), and for large \(n\) this support is separated from \(I\), differentiation under the integral sign is valid and gives
\[
  h_c'(t)
  =-c\int\frac{x}{(t-x)^2}\,\mu_{c,\sigma^2}(dx).
\]
For \(t\in I\) and \(x\in[a_c,b_c]\), the denominator is bounded below by a positive constant and \(x\) is bounded above by a constant depending only on \(\sigma^2\) and \(\varepsilon\).  Hence \(\sup_{t\in I}|h_c'(t)|=\mathrm O(c)\).
\end{proof}

\begin{lem}[Outlier identity for the identity bulk]\label{l:mpid}
Let \(\lambda>\sigma^2\) and let \(c\ge0\).  Assume that
\begin{equation}
  \sqrt c<\frac{\lambda-\sigma^2}{\sigma^2}.
  \label{e:l63super}
\end{equation}
For \(c>0\), assume also \(c<1\).  Define
\begin{equation}
  \rho_c(\lambda)
  =\lambda\left(1+c\frac{\sigma^2}{\lambda-\sigma^2}\right).
  \label{e:rhoc}
\end{equation}
Then
\begin{equation}
  \rho_c(\lambda)-\sigma^2(1+\sqrt c)^2
  =\frac{\bigl(\lambda-
  \sigma^2-\sigma^2\sqrt c\bigr)^2}{\lambda-\sigma^2}.
  \label{e:rhogap}
\end{equation}
In particular, under \eqref{e:l63super}, \(\rho_c(\lambda)>b_c=\sigma^2(1+\sqrt c)^2\).  Moreover,
\begin{equation}
  h_c(\rho_c(\lambda))
  =c\frac{\sigma^2}{\lambda-\sigma^2},
  \qquad
  \rho_c(\lambda)=\lambda\bigl(1+h_c(\rho_c(\lambda))\bigr).
  \label{e:outid}
\end{equation}
\end{lem}

\begin{proof}
Put \(d=\lambda-\sigma^2>0\) and \(\rho=\rho_c(\lambda)\).  First,
\[
\begin{split}
  \rho-\sigma^2(1+\sqrt c)^2
  &=d+\frac{c\sigma^2\lambda}{d}-2\sigma^2\sqrt c-c\sigma^2  \\
  &=d+\frac{c\sigma^4}{d}-2\sigma^2\sqrt c
   =\frac{(d-\sigma^2\sqrt c)^2}{d}.
\end{split}
\]
This proves \eqref{e:rhogap}.  Hence \(\rho>b_c\) under \eqref{e:l63super}.

If \(c=0\), then \(\rho_0(\lambda)=\lambda\) and \(h_0\equiv0\), so the claim is immediate.  Assume now \(0<c<1\).  Define
\[
  s=\frac{c}{d+c\sigma^2}.
\]
Then \(0<s<1/\sigma^2\), and
\[
  \rho s-c
  =\lambda\left(1+c\frac{\sigma^2}{d}\right)\frac{c}{d+c\sigma^2}-c
  =\frac{c\sigma^2}{d}.
\]
Also,
\[
  \frac{\sigma^2s}{1-\sigma^2s}
  =\frac{c\sigma^2}{d}.
\]
Thus
\[
  F_\rho(s)
  :=\rho s-c-\frac{\sigma^2s}{1-\sigma^2s}=0.
\]
We must identify this root with the Marchenko--Pastur root \(s_c(\rho)\).  By Proposition~\ref{a:mproot}, for \(\rho>b_c\), the root corresponding to \(s_c(\rho)\) is the root of \(F_\rho(s)=0\) at which \(F_\rho'(s)>0\).  Now
\[
  F_\rho'(s)=\rho-\frac{\sigma^2}{(1-\sigma^2s)^2}.
\]
Using \(1-\sigma^2s=d/(d+c\sigma^2)\) and \(\rho=\lambda(d+c\sigma^2)/d\), we get
\[
  F_\rho'(s)>0
  \quad\Longleftrightarrow\quad
  \lambda d>\sigma^2(d+c\sigma^2)
  \quad\Longleftrightarrow\quad
  d^2>c\sigma^4.
\]
This last condition is exactly \eqref{e:l63super}.  Hence \(s=s_c(\rho)\).  Therefore, using \(h_c(t)=ts_c(t)-c\),
\[
  h_c(\rho)=\rho s-c=\frac{c\sigma^2}{d}
  =c\frac{\sigma^2}{\lambda-\sigma^2}.
\]
This is the first identity in \eqref{e:outid}.  The second follows from the definition of \(\rho\), because
\[
  \frac{\rho}{\lambda}-1
  =c\frac{\sigma^2}{\lambda-\sigma^2}
  =h_c(\rho).
\]
\end{proof}

We can finally state the proof of Proposition \ref{p:dclt}:

\begin{proof}
We first give the main line of the proof. After that, we prove the only technical estimate used in the main line.

\medskip
\noindent\textbf{Step 1: Put the target spike in the first coordinate.}
Use a deterministic permutation to move coordinate \(\ell\) to coordinate \(1\). This only conjugates \(\Lambda_p\) and \(S_p\), so it does not change their eigenvalues. After the permutation, relabel the spike values again as
\(\lambda_{p,1},\ldots,\lambda_{p,K}\), with the target value in the first place. We write
\[
  \lambda_p=\lambda_{p,1},
  \qquad
  \lambda_*=\lambda_\ell^*,
  \qquad
  \widehat\lambda_p=\widehat\lambda_{p,\ell}^{\mathrm{det}},
\]
\[
  m=p-K,
  \qquad
  c_p=\frac mn=\frac{p-K}{n}.
\]
At the end, we return to the original order of the coordinates.

By Lemma~\ref{a:sepint}, there are fixed numbers \(\nu>0\), \(\varepsilon>0\), an inner interval \(J\), and a compact interval \(I\) such that, for all large \(p\),
\begin{equation}
  \lambda_p\in J\subset\operatorname{int}I,
  \qquad
  \operatorname{dist}(\lambda_p,\partial I)\ge\nu,
  \label{e:p64I1}
\end{equation}
\begin{equation}
  \operatorname{dist}\!\left(
    I,\{\sigma^2,\lambda_{p,2},\ldots,\lambda_{p,K}\}
  \right)\ge4\nu,
  \qquad
  I\subset(\sigma^2+\varepsilon,\infty).
  \label{e:p64I2}
\end{equation}

Write the data matrix as
\[
  Y_p=\Lambda_p^{1/2}Z_p,
  \qquad
  Z_p=\begin{pmatrix}z\\ Z_0\end{pmatrix},
\]
where \(z\in\mathbb F^{1\times n}\),
\(Z_0\in\mathbb F^{(p-1)\times n}\), and \(z\) is independent of \(Z_0\). Here \(\mathbb F=\mathbb R\) or \(\mathbb C\). Put
\[
  D_p=\operatorname{diag}(\lambda_{p,2},\ldots,\lambda_{p,K},
  \sigma^2,\ldots,\sigma^2).
\]
Then
\begin{equation}
  S_p=
  \begin{pmatrix}
    \lambda_pq_p&b_p^*\\
    b_p&C_p
  \end{pmatrix},
  \label{e:p64block}
\end{equation}
where
\[
  q_p=\frac{zz^*}{n}=\frac1n\sum_{k=1}^n|z_k|^2,
  \qquad
  b_p=\frac{\sqrt{\lambda_p}}{n}D_p^{1/2}Z_0z^*,
  \qquad
  C_p=\frac1nD_p^{1/2}Z_0Z_0^*D_p^{1/2}.
\]

Let \(Z_b\) be the last \(m\) rows of \(Z_0\), and set
\[
  C_b=\frac{\sigma^2}{n}Z_bZ_b^*.
\]
Fix \(0<\eta\le\nu\), and define
\begin{equation}
\begin{split}
  \Omega_p^{(0)}={}&
  \{\|C_p-D_p\|_{\mathrm{op}}\le\eta\}\\
  &\cap
  \{\|C_b-\sigma^2I_m\|_{\mathrm{op}}\le\eta\},
\end{split}
  \label{e:p64omega0def}
\end{equation}
\begin{equation}
  \Omega_p^{(1)}
  =\{\|S_p-\Lambda_p\|_{\mathrm{op}}\le\nu\}.
  \label{e:p64omega1def}
\end{equation}
The first event depends only on \(Z_0\). By Corollary~\ref{a:wishartop},
\begin{equation}
  \mathbb P(\Omega_p^{(0)})\to1,
  \qquad
  \mathbb P(\Omega_p^{(1)})\to1.
  \label{e:p64omegaprob}
\end{equation}
On \(\Omega_p^{(0)}\), Lemma~\ref{a:weyl} and \eqref{e:p64I2} give
\begin{equation}
  \operatorname{dist}(I,\operatorname{Spec}(C_p))\ge3\nu,
  \qquad
  \sup_{t\in I}\|(tI-C_p)^{-1}\|_{\mathrm{op}}
  \le\frac1{3\nu}.
  \label{e:p64resbd}
\end{equation}
Also, on this event,
\begin{equation}
  \|C_p\|_{\mathrm{op}}
  \le \|D_p\|_{\mathrm{op}}+\eta\le C.
  \label{e:p64Cpbd}
\end{equation}

The eigenvalue in \eqref{e:p64label} was defined by its rank, so it is defined on the whole probability space. On \(\Omega_p^{(1)}\), Lemma~\ref{a:weyl} gives
\[
  |\widehat\lambda_p-\lambda_p|\le\nu.
\]
Thus \(\widehat\lambda_p\in I\). Every other population eigenvalue is at distance at least \(4\nu\) from \(I\), so every other sample eigenvalue is at distance at least \(3\nu\) from \(I\). Hence
\begin{equation}
  \operatorname{Spec}(S_p)\cap I=\{\widehat\lambda_p\}
  \quad\text{on }\Omega_p^{(1)}.
  \label{e:p64unique}
\end{equation}

Define the deterministic center
\begin{equation}
  \rho_p
  =\lambda_p\left(
    1+c_p\frac{\sigma^2}{\lambda_p-\sigma^2}
  \right).
  \label{e:rhop}
\end{equation}
Because \(\lambda_p-\sigma^2\ge\delta_0\),
\(\lambda_p\) is bounded, and \(c_p\to0\), we have
\(\rho_p-\lambda_p\to0\). Therefore \(\rho_p\in I\) for all large \(p\). Also, the condition of Lemma~\ref{l:mpid} holds for all large \(p\), and
\begin{equation}
  \rho_p
  =\lambda_p\{1+h_{c_p}(\rho_p)\}.
  \label{e:rhoeq}
\end{equation}

\medskip
\noindent\textbf{Step 2: State the key estimate.}
On \(\Omega_p^{(0)}\), define, for \(t\in I\),
\[
  R_p(t)=(tI-C_p)^{-1},
  \qquad
  B_p(t)=b_p^*R_p(t)b_p.
\]
Outside \(\Omega_p^{(0)}\), set these quantities equal to zero. We will prove below that
\begin{equation}
  \Delta_p
  :=\mathbf1_{\Omega_p^{(0)}}
  \sup_{t\in I}|B_p(t)-\lambda_ph_{c_p}(t)|
  =o_{\mathbb P}(n^{-1/2}).
  \label{e:p64key}
\end{equation}
This estimate is uniform in \(t\), so it may later be used at the random point
\(t=\widehat\lambda_p\).

\medskip
\noindent\textbf{Step 3: Use the Schur equation.}
Let
\[
  \Omega_p=\Omega_p^{(0)}\cap\Omega_p^{(1)}.
\]
By \eqref{e:p64omegaprob}, \(\mathbb P(\Omega_p)\to1\). On \(\Omega_p\), we have
\(\widehat\lambda_p\in I\) and
\(I\cap\operatorname{Spec}(C_p)=\varnothing\). The Schur complement formula in Lemma~\ref{a:schur} applied to \eqref{e:p64block} gives the exact equation
\begin{equation}
  \widehat\lambda_p
  =\lambda_pq_p+B_p(\widehat\lambda_p)
  \quad\text{on }\Omega_p.
  \label{e:p64S}
\end{equation}
Define a random remainder on the whole probability space by
\[
  r_p=
  \begin{cases}
  B_p(\widehat\lambda_p)-\lambda_ph_{c_p}(\widehat\lambda_p),
  &\text{on }\Omega_p,\\
  0,&\text{on }\Omega_p^c.
  \end{cases}
\]
Then \(|r_p|\le\Delta_p\) on \(\Omega_p\), and therefore
\begin{equation}
  r_p=o_{\mathbb P}(n^{-1/2}).
  \label{e:p64rem}
\end{equation}
Subtracting \eqref{e:rhoeq} from \eqref{e:p64S} gives, on \(\Omega_p\),
\begin{equation}
  \widehat\lambda_p-\rho_p
  =\lambda_p(q_p-1)+r_p
  +\lambda_p\{h_{c_p}(\widehat\lambda_p)-h_{c_p}(\rho_p)\}.
  \label{e:p64main}
\end{equation}

Lemma~\ref{l:resb} gives
\[
  \sup_{t\in I}|h_{c_p}'(t)|\le Cc_p.
\]
Since \(\widehat\lambda_p,\rho_p\in I\) on \(\Omega_p\), the mean value theorem gives a random number \(\xi_p\) such that
\[
  h_{c_p}(\widehat\lambda_p)-h_{c_p}(\rho_p)
  =\xi_p(\widehat\lambda_p-\rho_p),
  \qquad
  |\xi_p|\le Cc_p
\]
on \(\Omega_p\). Set \(\xi_p=0\) on \(\Omega_p^c\). Then \eqref{e:p64main} becomes
\begin{equation}
  (1-\lambda_p\xi_p)(\widehat\lambda_p-\rho_p)
  =\lambda_p(q_p-1)+r_p
  \quad\text{on }\Omega_p.
  \label{e:p64absorb}
\end{equation}
Because \(\lambda_p\) is bounded and \(c_p\to0\), for all large \(p\),
\[
  |1-\lambda_p\xi_p|\ge\frac12
  \quad\text{on }\Omega_p,
\]
and \((1-\lambda_p\xi_p)^{-1}=1+o(1)\) uniformly on \(\Omega_p\). Also,
\[
  \sqrt n(q_p-1)=O_{\mathbb P}(1),
  \qquad
  \sqrt n\,r_p=o_{\mathbb P}(1).
\]
It follows from \eqref{e:p64absorb} that
\[
  \mathbf1_{\Omega_p}\left[
    \sqrt n(\widehat\lambda_p-\rho_p)
    -\lambda_p\sqrt n(q_p-1)
  \right]
  =o_{\mathbb P}(1).
\]
Since \(\mathbb P(\Omega_p^c)\to0\), the indicator can be removed. Hence
\begin{equation}
  \sqrt n(\widehat\lambda_p-\rho_p)
  =\lambda_p\sqrt n(q_p-1)+o_{\mathbb P}(1).
  \label{e:p64reduction}
\end{equation}

\medskip
\noindent\textbf{Step 4: Apply the one-dimensional central limit theorem.}
In the real case, \(|z_k|^2=z_k^2\) and
\(z_k\sim \mathcal N (0,1)\). In the complex case,
\(|z_k|^2\sim\operatorname{Exp}(1)\). Thus
\[
  \sqrt n(q_p-1)
  =\frac1{\sqrt n}\sum_{k=1}^n(|z_k|^2-1)
  \xrightarrow{\ \mathrm{law}\ } \mathcal N (0,\kappa).
\]
Since \(\lambda_p\to\lambda_*\), Slutsky's theorem and
\eqref{e:p64reduction} give
\begin{equation}
  \sqrt n(\widehat\lambda_p-\rho_p)
  \xrightarrow{\ \mathrm{law}\ } \mathcal N (0,\kappa\lambda_*^2).
  \label{e:p64rhoCLT}
\end{equation}
Finally,
\[
  \rho_p
  =\lambda_p+\frac{p-K}{n}
  \frac{\sigma^2\lambda_p}{\lambda_p-\sigma^2}.
\]
The ratio on the right is bounded by \eqref{e:p64gap}--\eqref{e:A1}. Therefore
\[
  \sqrt n\left|
    \rho_p-\lambda_p
    -\frac pn\frac{\sigma^2\lambda_p}{\lambda_p-\sigma^2}
  \right|
  \le \frac{C K}{\sqrt n} \to0.
\]
Together with \eqref{e:p64rhoCLT}, this proves \eqref{e:dclt}.

It remains to prove \eqref{e:p64key}.

\medskip
\noindent\textbf{Step 5: Remove the conditional quadratic-form fluctuation.}
On \(\Omega_p^{(0)}\), put
\[
  A_p(t)=Z_0^*D_p^{1/2}R_p(t)D_p^{1/2}Z_0.
\]
Then
\begin{equation}
  B_p(t)=\frac{\lambda_p}{n^2}zA_p(t)z^*.
  \label{e:p64QF}
\end{equation}
Conditionally on \(Z_0\), the matrix \(A_p(t)\) is deterministic and Hermitian, while \(z\) is an independent standard Gaussian row vector. By Lemma~\ref{a:qform},
\begin{equation}
\begin{split}
  G_p(t)
  &:=\mathbb E_z[B_p(t)\mid Z_0]\\
  &=\lambda_p\frac1n
  \operatorname{Tr}\{C_p(tI-C_p)^{-1}\}.
\end{split}
  \label{e:p64mean}
\end{equation}
Here \(\mathbb E_z\) means expectation with respect to \(z\), with \(Z_0\) fixed.

On \(\Omega_p^{(0)}\), \eqref{e:p64resbd}--\eqref{e:p64Cpbd} give
\[
  \|D_p^{1/2}Z_0\|_{\mathrm{op}}^2
  =n\|C_p\|_{\mathrm{op}}\le Cn.
\]
Hence, uniformly for \(t\in I\),
\[
  \|A_p(t)\|_{\mathrm{op}}\le Cn,
  \qquad
  \operatorname{rank}A_p(t)\le p-1,
\]
and therefore
\[
  \operatorname{Tr}(A_p(t)^2)\le Cpn^2.
\]
Using the variance formula in Lemma~\ref{a:qform}, we get
\[
  \mathbf1_{\Omega_p^{(0)}}
  \operatorname{Var}_z(B_p(t)\mid Z_0)
  \le C\frac p{n^2}.
\]
Let \(F_p(t)=B_p(t)-G_p(t)\). Since \(\Omega_p^{(0)}\) depends only on \(Z_0\), conditional Cauchy--Schwarz and the tower property give, for each fixed \(t\in I\),
\begin{equation}
  \mathbb E\left[
    \mathbf1_{\Omega_p^{(0)}}|F_p(t)|
  \right]
  \le C\frac{\sqrt p}{n}
  =o(n^{-1/2}).
  \label{e:p64Qpoint}
\end{equation}

We now make the estimate uniform in \(t\). On \(\Omega_p^{(0)}\),
\[
  R_p'(t)=-R_p(t)^2,
\]
so
\[
  A_p'(t)
  =-Z_0^*D_p^{1/2}R_p(t)^2D_p^{1/2}Z_0,
  \qquad
  B_p'(t)=\frac{\lambda_p}{n^2}zA_p'(t)z^*.
\]
We do not need to exchange differentiation and conditional expectation. We can compute both sides directly from \eqref{e:p64mean}:
\[
  G_p'(t)
  =-\lambda_p\frac1n
  \operatorname{Tr}\{C_pR_p(t)^2\}
  =\mathbb E_z[B_p'(t)\mid Z_0].
\]
Thus
\[
  F_p'(t)=B_p'(t)-\mathbb E_z[B_p'(t)\mid Z_0].
\]
The same bounds as above, now with \(R_p(t)^2\), give
\begin{equation}
  \sup_{t\in I}
  \mathbb E\left[
    \mathbf1_{\Omega_p^{(0)}}|F_p'(t)|
  \right]
  \le C\frac{\sqrt p}{n}.
  \label{e:p64Qder}
\end{equation}
Fix \(t_0\in I\). The fundamental theorem of calculus gives, on \(\Omega_p^{(0)}\),
\[
  \sup_{t\in I}|F_p(t)|
  \le |F_p(t_0)|+\int_I|F_p'(u)|\,du.
\]
Using \eqref{e:p64Qpoint}, \eqref{e:p64Qder}, Fubini's theorem, and Lemma~\ref{l:localization}, we obtain
\begin{equation}
  \mathbf1_{\Omega_p^{(0)}}
  \sup_{t\in I}|B_p(t)-G_p(t)|
  =o_{\mathbb P}(n^{-1/2}).
  \label{e:p64Q}
\end{equation}

\medskip
\noindent\textbf{Step 6: Remove the other \(K-1\) spike coordinates.}
Write
\[
  Z_0=\begin{pmatrix}Z_s\\ Z_b\end{pmatrix},
  \qquad
  Z_s\in\mathbb F^{(K-1)\times n},
  \qquad
  Z_b\in\mathbb F^{m\times n},
\]
and embed \(C_b\) into a \((p-1)\times(p-1)\) matrix:
\[
  \widetilde C_b
  =\operatorname{diag}(0_{K-1},C_b).
\]
The matrix
\[
  E_p=C_p-\widetilde C_b
\]
has nonzero entries only in the first \(K-1\) rows or the first \(K-1\) columns. Hence
\[
  \operatorname{rank}(E_p)\le2(K-1).
\]
Also, \(\|E_p\|_{\mathrm{op}}\le C\) on \(\Omega_p^{(0)}\). On the same event, the resolvents of both \(C_p\) and \(\widetilde C_b\) are uniformly bounded for \(t\in I\). Therefore Lemma~\ref{a:frtrace} gives
\begin{equation}
  \mathbf1_{\Omega_p^{(0)}}
  \sup_{t\in I}
  \left|
    \frac1n\operatorname{Tr}\{(tI-C_p)^{-1}
    -(tI-\widetilde C_b)^{-1}\}
  \right|
  \le\frac Cn.
  \label{e:p64ResTr}
\end{equation}
Using
\[
  C(tI-C)^{-1}=t(tI-C)^{-1}-I
\]
for the two matrices of the same size, the identity terms cancel. Since the zero block of \(\widetilde C_b\) contributes nothing, \eqref{e:p64ResTr} gives
\begin{equation}
\begin{split}
  &\mathbf1_{\Omega_p^{(0)}}
  \sup_{t\in I}
  \left|
    \frac1n\operatorname{Tr}\{C_p(tI-C_p)^{-1}\}
    -\frac1n\operatorname{Tr}\{C_b(tI_m-C_b)^{-1}\}
  \right|\\
  &\hspace{40mm}\le\frac Cn.
\end{split}
  \label{e:p64finite}
\end{equation}

Now apply Lemma~\ref{l:resb} to \(C_b\). We have \(m/n=c_p\to0\) and \(I\subset(\sigma^2+\varepsilon,\infty)\), so
\[
  \sup_{t\in I}
  \left|
    \frac1n\operatorname{Tr}\{C_b(tI_m-C_b)^{-1}\}
    -h_{c_p}(t)
  \right|
  =o_{\mathbb P}(n^{-1/2}).
\]
Since \(n^{-1}=o(n^{-1/2})\), \eqref{e:p64finite} implies
\begin{equation}
  \mathbf1_{\Omega_p^{(0)}}
  \sup_{t\in I}
  \left|
    \frac1n\operatorname{Tr}\{C_p(tI-C_p)^{-1}\}
    -h_{c_p}(t)
  \right|
  =o_{\mathbb P}(n^{-1/2}).
  \label{e:p64T}
\end{equation}
Combining \eqref{e:p64mean}, \eqref{e:p64Q}, and \eqref{e:p64T}, and using the boundedness of \(\lambda_p\), proves \eqref{e:p64key}. This completes the proof.

\end{proof}
\section*{Acknowledgments}
This work was supported in part by JSPS KAKENHI Grant Numbers 25H00593, 23K03133 and 26K06828, and by a JSPS Bilateral Joint Research Project 120249936. 

ChatGPT (OpenAI) was used to assist with language editing, literature searches, the examination of mathematical arguments during the preparation of this manuscript, and the \LaTeX{} formatting and preparation of the submission files. The authors take full responsibility for the content of the paper.
\bibliographystyle{plain}   
\bibliography{bibliography}
\end{document}